\documentclass{article} 
\usepackage{PRIMEarxiv}

\usepackage{amsmath}
\usepackage{amsthm}
\usepackage{amsopn}
\usepackage{amsfonts}
\usepackage{amssymb}
\usepackage{mathrsfs}
\usepackage{bbm}
\usepackage{enumitem}
\usepackage{hyperref}
\usepackage{url}
\usepackage{graphicx}
\usepackage{subcaption}
\usepackage{algorithm}
\usepackage{algorithmic}
\usepackage{natbib}
\usepackage{booktabs}
\usepackage{cleveref}
\usepackage{todonotes}
\usepackage{pifont}
\usepackage{dsfont}

\usepackage{localmacros}

\newtheorem{theorem}{Theorem}
\newtheorem{corollary}{Corollary}[theorem]
\newtheorem{lemma}{Lemma}
\newtheorem{proposition}{Proposition}
\newtheorem{definition}{Definition}
\newtheorem{remark}{Remark}

\title{Iso-Riemannian Optimization on \\
Learned Data Manifolds}
\author{Willem Diepeveen \\
Department of Mathematics\\
University of California, Los Angeles\\
Los Angeles, CA 90095, USA \\
\texttt{wdiepeveen@math.ucla.edu} \\
\And
Melanie Weber \\
School of Engineering and Applied Sciences\\
Harvard University \\
Cambridge, MA 02138, USA \\
\texttt{mweber@seas.harvard.edu}
}

\allowdisplaybreaks

\begin{document}

\maketitle

\begin{abstract}
    We develop a theory of iso-Riemannian optimization for problems constrained to learned data manifolds, a setting in which classical Riemannian optimization -- and Riemannian gradient descent in particular -- can be poorly suited. That is, favorable Euclidean properties of an objective need not translate into geodesic convexity or \(L\)-smoothness, and the Riemannian gradient may provide an unsuitable search direction. We instead combine the manifold mappings induced by the iso-connection, whose geodesics have constant Euclidean speed, with the \(\ell^2\)-projected Euclidean gradient, resulting in \(\ell^2\)-projected gradient iso-Riemannian descent, in an attempt to alleviate these issues. We analyze this scheme from two complementary perspectives. First, we introduce iso-g-convexity and iso-\(L\)-smoothness, relate strong iso-g-convexity to a Polyak-\L{}ojasiewicz-type condition, and establish general convergence guarantees. This function-based theory addresses the conditioning and search-direction issues that motivate our approach. Second, we introduce iso-monotonicity and iso-Lipschitzness for vector fields. However, while these notions yield analogous convergence results in one dimension and enable applications such as the computation of iso-Riemannian barycentres, the resulting theory need not extend meaningfully to higher dimensions. So unlike in the classical Levi-Civita setting, these function- and vector field-based perspectives need not be equivalent in the iso-Riemannian setting. Consequently, distinct extensions of classical Riemannian optimization lead to different assumptions and convergence guarantees. The theory developed here suggests that, for function optimization, the function-based perspective provides the more general, tractable and overall suitable framework, whereas the vector field perspective is naturally reserved for problems without a direct function-based formulation.

\end{abstract}

\blfootnote{Our code is available at \href{https://github.com/Weber-GeoML/Iso-Riem-Opt}{https://github.com/Weber-GeoML/Iso-Riem-Opt}.}

\section{Introduction}
\label{sec:introduction}

Many real-valued datasets are widely believed to concentrate near a low-dimensional geometric structure, commonly referred to as the \emph{manifold hypothesis} \cite{fefferman2016testing,horvat2022intrinsic,kamkari2024geometric,loaiza2024deep,pope2021intrinsic,stanczuk2022your}. If this structure is modeled as a smooth manifold \(\manifold \subset \Real^\dimInd\), it is natural to formulate optimization problems -- such as those arising in inverse problems -- subject to the constraint that the solution lies on \(\manifold\). In particular, one seeks to solve
\begin{equation}
    \inf_{\Vector \in \manifold} \function(\Vector),
    \label{eq:mfld-opt-problem-init}
\end{equation}
where \(\function : \Real^\dimInd \to \Real\) is the objective function. Despite the apparent simplicity of \eqref{eq:mfld-opt-problem-init}, a central practical challenge is that the manifold \(\manifold\) is typically not available in closed form for generic data modalities, such as natural images.

Only recently has it become feasible not only to learn the manifold \(\manifold\), but also to equip it with a Riemannian structure admitting closed-form manifold mappings. In particular, recent work has shown that generative modeling can efficiently learn pullback manifold structures on the entire ambient space $\mathbb{R}^{\dimInd}$, together with manifold mappings that are fast to evaluate \cite{diepeveen2025scorebased}. This construction can then be leveraged to extract a lower-dimensional submanifold $\manifold$ approximating the dataset, which is totally geodesic; consequently, the Riemannian metric on $\manifold$ inherited from the ambient pullback structure induces the same manifold mappings, with fast evaluation retained by construction\footnote{We note that there is already a substantial body of work on learning Riemannian structure from data \cite{arvanitidis2016locally,diepeveen2024pulling,hauberg2012geometric,peltonen2004improved,Scarvelis2023,sorrenson2025learning,sun2024geometryaware}; see \cite{gruffaz2025riemannian} for an overview. However, to the best of our knowledge, \cite{diepeveen2025scorebased} is the only approach that is both efficiently trainable and permits efficient evaluation of the learned geometric structure, namely the manifold mappings.}.

While one can now formulate \eqref{eq:mfld-opt-problem-init} either as a constrained Euclidean optimization problem or as an unconstrained Riemannian one, the difficulty lies in the fact that standard methods in both settings face fundamental challenges. From the Euclidean perspective, the constraint set is known only implicitly, so standard techniques such as projection, Lagrange multipliers, and penalty methods are not directly applicable. From the Riemannian perspective, even when \(\function\) has favorable Euclidean properties -- such as convexity in the classical inverse problem \(\function(\Vector) := \|\forward \Vector - \mathbf{b}\|_2^2\), with \(\forward \in \Real^{\dimIndB \times \dimInd}\) and \(\mathbf{b} \in \Real^{\dimIndB}\) -- these properties do not generally imply geodesic convexity, which underlies the theory of Riemannian optimization. 

\begin{figure}[p]
    \centering
    \captionsetup[subfigure]{justification=centering}

    \begin{subfigure}{\linewidth}
        \centering
        \includegraphics[width=0.35\linewidth]{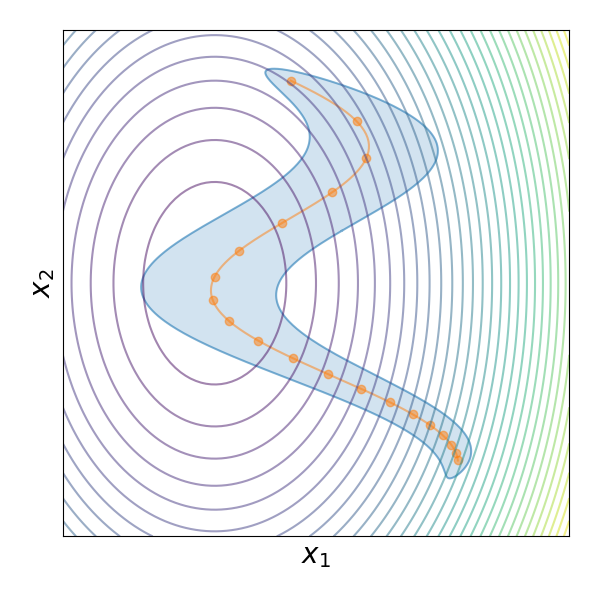}
        \hspace{0.03\linewidth}
        \includegraphics[width=0.35\linewidth]{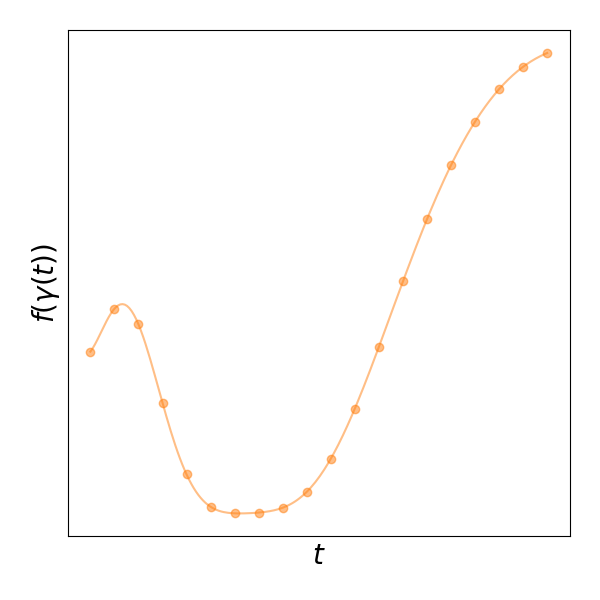}
        \caption{Geodesics on manifolds generally exhibit non-constant
        Euclidean speed, as illustrated by the orange curve on the blue
        manifold (left), leading to poor apparent conditioning of the objective
        function along this geodesic (right).}
        \label{fig:geodesic-convexity-issue}
    \end{subfigure}

    \vspace{0.6em}

    \begin{subfigure}{\linewidth}
        \centering
        \includegraphics[width=0.35\linewidth]{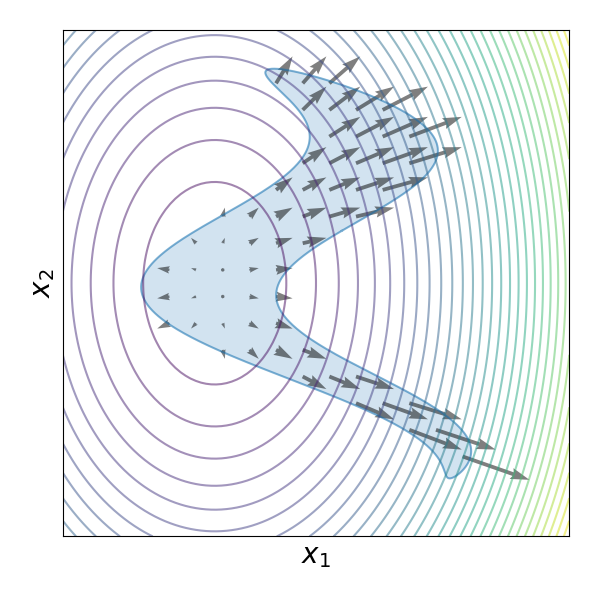}
        \hspace{0.03\linewidth}
        \includegraphics[width=0.35\linewidth]{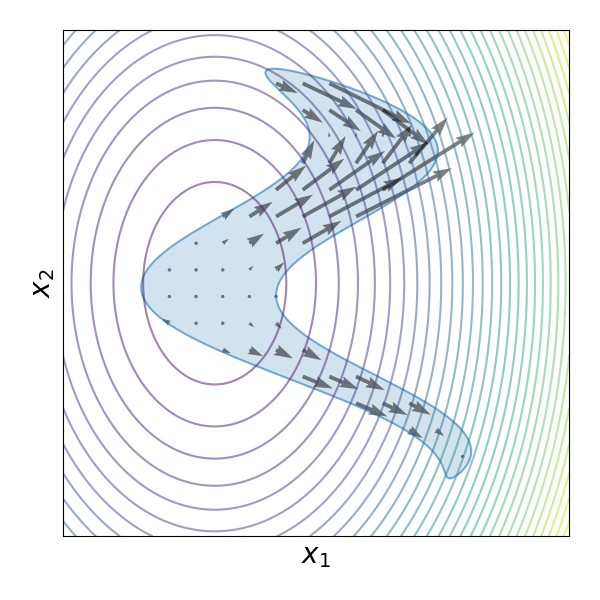}
        \caption{The Euclidean gradient (left) and Riemannian gradient (right) can lead to
        different descent behavior. Here, descending along the negative
        Euclidean gradient moves toward the minimizer, while Riemannian-gradient
        descent may instead approach a poor local minimum on the manifold.}
        \label{fig:eucliden-vs-riemannian-gradients-issue}
    \end{subfigure}

    \vspace{0.6em}

    \begin{subfigure}{\linewidth}
        \centering
        \includegraphics[width=0.35\linewidth]{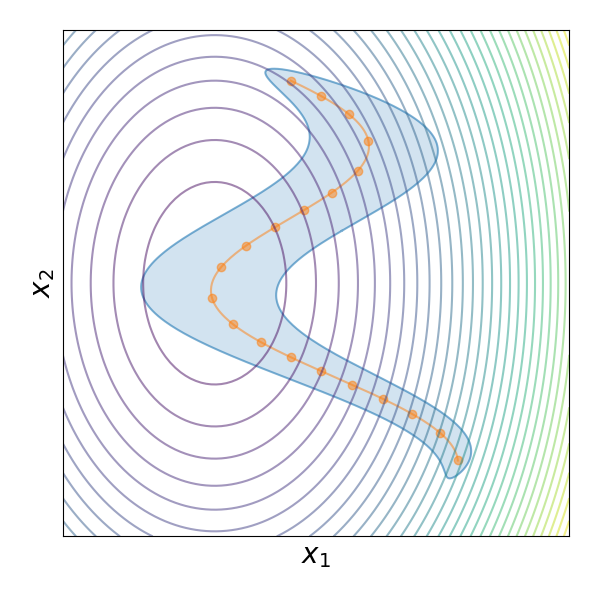}
        \hspace{0.03\linewidth}
        \includegraphics[width=0.35\linewidth]{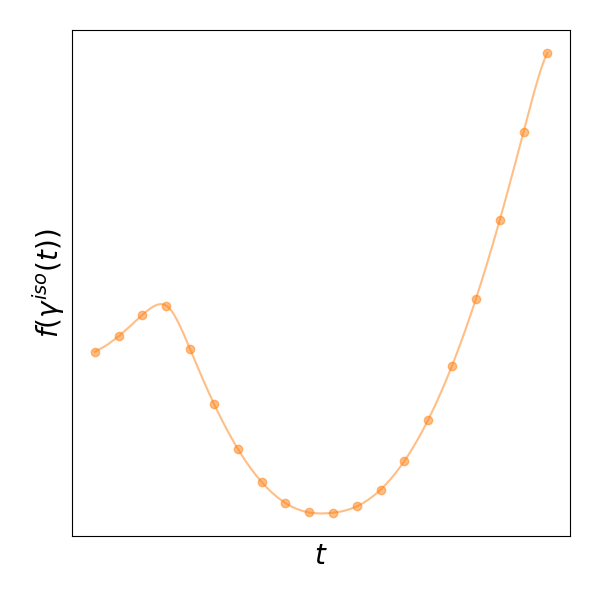}
        \caption{Iso-geodesics exhibit constant Euclidean speed, as illustrated by the orange curve on the blue manifold (left), leading to
        improved apparent conditioning of the objective function along the
        iso-geodesic (right).}
        \label{fig:iso-geodesic-convexity-issue-sol}
    \end{subfigure}

    \caption{
Limitations of standard Riemannian geometry for optimization in (a) and (b), and a first motivation for iso-Riemannian geometry as a path
forward in (c).
}
    \label{fig:river-cross}
\end{figure}

Since Euclidean methods provide limited leverage if even basic algorithms cannot be evaluated reliably, we revisit Riemannian optimization in the learned-manifold setting, where two principal difficulties arise, as illustrated in \Cref{fig:geodesic-convexity-issue,fig:eucliden-vs-riemannian-gradients-issue}. \Cref{fig:geodesic-convexity-issue} shows that geodesics under the learned pullback structure typically have non-constant Euclidean speed: although their velocity has constant norm in the learned Riemannian metric, they may accelerate or decelerate from a Euclidean perspective. Consequently, the intrinsic conditioning of the optimization problem may be effectively exacerbated, leading to step sizes that are smaller than necessary. \Cref{fig:eucliden-vs-riemannian-gradients-issue} illustrates that the Riemannian gradient may be a less suitable search direction than the Euclidean gradient\footnote{Or the projected Euclidean gradient, if we interpret this as a two-dimensional slice of a higher-dimensional problem.}, which is the more natural descent direction for the Euclidean convex function shown there. In particular, the Riemannian gradient may even steer the iterates toward an undesired local minimizer. These two effects typically occur simultaneously in practice. Thus, although the stationary points of \(\function\) are independent of the particular Riemannian metric chosen on \(\manifold\), whether these points can be reached efficiently depends crucially on how Riemannian optimization is formulated and implemented in this learned geometric setting.

From a theoretical standpoint, changing the metric can partially alleviate these issues, but only at a cost. If $\manifold$ is endowed with the ambient $\ell^2$ metric rather than the learned one, geodesics have constant $\ell^2$-speed and the Riemannian gradient becomes the projected Euclidean gradient. However, $\manifold$ is no longer necessarily a totally geodesic submanifold, and closed-form access to the manifold mappings is lost. The former is problematic because, on bounded submanifolds such as in \Cref{fig:geodesic-convexity-issue}, the geodesic connecting two endpoints may exit the manifold in order to cut across the interior bends as directly as possible. The latter is problematic because it compromises efficiency. Taken together, these drawbacks make this a poor foundation for both theory and algorithms.

Since optimizing on $\manifold$ under either the learned metric or the ambient $\ell^2$ metric runs into both theoretical and practical issues, we are led to depart from standard Riemannian optimization. At a high level, this paper seeks a middle ground by incorporating Euclidean geometry into the notions of geodesic convexity and Riemannian gradient descent, thereby developing a new framework for optimization on learned data manifolds. In particular, reparametrized geodesics with constant Euclidean speed may help alleviating the first issue and enable a new form of convergence analysis, in which geodesic convexity, $L$-smoothness, and more generally conditioning are considered after reparametrization, as shown in \Cref{fig:iso-geodesic-convexity-issue-sol}\footnote{Although we do not expect this to eliminate all non-convexity induced by the manifold constraint, as also shown in \Cref{fig:iso-geodesic-convexity-issue-sol}.}. For the second issue, our discussion above suggests that using the projected Euclidean gradient as the search direction may be more effective, and we aim to understand this choice through the same reparametrized-geodesic perspective. Since reparametrized geodesics give rise to iso-geodesics, namely geodesics under the recently proposed iso-connection \cite{diepeveen2025manifold}, this naturally raises the question of whether the resulting iso-Riemannian geometry can provide a suitable framework for generalizing Riemannian optimization and for replacing the Riemannian gradient with the projected Euclidean gradient. We aim to answer this by developing a theory of iso-Riemannian optimization and deriving conditions guaranteeing convergence of the following first-order method:
\begin{equation}
    \Vector^{(\sumIndC+1)} := \exp^{\iso}_{\Vector^{(\sumIndC)}} \bigl( - r \mathbf{P}_{\Vector^{(\sumIndC)}} \nabla \function (\Vector^{(\sumIndC)}) \bigr), \quad r>0,\;\Vector^{(0)}\in \manifold,
    \label{eq:iso-first-order-scheme}
\end{equation}
where \(\exp^{\iso}_{\Vector} : \tangent_{\Vector} \Real^\dimInd \to \Real^\dimInd\) is the iso-exponential map \cite{diepeveen2025manifold} generated by \((\manifold, (\cdot,\cdot))\), and \(\mathbf{P}_{\Vector} : \tangent_\Vector \Real^\dimInd \to \tangent_\Vector \manifold\) is the \(\ell^2\)-projection onto the tangent space at \(\Vector \in \manifold\), which we will refer to throughout the paper as $\ell^2$-projected gradient iso-Riemannian descent ($\ell^2$-PG-IRD).

\subsection{Contributions}

To the best of our knowledge, this work is the first to study first-order optimization in an iso-Riemannian setting or, more generally, with a non-Levi-Civita connection. Our approach departs from both the classical literature on optimization over abstract manifolds, which, following early projection-based work \cite{luenberger1972gradient}, is predominantly formulated with respect to the Levi-Civita connection \cite{ablin2022fast,ablin2024infeasible,absil2007trust,adachi2022riemannian,agarwal2021adaptive,bacak2014computing,bergmann2016parallel,bergmann2021fenchel,bergmann2024difference,bergmann2024riemannian,diepeveen2021inexact,ferreira1998subgradient,ferreira2002proximal,goyens2025geometric,hoseini2023proximal,huang2015broyden,huang2018riemannian,lai2024riemannian,liu2020simple,smith1994optimization,souza2015proximal,udriste1994convex,vary2024optimization,weber2021projection,weber2023riemannian}; see \cite{absil2008optimization,boumal2023introduction} for overviews, and existing approaches to learned data manifolds, which typically optimize over linear subspaces or through global or local learned charts \cite{cartis2024learning,cosson2023low,daras2021intermediate,hand2018global,hand2018phase,lei2019inverting,robinett2025manifold,shustin2022manifold,sober2020manifold,sritharan2021computing}. See \Cref{sec:related-work} for a detailed discussion.

A central theme of this work is that several formulations of classical Riemannian optimization that are equivalent in the Levi-Civita setting can lead to distinct assumptions, algorithmic interpretations, and convergence guarantees in the iso-Riemannian setting. We develop two such perspectives: a function-based approach, formulated through iso-\(g\)-convexity and iso-\(L\)-smoothness, and a vector-field-based approach, formulated through iso-monotonicity and iso-Lipschitzness.


In classical Riemannian optimization, convergence of gradient-based methods can be established from either perspective. From a function-based viewpoint, strong geodesic convexity -- or, more generally, a Polyak-\L{}ojasiewicz (PL) condition -- together with \(L\)-smoothness of the objective yields convergence. From a vector field viewpoint, convergence follows from strong monotonicity and Lipschitz continuity of the descent field. In the Levi-Civita setting, these formulations are closely related: for example, strongly geodesically convex functions have strongly monotone Riemannian gradients, and \(L\)-smooth functions have Lipschitz Riemannian gradients. In contrast, these implications need not persist in the iso-Riemannian setting.

For clarity, and in view of the learned data manifolds of primary interest, we focus on unbounded manifolds. This allows us to isolate the principal geometric and optimization-theoretic issues without the additional technicalities associated with boundaries. The phenomena studied here are not restricted to unbounded settings, and we expect analogous considerations to arise on bounded manifolds; a systematic treatment of the latter, however, lies beyond the scope of this work.

Our contributions are as follows.

\paragraph{Convergence from a function-based perspective}
We extend geodesic convexity and smoothness to the iso-Riemannian setting by introducing iso-g-convexity and iso-\(L\)-smoothness. These notions are obtained by replacing manifold mappings with iso-manifold mappings, Riemannian inner products with \(\ell^2\)-inner products, and Riemannian gradients with projected Euclidean gradients. 

We establish that they recover standard Euclidean convexity and \(L\)-smoothness on \((\mathbb{R},(\cdot,\cdot))\), relate them to their Euclidean counterparts through a geometric correction on Riemannian manifolds, and show that strong iso-g-convexity is closely connected to a PL-type condition. This yields a general convergence result for \(\ell^2\)-PG-IRD and resolves the conditioning and search-direction issues illustrated in \Cref{fig:geodesic-convexity-issue,fig:iso-geodesic-convexity-issue-sol,fig:eucliden-vs-riemannian-gradients-issue}. 

We further consider least-squares problems arising in inverse problems, for which the proposed methodology can improve on the unconstrained formulation, and demonstrate the practical utility of the scheme for both modeled and learned pullback structures.

\paragraph{Convergence from a vector field perspective}
We introduce iso-monotonicity and iso-Lipschitzness for vector fields and investigate the convergence of iterations evolving along such fields. On \((\mathbb{R},(\cdot,\cdot))\) and on one-dimensional manifolds, these notions are equivalent to iso-g-convexity and iso-\(L\)-smoothness, respectively, and yield an analogous convergence theory. 

Beyond one dimension, however, we do not presently obtain comparably general results; instead, iso-monotonicity and iso-Lipschitzness may be overly restrictive in higher dimensions, suggesting that the function-based perspective provides a more suitable framework for studying optimization in the iso-Riemannian setting. 

Nonetheless, the vector field perspective applies naturally to problems posed as finding zeros of a vector field, such as the iso-Riemannian barycentre problem, for which no direct function-based formulation is available. Despite the currently limited theory, the resulting algorithms remain effective in practice for modeled and learned pullback structures and support downstream tasks such as iso-\(K\)-means clustering, which we also introduce.




\subsection{Outline}

In \Cref{sec:prelims}, we introduce the geometric preliminaries used in the paper. In particular, we review Euclidean pullback manifolds, outline the learning of Riemannian submanifolds from data -- thereby motivating the setting considered in this work -- and present the main ingredients of iso-Riemannian geometry. Next, in \Cref{sec:function-convergence}, we develop a function-based convergence theory based on iso-g-convexity and iso-\(L\)-smoothness, and establish convergence guarantees for the first-order scheme in \cref{eq:iso-first-order-scheme}. In \Cref{sec:vectorfield-convergence}, we study the complementary vector field perspective through iso-monotonicity and iso-Lipschitzness, discuss its more limited higher-dimensional theory, and apply it to the computation of iso-Riemannian barycentres. Finally, \Cref{sec:conclusions} summarizes our findings and discusses directions for future work, including extensions of the proposed framework to retraction-based methods. Although such extensions lie beyond the scope of the present work, we expect the framework to extend naturally beyond the iso-exponential map setting considered here.
\section{Preliminaries}
\label{sec:prelims}
We first review some standard notation from differential and Riemannian geometry, referring the reader to \cite{boothby2003introduction,carmo1992riemannian,lee2013smooth,sakai1996riemannian} for further background, and set the stage for Euclidean pullback manifolds on $\Real^\dimInd$, following \cite{diepeveen2024pulling}. Next, we outline, at a high level, the procedure for learning a Riemannian submanifold $(\manifold, (\cdot, \cdot))$ of $\mathbb{R}^d$ that approximates a given data set, as developed in \cite{diepeveen2025scorebased,diepeveen2025manifold}. Finally, we present the key ideas underlying iso-Riemannian geometry on such $(\manifold, (\cdot, \cdot))$, with additional details provided in \cite{diepeveen2025manifold}.

\subsection{Riemannian geometry}

\paragraph{Smooth manifolds and tangent spaces} 
Let $\manifold$ be a \emph{$\dimInd$-dimensional smooth manifold}, i.e., a topological manifold of dimension $\dimInd$ equipped with a \emph{maximal smooth atlas}, meaning a collection of charts whose transition functions are all smooth, making the manifold locally diffeomorphic to $\Real^\dimInd$. We write $C^\infty(\manifold)$ for the space of smooth functions over $\manifold$. The \emph{tangent space} at $\mPoint \in \manifold$, which is defined as the space of all \emph{derivations} at $\mPoint$, is denoted by $\tangent_\mPoint \manifold$ and for \emph{tangent vectors} we write $\tangentVector_\mPoint \in \tangent_\mPoint \manifold$. For the \emph{tangent bundle} we write $\tangent\manifold$ and smooth vector fields, which are defined as \emph{smooth sections} of the tangent bundle, are written as $\vectorfield(\manifold) \subset \tangent\manifold$.

\paragraph{Riemannian manifolds} 
A smooth manifold $\manifold$ becomes a \emph{Riemannian manifold} if it is equipped with a smoothly varying \emph{metric tensor field} $(\cdot, \cdot) : \vectorfield(\manifold) \times \vectorfield(\manifold) \to C^\infty(\manifold)$. This tensor field induces a \emph{(Riemannian) distance} $\distance_{\manifold} : \manifold\times\manifold\to\Real$. The metric tensor can also be used to construct a unique affine connection, the \emph{Levi-Civita connection}, that is denoted by $\nabla_{(\,\cdot\,)}(\,\cdot\,) : \vectorfield(\manifold) \times \vectorfield(\manifold) \to \vectorfield(\manifold)$. 
This connection is in turn the cornerstone of a myriad of manifold mappings.

One is the notion of a \emph{geodesic}, which for two points $\mPoint,\mPointB \in \manifold$ is defined as a curve $\geodesic_{\mPoint,\mPointB} : [0,1] \to \manifold$ with minimal length that connects $\mPoint$ with $\mPointB$ -- that is, if such a curve exists. 
To ensure existence, we often consider \emph{(geodesically) convex} subsets, i.e., sets $\mathcal{D}\subset \manifold$ such that $\geodesic_{\mPoint,\mPointB}\subset \mathcal{D}$ for any $\mPoint,\mPointB\in \mathcal{D}$. In addition, when geodesics are also unique on $\mathcal{D}$, we call $\mathcal{D}$ \emph{strongly (geodesically) convex}. 
Another closely related notion to geodesics is the curve $t \mapsto \geodesic_{\mPoint,\tangentVector_\mPoint}(t)$  for a geodesic starting from $\mPoint\in\manifold$ with velocity $\dot{\geodesic}_{\mPoint,\tangentVector_\mPoint} (0) = \tangentVector_\mPoint \in \tangent_\mPoint\manifold$. This can be used to define the \emph{exponential map} $\exp_\mPoint : \mathcal{D}_\mPoint \to \manifold$ at $\mPoint$ as \(\exp_\mPoint(\tangentVector_\mPoint) := \geodesic_{\mPoint,\tangentVector_\mPoint}(1),\) where \(\mathcal{D}_\mPoint \subset \tangent_\mPoint \manifold\) is the set on which \(\geodesic_{\mPoint,\tangentVector_\mPoint}(1)\) is defined. The manifold $\manifold$ is said to be \emph{(geodesically) complete} whenever $\mathcal{D}_{\mPoint}=\tangent_{\mPoint} \manifold$ for all $\mPoint \in \manifold$. Furthermore, the \emph{logarithmic map} $\log_\mPoint : \exp_\mPoint(\mathcal{D}'_\mPoint ) \to \mathcal{D}'_\mPoint$ at $\mPoint$ is defined as the inverse of $\exp_\mPoint$, so it is well-defined on  $\mathcal{D}'_{\mPoint} \subset \mathcal{D}_{\mPoint}$ where $\exp_\mPoint$ is a diffeomorphism. Moreover, for \emph{parallel transport} $\partransport_{\mPointB \leftarrow \mPoint}: \tangent_\mPoint \manifold \to \tangent_\mPointB \manifold$ of a vector $\tangentVector_\mPoint \in \tangent_\mPoint \manifold$ along a geodesic from $\mPoint$ to
$\mPointB$ we write $\partransport_{\mPointB \leftarrow \mPoint} \tangentVector_\mPoint$. Finally, the \emph{Riemannian gradient} of a smooth function $\function \in C^\infty(\manifold)$ denotes the unique vector field $\Grad \function \in \vectorfield (\manifold)$ such that $(\Grad \function, \tangentVector_\mPoint)_{\mPoint} = \tangentVector_\mPoint \function := D_{\mPoint} \function[\tangentVector_\mPoint]$ holds for any $\tangentVector \in \vectorfield (\manifold)$ and $\mPoint \in \manifold$, where
$D_{\mPoint} \function: \tangent_{\mPoint} \manifold \to \Real$ denotes the differential of $\function$.

\paragraph{Pullback manifolds} 
If $(\manifold, (\cdot,\cdot))$ is a $\dimInd$-dimensional Riemannian manifold, $\manifoldB$ is a $\dimInd$-dimensional smooth manifold and $\diffeo:\manifoldB \to \manifold$ is a diffeomorphism, the \emph{pullback metric}
\begin{equation}
    (\tangentVector, \tangentVectorB)_\mPoint^\diffeo := (D_{\mPoint}\diffeo[\tangentVector_{\mPoint}], D_{\mPoint}\diffeo[\tangentVectorB_{\mPoint}])_{\diffeo(\mPoint)}, \quad \mPoint \in \manifoldB, \tangentVector, \tangentVectorB \in \vectorfield(\manifoldB)
    \label{eq:pull-back-metric}
\end{equation}
where $D_{\mPoint}\diffeo: \tangent_\mPoint \manifoldB \to \tangent_{\diffeo(\mPoint)} \manifold$ denotes the differential of $\diffeo$,
defines a Riemannian structure on $\manifoldB$, which we denote by $(\manifoldB, (\cdot,\cdot)^\diffeo)$. 
Pullback mappings are denoted similarly to (\ref{eq:pull-back-metric}) with the diffeomorphism $\diffeo$ as a superscript, i.e., we write $\distance^\diffeo_{\manifoldB}(\mPoint, \mPointB)$, $\geodesic^\diffeo_{\mPoint, \mPointB}$, $\exp^\diffeo_\mPoint (\tangentVector_\mPoint)$, $\log^\diffeo_{\mPoint} \mPointB$, and $\partransport^{\diffeo}_{\mPointB \leftarrow \mPoint}$ for $\mPoint,\mPointB \in \manifoldB$ and $\tangentVector_\mPoint \in \tangent_\mPoint \manifoldB$. Pullback metrics literally pull back all geometric information from the Riemannian structure on $\manifold$. 
In particular, closed-form manifold mappings on $(\manifold, (\cdot,\cdot))$ yield under mild assumptions closed-form manifold mappings on $(\manifoldB, (\cdot,\cdot)^\diffeo)$. 

\paragraph{Euclidean pullback manifolds}
For Euclidean pullback manifolds $(\Real^\dimInd,(\cdot,\cdot)^\diffeo)$ generated by a diffeomorphism $\diffeo:\Real^\dimInd\to \Real^\dimInd$ pulling back the standard Euclidean structure $(\Real^\dimInd, (\cdot, \cdot)_2)$, we have
\begin{align}
    \distance_{\Real^{\dimInd}}^{\diffeo}(\Vector, \VectorB) &= \|\diffeo(\Vector) - \diffeo(\VectorB)\|_2,
    \label{eq:thm-distance-remetrized}\\
    \geodesic^{\diffeo}_{\Vector, \VectorB}(t) &= \diffeo^{-1}((1 - t)\diffeo(\Vector) + t \diffeo(\VectorB)),
    \label{eq:thm-geodesic-remetrized}\\
    \exp^{\diffeo}_\Vector (\tangentVector_\Vector) &= \diffeo^{-1}(\diffeo(\Vector) + D_{\Vector} \diffeo[\tangentVector_\Vector]),
    \label{eq:thm-exp-remetrized}\\
    \log^{\diffeo}_{\Vector} (\VectorB) &= D_{\diffeo(\Vector)}\diffeo^{-1}[\diffeo(\VectorB) - \diffeo(\Vector)],
    \label{eq:thm-log-remetrized}\\
    \partransport^{\diffeo}_{\VectorB \leftarrow \Vector} \tangentVector_\Vector &= D_{\diffeo(\VectorB)}\diffeo^{-1}[D_{\Vector} \diffeo[\tangentVector_\Vector]],
    \label{eq:thm-exp-remetrized}
\end{align}
where $\Vector, \VectorB\in \Real^\dimInd$ and $\tangentVector_\Vector \in \tangent_\Vector \Real^\dimInd \cong \Real^\dimInd$, and have \cite[Prop~3.7]{diepeveen2024pulling}
\begin{equation}
    \operatorname{argmin}_{\Vector\in \Real^\dimInd} \sum_{\sumIndA=1}^\dataPointNum \distance^\diffeo_{\Real^\dimInd}(\Vector, \Vector^\sumIndA)^2 = \diffeo^{-1} (\frac{1}{\dataPointNum} \sum_{\sumIndA=1}^\dataPointNum \diffeo(\Vector^\sumIndA)),
    \label{eq:thm-bary-remetrized}
\end{equation}
for the Riemannian barycentre \cite{karcher1977riemannian}, where $\Vector^1, \ldots, \Vector^\dataPointNum\in \Real^\dimInd$. 

\subsection{Learned data manifolds}

\paragraph{Data-driven Euclidean pullback manifolds}
The quality of the Euclidean pullback structure induced by a diffeomorphism $\diffeo:\mathbb{R}^\dimInd\to\mathbb{R}^\dimInd$ for a data set $\Vector^1, \ldots, \Vector^n\sim \density_{\text{data}}$ can be measured by how well the associated geodesics traverse high-likelihood regions of the data distribution. Concretely, one seeks a $\diffeo$ such that the function
\begin{equation}
    t \mapsto - \log\bigl(\density_{\text{data}} (\geodesic^{\diffeo}_{\Vector, \VectorB}(t))\bigr)
\end{equation}
remains small for all $t\in $ and for all pairs $\Vector, \VectorB$ in the data set.

In practice, such a diffeomorphism is learned by minimizing the negative log-likelihood of a distribution $\density_{\networkParams}: \mathbb{R}^\dimInd \to \mathbb{R}$ of the form
\begin{equation}
\density_{\networkParams}(\Vector)\propto \exp \Bigl(-\frac{1}{2} \|\diffeo_\networkParams(\Vector)\|^2 \Bigr),
\label{eq:stat-model-geom}
\end{equation}
where $\diffeo_\networkParams:\mathbb{R}^\dimInd\to\mathbb{R}^\dimInd$ is a diffeomorphic neural network with parameters $\networkParams$.

To build intuition for why this objective yields a suitable diffeomorphism, observe first that for any choice of endpoints $\Vector,\VectorB\in \mathbb{R}^\dimInd$ and parameters $\networkParams$, the function
\begin{equation}
    t \mapsto -\log\bigl(\density_{\networkParams}(\geodesic^{\diffeo_{\networkParams}}_{\Vector, \VectorB}(t))\bigr)
\end{equation}
is convex. This can be interpreted as the geodesics preferentially passing through high-density regions of $\density_{\networkParams}$, which aligns with our goal provided that $\density_{\networkParams} \approx \density_{\text{data}}$. If the data distribution is feasible in the model class -- that is, if there exists $\networkParams^*$ such that $\density_{\networkParams^*} = \density_{\text{data}}$ -- then minimizing the negative log-likelihood recovers these parameters $\networkParams^*$. Equivalently, minimizing the negative log-likelihood implicitly encourages geodesics between data points to concentrate in high-density regions of $\density_{\text{data}}$, without requiring any explicit path-based or contrastive loss.




\paragraph{Extracting data manifolds from data-driven Euclidean pullback manifolds}
Using the manifold mappings on $(\mathbb{R}^\dimInd, (\cdot,\cdot)^\diffeo)$ for $\diffeo := \diffeo_{\networkParams^*}$, we can extract a lower-dimensional manifold $\manifold$ of the form
\begin{equation}
\manifold := \exp^{\diffeo}_{\bar{\Vector}} \bigl(\operatorname{span} \{\tangentVector_{\bar{\Vector}}^{1}, \ldots, \tangentVector_{\bar{\Vector}}^{\dimIndC}\}\bigr), \quad \tangentVector_{\bar{\Vector}}^{1}, \ldots, \tangentVector_{\bar{\Vector}}^{\dimIndC} \in \tangent_{\bar{\Vector}}\mathbb{R}^\dimInd,
\end{equation}
where ${\bar{\Vector}}\in \Real^\dimInd$ is the data barycentre and the tangent vectors $\tangentVector_{\bar{\Vector}}^{1}, \ldots, \tangentVector_{\bar{\Vector}}^{\dimIndC}$ are obtained via tangent space principal component analysis, or more concretely a low-rank approximation of the linearized data points $\log^{\diffeo}_{\bar{\Vector}} (\Vector^1), \ldots, \log^{\diffeo}_{\bar{\Vector}} (\Vector^\dataPointNum)$.
Crucially, this low-rank approximation must be performed with respect to the $\ell^2$-inner product rather than the Riemannian inner product.

Equipping $\manifold$ with the ambient pullback structure endows it with several favorable geometric properties. By construction, $\manifold$ is a totally geodesic submanifold: for any two points in $\manifold$, every geodesic segment joining them within $\manifold$ is also a geodesic of the ambient space. Consequently, the Riemannian metric on $\manifold$ induced by the ambient pullback structure yields the same manifold mappings. Moreover, geodesics are unique and the exponential map at any point is well-defined for any tangent vector and a global diffeomorphism onto $\manifold$.

\subsection{Iso-Riemannian geometry}

For any Riemannian manifold $(\manifold, (\cdot, \cdot))$ that is embedded in $\Real^\dimInd$, the \emph{iso-connection} $\nabla^{\iso}_{(\cdot)} (\cdot): \vectorfield(\manifold) \times \vectorfield(\manifold)\to \vectorfield(\manifold)$ under $(\manifold, (\cdot, \cdot))$ is defined through the Levi-Civita connection as 
\begin{equation}
    \nabla^{\iso}_{\tangentVectorB_\Vector} \tangentVector := \frac{1}{\|\tangentVectorB_\Vector\|_2} \nabla_{\tangentVectorB_\Vector} \|\mathcal{P}_{(\cdot) \leftarrow \Vector} \tangentVectorB_\Vector\|_2 \tangentVector ,\quad \Vector \in \manifold, \tangentVector, \tangentVectorB \in \vectorfield(\manifold),
    \label{eq:iso-connection}
\end{equation}
and is notably only 1-homogeneous in $\tangentVectorB$ rather than linear.

Similarly to how the Levi-Civita connection generates Riemannian geometry, the iso-connection -- despite being non-linear -- generates iso-Riemannian geometry, i.e., it comes with geodesics $\geodesic^{\iso}_{\Vector, \VectorB}$, exponential mapping $\exp^{\iso}_\Vector$, logarithmic mapping $\log^{\iso}_{\Vector}$ and parallel transport $\partransport^{\iso}_{\VectorB \leftarrow \Vector}$. Notably, these can be written as functions of the manifold mappings under the Levi-Civita connection. In particular, assuming that $(\manifold, (\cdot, \cdot))$ is complete and the exponential map -- the one under the Levi-Civita connection -- is a global diffeomorphism, we have
\begin{align}
    \geodesic^{\iso}_{\Vector, \VectorB}(t) &= \geodesic_{\Vector,\VectorB}(\timechange_{\Vector,\VectorB}(t)),
    \label{eq:thm-geodesic-iso}\\
    \exp^{\iso}_\Vector (\tangentVector_\Vector) &= \exp_\Vector (\vectorchange_\Vector (\tangentVector_\Vector) \tangentVector_\Vector), 
    \label{eq:thm-exp-iso}\\
    \log^{\iso}_{\Vector} (\VectorB) &= \frac{\int_0^1 \|\dot{\geodesic}_{\Vector,\VectorB}(s)\|_2 \, \mathrm{d}s}{\|\log_\Vector (\VectorB)\|_2 }\log_\Vector (\VectorB)
    \label{eq:thm-log-iso}\\
    \partransport^{\iso}_{\VectorB \leftarrow \Vector} \tangentVector_\Vector &= \frac{\|\log_{\Vector}(\VectorB)\|_2}{\|\log_{\VectorB}(\Vector)\|_2} \mathcal{P}_{\VectorB\leftarrow \Vector} ( \tangentVector_\Vector),
    \label{eq:thm-exp-iso}
\end{align}
where the mapping $\timechange_{\Vector,\VectorB}:[0,1]\to[0,1]$ is defined as
\begin{equation}
    \timechange_{\Vector,\VectorB} (t) := \inf\{t' \in [0,1] \mid t \int_0^1 \|\dot{\geodesic}_{\Vector,\VectorB}(s)\|_2 \, \mathrm{d}s = \int_0^{t'} \|\dot{\geodesic}_{\Vector,\VectorB}(s)\|_2 \, \mathrm{d}s\},
\end{equation}
and the mapping $\vectorchange_\Vector : \tangent_\Vector \Real^\dimInd \to \Real$ is defined as
\begin{equation}
    \vectorchange_\Vector(\tangentVector_\Vector) := \inf \{t' \geq 0 \mid \int_{0}^{1}\|\dot{\geodesic}_{\Vector, \exp_\Vector (t' \tangentVector_\Vector)}(s)\|_2 \; \mathrm{d}s = \|\tangentVector_\Vector\|_{2}\}.
\end{equation}
Under these assumptions, we also have
\begin{equation}
    \log^{\iso}_{\Vector} (\exp^{\iso}_\Vector (\tangentVector_\Vector)) = \tangentVector_\Vector, \quad \text{and} \quad \exp^{\iso}_\Vector(\log^{\iso}_{\Vector} (\VectorB)) = \VectorB,
\end{equation}
and 
\begin{equation}
    \dot{\geodesic}^{\iso}_{\Vector, \VectorB}(t) = \frac{\int_0^1 \|\dot{\geodesic}_{\Vector,\VectorB}(s)\|_2 \, \mathrm{d}s}{\|\dot{\geodesic}_{\Vector, \VectorB}(t)\|_2 }\dot{\geodesic}_{\Vector, \VectorB}(t),
    \label{eq:iso-geo-speed}
\end{equation}
which implies that $\dot{\geodesic}^{\iso}_{\Vector, \VectorB}(0) = \log^{\iso}_{\Vector} (\VectorB)$ and that iso-geodesics have a constant $\ell^2$-speed.

Besides the above-mentioned manifold mappings that are directly generated by the iso-connection, there is also the \emph{iso-distance} mapping $\distance^{\iso}_{\manifold}: \manifold\times \manifold\to \Real$ under $(\manifold, (\cdot, \cdot))$, which is defined as
\begin{equation}
    \distance_{\manifold}^{\iso}(\Vector, \VectorB) :=  
\int_0^1 \|\dot{\geodesic}_{\Vector,\VectorB}(s)\|_2 \, \mathrm{d}s
\label{eq:iso-distance}
\end{equation}
and naturally fits in the iso-Riemannian framework as it satisfies
\begin{equation}
     \distance_{\manifold}^{\iso}(\Vector, \VectorB) = \|\log_\Vector^{\iso}(\VectorB)\|_2 \quad \text{and} \quad \distance_{\manifold}^{\iso}(\Vector, \exp_\Vector^{\iso}(\tangentVector_\Vector)) = \|\tangentVector_\Vector\|_{2}.
\end{equation}
However, despite the name, the iso-distance mapping is generally not an actual metric. 

The connection to the Levi-Civita manifold mappings is also of pracitcal importance. In particular, closed-form manifold mappings -- which are available for learned data manifolds $(\manifold, (\cdot, \cdot)) := (\manifold, (\cdot, \cdot)^\diffeo)$ -- can be transformed into computationally efficient iso-manifold mappings, either directly or via additional numerical approximation.



\section{Convergence from a function perspective}
\label{sec:function-convergence}
To establish that problems of the form \cref{eq:mfld-opt-problem-init} can be solved by first-order schemes such as $\ell^2$-PG-IRD in \cref{eq:iso-first-order-scheme}, we extend the classical notions of g-convexity and $L$-smoothness to the iso-Riemannian setting. Specifically, we introduce iso-g-convexity and iso-$L$-smoothness by systematically replacing manifold mappings with iso-manifold mappings, Riemannian inner products with $\ell^2$-inner products, and Riemannian gradients with projected Euclidean gradients. Before stating the main convergence theorem, we confirm that these definitions behave as intended: on $(\mathbb{R}, (\cdot,\cdot))$ they reduce to standard Euclidean convexity and $L$-smoothness; on $(\manifold, (\cdot,\cdot))$ they coincide with the Euclidean notions up to a geometric correction; and strong iso-g-convexity is closely linked to the Polyak–\L{}ojasiewicz (PL) property. We conclude this section by examining a class of least squares problems that arise in inverse problems, and we discuss how the proposed methodology yields improvements over the unconstrained setting. We defer all proofs to \Cref{app:supplement-function-convergence-proofs}.

We note that, in doing so, we naturally show that the conditioning issue illustrated in \Cref{fig:iso-geodesic-convexity-issue-sol} can be resolved within the proposed iso-g-convexity and iso-$L$-smoothness framework, and that our desire to use a different search direction, as depicted in \Cref{fig:eucliden-vs-riemannian-gradients-issue}, can be addressed through the connection between iso-g-convex functions and PL functions. In addition to these theoretical developments, we present numerical experiments illustrating the practical utility of the proposed scheme and considerations relevant to its implementation. We consider totally geodesic submanifolds of both modeled and learned pullback structures; for the latter, we use the off-the-shelf structure and 20-dimensional manifold representation for MNIST from \cite{diepeveen2025manifold}. We defer all details of the numerical experiments to \Cref{app:numerics-function-convergence}.

Finally, we emphasize that we develop our theory for the general case of a Riemannian manifold $(\manifold, (\cdot,\cdot))$ embedded in $\Real^\dimInd$ on which the iso-exponential map is both globally defined (i.e., on the entire tangent bundle) and a global diffeomorphism; this choice keeps the notation simple and ensures that the definitions are well-posed. We remind the reader that the applications of interest are totally geodesic submanifolds of learned Euclidean pullback structures on the entire ambient space. For such spaces, these assumptions are satisfied, and consequently the iterates produced by $\ell^2$-PG-IRD in \cref{eq:iso-first-order-scheme} are always well-defined.

\subsection{Iso-g-convex and iso-L-smooth functions}
\label{sec:iso-g-convex}

We begin by proposing a new notion of geodesic convexity, which will underlie subsequent analysis, and proceed in the same manner for $L$-smoothness.
\begin{definition}[Iso-g-convex functions]
\label{def:iso-g-convex}
    A function $\function:\Real^\dimInd\to \Real$ is \emph{iso-g-convex} (or \emph{iso-geodesically convex}) on a Riemannian manifold $(\manifold, (\cdot,\cdot))$ embedded in $\Real^\dimInd$ on which the iso-exponential map is a global diffeomorphism, if for every $\Vector \neq \VectorB\in \manifold$
    \begin{equation}
        \function(\VectorB) \geq \function(\Vector) + (\nabla \function(\Vector), \log^{\iso}_{\Vector}(\VectorB) )_2,
    \end{equation}
    and \emph{$\alpha$-strongly iso-g-convex} if 
    \begin{equation}
        \function(\VectorB) \geq \function(\Vector) + (\nabla \function(\Vector), \log^{\iso}_{\Vector}(\VectorB) )_2 + \frac{\alpha}{2} \distance^{\iso}_{\manifold} (\Vector, \VectorB)^2.
    \end{equation}
\end{definition}

\begin{definition}[Iso-$L$-smooth functions]
\label{def:iso-g-lip}
    A function $\function:\Real^\dimInd\to \Real$ is \emph{iso-$L$-smooth} on a Riemannian manifold $(\manifold, (\cdot,\cdot))$ embedded in $\Real^\dimInd$ on which the iso-exponential map is globally defined, if for every $\Vector\in \manifold$ and every $\tangentVector_\Vector\in \tangent_\Vector\manifold$
    \begin{equation}
        \function(\exp^{\iso}_\Vector(\tangentVector_\Vector) ) \leq \function(\Vector) + (\nabla \function(\Vector), \tangentVector_\Vector )_2 + \frac{L}{2}\|\tangentVector_\Vector\|_2^2.
    \end{equation}
\end{definition}

As a direct sanity check of the definitions, Euclidean convexity on $\Real$ is equivalent to iso-g-convexity and a similar equivalence holds for $L$-smoothness, as one would naturally expect.

\begin{proposition}[Equivalences on $\Real$]
\label{prop:equivalence-iso-convex-smooth-R}
    Let $f:\Real \to \Real$ be a real-valued function and consider any Riemannian structure $(\Real,(\cdot, \cdot))$. 

    Then:
    \begin{enumerate}[label=(\roman*)]
        \item $\function$ is ($\alpha$-strongly) convex on $\Real$ if and only if it is ($\alpha$-strongly) iso-g-convex on $(\Real, (\cdot, \cdot))$.
        \item $\function$ is $L$-smooth on $\Real$ if and only if it is iso-$L$-smooth on $(\Real, (\cdot, \cdot))$.
    \end{enumerate}
\end{proposition}

Since we are in the end interested in considering functions on manifolds, it will be helpful to have a different way of checking iso-g-convexity and iso-$L$-smoothness than directly from the definitions. In particular, when the function $\function$ is twice continuously differentiable, there is a more convenient criterion.

\begin{theorem}[Iso-g-convex and iso-$L$-smooth differentiable functions]
\label{thm:iso-convex-smooth-impl}
    Let $\function:\Real^\dimInd\to \Real$ be a function that is twice continuously differentiable and let $(\manifold, (\cdot,\cdot))$ be a Riemannian manifold embedded in $\Real^\dimInd$ on which the iso-exponential map is both globally defined and a global diffeomorphism. The following implications hold:
    \begin{enumerate}[label=(\roman*)]
        \item If for every $\Vector \neq \VectorB\in \manifold$
        \begin{equation}
            \frac{\mathrm{d}^2}{\mathrm{d}t^2} \function (\geodesic^{\iso}_{\Vector,\VectorB}(t)) \geq 0, \quad t\in [0,1],
        \end{equation}
        then $\function$ is iso-g-convex on $(\manifold, (\cdot,\cdot))$.
    \item If for every $\Vector \neq \VectorB\in \manifold$
    \begin{equation}
        \frac{\mathrm{d}^2}{\mathrm{d}t^2} \function (\geodesic^{\iso}_{\Vector,\VectorB}(t)) \geq \alpha \distance^{\iso}_{\manifold} (\Vector, \VectorB)^2, \quad t\in [0,1],
    \end{equation}
    then $\function$ is $\alpha$-strongly iso-g-convex on $(\manifold, (\cdot,\cdot))$.
    \item If for every $\Vector\neq \VectorB \in \manifold$ 
    \begin{equation}
        \frac{\mathrm{d}^2}{\mathrm{d}t^2} \function (\geodesic^{\iso}_{\Vector,\VectorB}(t)) \leq L \distance^{\iso}_{\manifold} (\Vector, \VectorB)^2, \quad t\in [0,1],
    \end{equation}
    then $\function$ is iso-$L$-smooth on $(\manifold, (\cdot,\cdot))$.
    \end{enumerate}
\end{theorem}

We can now use the preceding result to examine more closely the relationship between Euclidean convex and $L$-smooth functions and their iso-Riemannian counterparts beyond the one-dimensional case in \Cref{prop:equivalence-iso-convex-smooth-R}. To that end, we begin by evaluating the second derivative along an iso-geodesic, as summarized in the following result.

\begin{lemma}
\label{lem:second-deriv-along-iso-geo}
    Let $\function:\Real^\dimInd\to \Real$ be a function that is twice continuously differentiable and let $(\manifold, (\cdot,\cdot))$ be a Riemannian manifold embedded in $\Real^\dimInd$ on which the iso-exponential map is both globally defined and a global diffeomorphism.

    Then, for every  $\Vector \neq \VectorB\in \manifold$
    \begin{align}
        \frac{\mathrm{d}^2}{\mathrm{d}t^2} & \function (\geodesic^{\iso}_{\Vector,\VectorB}(t))=
        \\
        &\frac{D_{\VectorC_t}^2 \function [\tangentVector_{\VectorC_t}, \tangentVector_{\VectorC_t}] - \bigl(\bigl(\mathbf{I} - \frac{\tangentVector_{\VectorC_t}}{\|\tangentVector_{\VectorC_t}\|_2} \otimes \frac{\tangentVector_{\VectorC_t}}{\|\tangentVector_{\VectorC_t}\|_2} \bigr) \nabla \function (\VectorC_t), \Gamma_{\VectorC_t} [\tangentVector_{\VectorC_t}, \tangentVector_{\VectorC_t}]\bigr)_2}{\|\tangentVector_{\VectorC_t}\|_2^2}\distance^{\iso}_{\manifold} (\Vector,\VectorB)^2, \nonumber
    \end{align}
    where $\tangentVector_{\VectorC_t} := \dot{\geodesic}^{\iso}_{\Vector,\VectorB}(t)$ and $ \VectorC_t = \geodesic^{\iso}_{\Vector,\VectorB}(t)$.
\end{lemma}



Next, using \Cref{thm:iso-convex-smooth-impl} through \Cref{lem:second-deriv-along-iso-geo}, it naturally follows that Euclidean convex and $L$-smooth functions are related to iso-g-convex and iso-$L$-smooth functions via a geometric correction; conversely, with an appropriate geometric correction, certain non-convex functions may become iso-g-convex. 


\begin{proposition}[Implications on $\manifold$]
    \label{prop:implications-on-manifold}
    Let $\function:\Real^\dimInd\to \Real$ be a function that is twice continuously differentiable and let $(\manifold, (\cdot,\cdot))$ be a Riemannian manifold embedded in $\Real^\dimInd$ on which the iso-exponential map is both globally defined and a global diffeomorphism. Furthermore, assume that 
    \begin{equation}
        \frac{D_{\Vector}^2 \function [\tangentVector_{\Vector}, \tangentVector_{\Vector}]}{\|\tangentVector_{\Vector}\|_2^2} \in [\beta, \eta],
    \end{equation}
    and
    \begin{equation}
    - \frac{( (\mathbf{I} - \frac{\tangentVector_{\Vector}}{\|\tangentVector_{\Vector}\|_2} \otimes \frac{\tangentVector_{\Vector}}{\|\tangentVector_{\Vector}\|_2})\nabla \function (\Vector), \Gamma_{\Vector} (\tangentVector_{\Vector}, \tangentVector_{\Vector}) )_2}{\|\tangentVector_{\Vector}\|_2^2} \in [\delta, \zeta],
    \label{eq:thm-geo-factor}
    \end{equation}
    where $\Gamma_{\Vector}:\tangent_{\Vector} \manifold \times \tangent_{\Vector} \manifold \to \tangent_{\Vector} \Real^{\dimInd}$ is the Christoffel operator on $(\manifold, (\cdot, \cdot))$. Finally, assume that $\beta + \delta >0$. 
    
    Then, the following hold:
    \begin{enumerate}[label=(\roman*)]
        \item $\function$ is $(\beta + \delta)$-strongly iso-g-convex on $(\manifold, (\cdot,\cdot))$.
        \item $\function$ is iso-$(\eta + \zeta)$-smooth on $(\manifold, (\cdot,\cdot))$.
    \end{enumerate}
\end{proposition}

\begin{remark}
\label{rem:non-convex-iso-g-convex}
    We indeed see that strongly convex functions with $\beta > 0$ remain strongly iso-convex as long as $\beta + \delta > 0$; more strikingly, under an appropriate geometric correction, functions that are not strongly convex or even non-convex -- with $\beta \leq 0$ -- can become strongly iso-g-convex provided $\delta > |\beta| \geq 0$.
\end{remark}

\begin{remark}
    We also note that this result generalizes aspects of the univariate real-valued case in \Cref{prop:equivalence-iso-convex-smooth-R}. In particular, $\delta=\zeta=0$ for any Riemannian structure on $\Real$, because the projection matrix always vanishes: since $\mathbf{I} = 1$ and $\xi_{x} \in \tangent_x \Real \cong \Real$
    \begin{equation}
        \mathbf{I} - \frac{\xi_x}{\|\xi_x\|_2} \otimes \frac{\xi_x}{\|\xi_x\|_2} = 1 - \frac{\xi_x\cdot \xi_x}{\|\xi_x\|_2^2} = 1-1 = 0.
    \end{equation}
    So is-g-convexity and iso-$L$-smoothness is fully determined by the second derivative, which yields one direction of \Cref{prop:equivalence-iso-convex-smooth-R}, i.e., that convex and $L$-smooth functions are iso-g-convex and iso-$L$-smooth on $(\Real, (\cdot, \cdot))$.
\end{remark}

\begin{remark}
\label{rem:christoffel-interaction}
    It is worth noting what the factor in \cref{eq:thm-geo-factor} encodes and how it relates to the issue highlighted in \Cref{fig:geodesic-convexity-issue}. Since the Christoffel operator precisely quantifies the ambient acceleration -- or, equivalently, the change in geodesic speed -- at a given point, we can interpret this geometric correction term as the interaction between geodesic acceleration and the gradient of the objective, with the component of acceleration in the direction of the geodesic velocity excluded. Consequently, speed-ups or slow-downs along a geodesic do not affect the conditioning of the problem; instead, the geometry remains as faithful as possible to the underlying Euclidean convexity and $L$-smoothness, which is precisely the discrepancy illustrated by \Cref{fig:geodesic-convexity-issue,fig:iso-geodesic-convexity-issue-sol}.
    
\end{remark}

\subsection{Finding minimizers of iso-g-convex and iso-L-smooth functions}
\label{sec:finding-mins-g-convex}

Next, to establish convergence of $\ell^2$-PG-IRD in \cref{eq:iso-first-order-scheme} for solving \cref{eq:mfld-opt-problem-init}, we generalize the classical proof for functions satisfying the PL property and $L$-smoothness. To that end, we first introduce an appropriate notion of the PL property based on a projected Euclidean gradient rather than a Riemannian gradient.

\begin{definition}[PL-functions]
    A function $\function:\Real^\dimInd\to \Real$ is \emph{\(\mu\)-PL} on a smooth manifold $\manifold \subset \Real^\dimInd$, if for every $\Vector\in \manifold$
    \begin{equation}
    \function(\Vector) - \function^* \leq \frac{1}{2\mu} \|\mathbf{P}_\Vector \nabla \function (\Vector)\|_2^2,
    \end{equation}
    where $\function^* = \inf_{\Vector\in\manifold} \function(\Vector)$ and  $\mathbf{P}_{\Vector}:\tangent_{\Vector} \Real^{\dimInd}\to\tangent_{\Vector} \manifold$ is the ($\ell^2$-)orthogonal projection operator onto the tangent space $\tangent_{\Vector} \manifold$ at $\Vector\in \manifold$.
\end{definition}

Although our definition does not involve iso-Riemannian geometry at all, $\alpha$-strongly iso-g-convex functions satisfy the PL property.

\begin{theorem}[Strongly iso-g-convex functions are PL]
\label{lem:iso-g-convex-implies-PL}
Let $\function:\Real^\dimInd\to \Real$ be a \(\alpha\)-strongly iso-g-convex function on a Riemannian manifold $(\manifold, (\cdot,\cdot))$ embedded in $\Real^\dimInd$ on which the iso-exponential map is both globally defined and a global diffeomorphism. 

Then, $\function$ is \(\alpha\)-PL.
\end{theorem}

With these definitions in place, we obtain sufficient conditions for convergence. Since our focus is on the convergence of iterates, we state the following result under the assumptions of strong iso-g-convexity and iso-$L$-smoothness. The result also characterizes when our scheme converges to a minimizer of a Euclidean convex and $L$-smooth function, while its proof further yields convergence in $\function$ under the weaker assumptions that $\function$ satisfies the PL property and is $L$-smooth. Both consequences are discussed in the subsequent corollaries.

\begin{theorem}[Convergence of $\ell^2$-projected gradient iso-Riemannian descent]
\label{thm:conv-pg-ird}
    Let $\function:\Real^\dimInd\to \Real$ be a $\alpha$-strongly iso-g-convex and iso-$L$-smooth function on a Riemannian manifold $(\manifold, (\cdot,\cdot))$ embedded in $\Real^\dimInd$ on which the iso-exponential map is both globally defined and a global diffeomorphism. Furthermore, assume that there is a point $\bar{\Vector} \in \manifold$ that satisfies $\mathbf{P}_{\bar{\Vector}} \nabla f (\bar{\Vector}) = \mathbf{0}$. 
    
    Then, the sequence $\{\Vector^{(\sumIndC)}\}_{\sumIndC=0}^\infty \in \manifold$ defined by the algorithm
    \[
        \Vector^{(\sumIndC+1)}=\exp^{\iso}_{\Vector^{(\sumIndC)}}(-r \mathbf{P}_{\Vector^{(\sumIndC)}} \nabla f (\Vector^{(\sumIndC)})), \quad \Vector^{(0)} \in \manifold,
    \]
    where \(0 < r \leq \frac{1}{L}\) converges linearly in \(\ell^2\) to $\bar{\Vector} \in \manifold$.
\end{theorem}

\begin{corollary}
\label{cor:iso-convex}
    Under the assumptions of \Cref{prop:implications-on-manifold} -- notably assuming that $\beta + \delta >0$ --, suppose additionally that there exists a point $\bar{\Vector} \in \manifold$ satisfying $\mathbf{P}_{\bar{\Vector}} \nabla f(\bar{\Vector}) = \mathbf{0}.$
    
    Then, the sequence $\{\Vector^{(\sumIndC)}\}_{\sumIndC=0}^\infty \in \manifold$ defined by the algorithm
    \[
        \Vector^{(\sumIndC+1)}=\exp^{\iso}_{\Vector^{(\sumIndC)}}(-r \mathbf{P}_{\Vector^{(\sumIndC)}} \nabla f (\Vector^{(\sumIndC)})), \quad \Vector^{(0)} \in \manifold,
    \]
    where \(0 < r \leq \frac{1}{\eta + \zeta}\) converges linearly in \(\ell^2\) to $\bar{\Vector} \in \manifold$.    
\end{corollary}

\begin{corollary}
\label{cor:iso-pl}
    Let $\function:\Real^\dimInd\to \Real$ be a $\mu$-PL and iso-$L$-smooth function on a Riemannian manifold $(\manifold, (\cdot,\cdot))$ embedded in $\Real^\dimInd$ on which the iso-exponential map is both globally defined and a global diffeomorphism. 
    
    Then, the sequence $\{\function(\Vector^{(\sumIndC)})\}_{\sumIndC=0}^\infty \in \Real$ defined by the algorithm
    \[
        \Vector^{(\sumIndC+1)}=\exp^{\iso}_{\Vector^{(\sumIndC)}}(-r \mathbf{P}_{\Vector^{(\sumIndC)}} \nabla f (\Vector^{(\sumIndC)})), \quad \Vector^{(0)} \in \manifold,
    \]
    where \(0 < r \leq \frac{1}{L}\) converges linearly to the minimum value $\function^*$.
\end{corollary}

\begin{remark}
\label{rem:non-convex-2}
    It is useful to interpret \Cref{rem:non-convex-iso-g-convex} in the context of \Cref{cor:iso-convex}: our algorithm converges whenever the geometric correction term preserves strong convexity, or, more generally, whenever the resulting function is strongly iso-g-convex -- even if $\function$ is not strongly convex in the Euclidean sense. Conversely, a Euclidean convex function may fail to be iso-convex when the geometric correction term acts unfavourably.
\end{remark}

\begin{remark}
    Although the constant \(\eta\) in \Cref{cor:iso-convex} may be known or readily estimated for a given problem, obtaining a practical estimate of \(\zeta\) is generally more challenging. Consequently, to make \(\ell^2\)-PG-IRD a practical algorithm, it may be necessary to incorporate a line-search procedure.
\end{remark}

\begin{remark}
    In the context of \Cref{cor:iso-pl}, it is important to emphasize that the class of PL functions is substantially broader than the class of iso-g-convex functions. In particular, the objective function in \Cref{fig:eucliden-vs-riemannian-gradients-issue}, considered with its Euclidean gradient, satisfies the PL property globally and is iso-$L$-smooth, despite not being iso-g-convex everywhere, as shown in \Cref{fig:iso-geodesic-convexity-issue-sol}. Thus, the discrepancy illustrated in \Cref{fig:eucliden-vs-riemannian-gradients-issue} fits naturally within our framework, even though the underlying objective is not iso-g-convex.
    
    At the same time, this example highlights a limitation of the current theory and suggests directions for future work. Because the manifold in the figure is open, the iso-exponential map is not globally defined, and an iterate may leave the manifold. Consequently, the algorithm may fail to converge -- even in $\function$ -- unless the relevant conditions are carefully controlled. This issue can arise even when the objective is iso-g-convex and is therefore not merely a consequence of the absence of iso-g-convexity. Addressing this situation would require incorporating constraints that keep the iterates on the manifold, which lies beyond the scope of the current framework. The figure is therefore primarily illustrative and points toward future directions. By contrast, the motivating setting of this work -- namely, a totally geodesic submanifold of an ambient Euclidean pullback structure for which the iso-exponential map is globally defined and is a global diffeomorphism -- is covered by our result.
\end{remark}



\subsection{Implications for solving inverse problems}
\label{sec:inverse-problems-appl}

As an application of the above theory, we will consider solving inverse problems constrained to a manifold. That is, we will focus on optimization problems of the form
\begin{equation}
    \inf_{\Vector\in \manifold} \frac{1}{2} \|\forward \Vector -\source \|_2^2, 
\end{equation}
where the matrix $\forward\in \Real^{\dimIndB\times \dimInd}$ is the forward operator, the vector $\source \in \Real^{\dimIndB}$ is the observed measurement or corrupted data, and $\manifold \subset \Real^\dimInd$ is a smooth manifold, which we assume to be a Riemannian manifold $(\manifold, (\cdot,\cdot))$ on which the iso-exponential map is both globally defined and a global diffeomorphism.

\paragraph{Practical considerations}
First, to see the above ideas in practice, we will consider experiments with both modeled and learned pullback structures to highlight some the discussion in \Cref{rem:non-convex-2}. 

\begin{figure}[h!]
    \centering
    \setlength{\tabcolsep}{2pt}
    \renewcommand{\arraystretch}{0}
    \begin{tabular}{ccc}
        \includegraphics[width=0.32\linewidth]{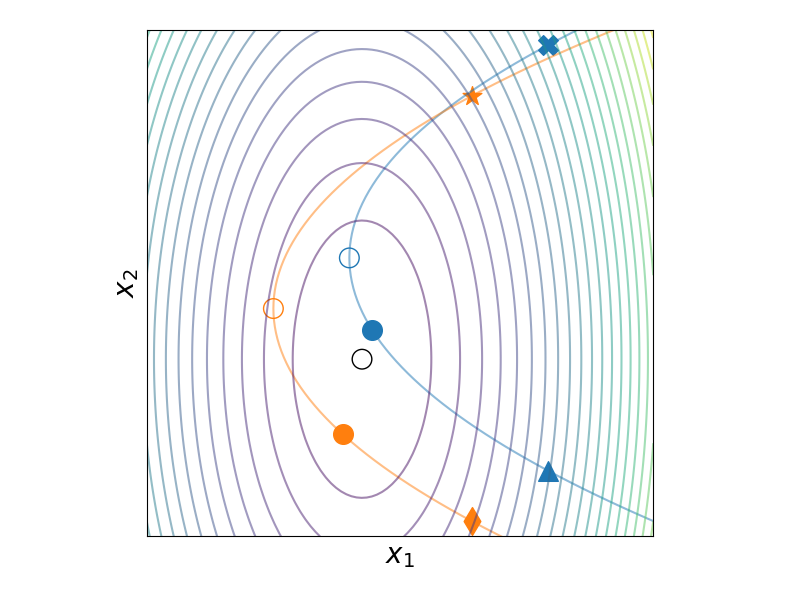} &
        \includegraphics[width=0.32\textwidth]{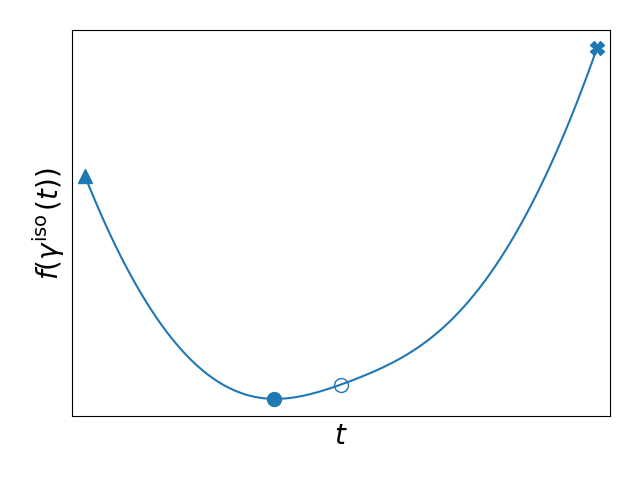} &
        \includegraphics[width=0.32\textwidth]{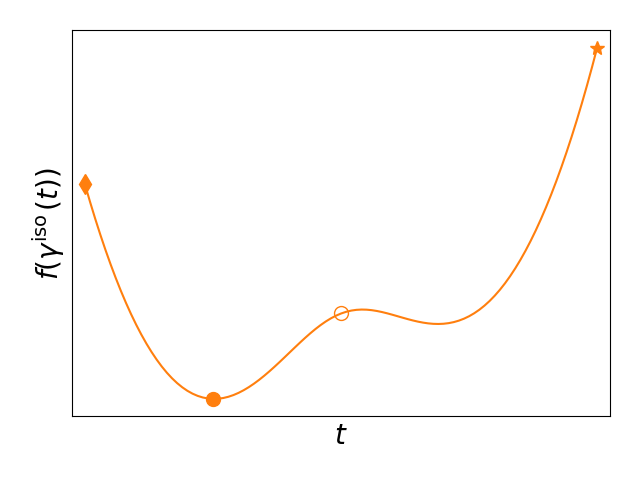} \\
        \makebox[0.29\textwidth]{\small Optimization landscape} &
        \makebox[0.29\textwidth]{\small Objective on blue manifold} &
        \makebox[0.29\textwidth]{\small Objective on orange manifold}
    \end{tabular}
    
    \caption{Optimization results of an inverse problem (left) -- visualized by its level sets -- under two different modeled manifold constraints (blue and orange) and the objective along iso-geodesics between respective marked end points (middle and right). For both cases $\ell^2$PG-IRD \cref{eq:iso-first-order-scheme} initialized at the respective unfilled circle markers converge to global minimizers denoted by the filled circle markers, which are notably distinct from the unconstraint solution denoted by the black unfilled circle marker. Notably, despite the identical shapes of the constraint manifolds, the optimization problem over the orange manifold is not iso-g-convex, in contrast to the blue manifold.}
    \label{fig:opt-modeled-r2}
\end{figure}

For a toy inverse problem on a one-dimensional manifold in the modeled pullback setting, \Cref{fig:opt-modeled-r2} shows that $\ell^2$-PG-IRD converges to the global minimizer under both manifold constraints. In a neighborhood of the minimizers, the convergence is linear. 

It is worth highlighting that although the two manifolds have the same shape and differ only by an offset, the corresponding optimization problems have fundamentally different properties. In particular, the blue manifold preserves the convexity of the original objective, since $\beta+\delta>0$, whereas the orange manifold does not; consequently, the objective is no longer iso-convex on the latter, and because of which we observe linear convergence only locally.

\begin{figure}[h!]
    \centering
    \setlength{\tabcolsep}{2pt}
    \renewcommand{\arraystretch}{0}
    \begin{tabular}{ccc}
        \includegraphics[width=0.29\textwidth]{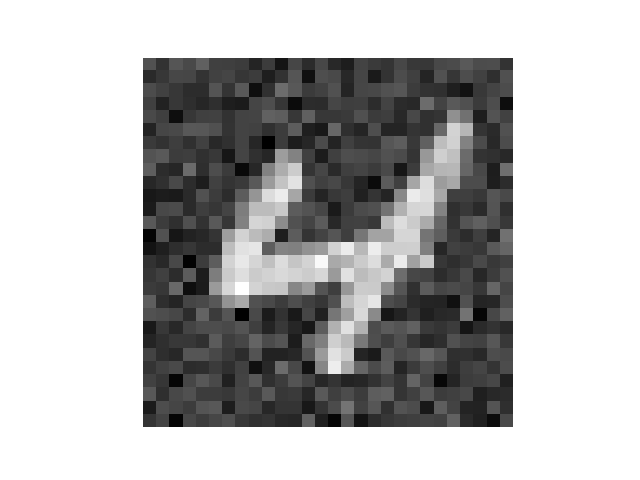} &
        \includegraphics[width=0.29\textwidth]{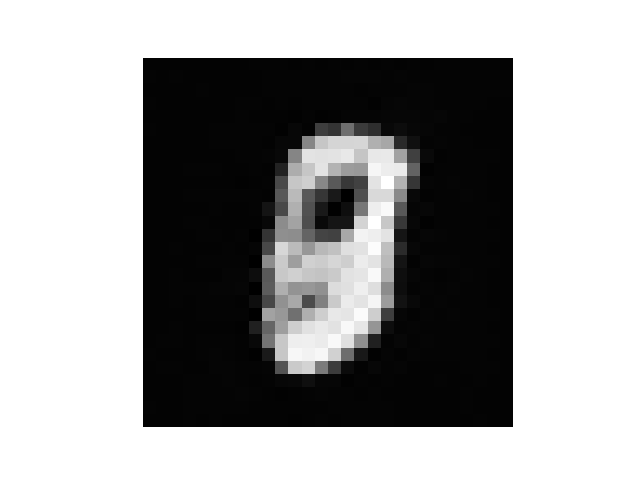} &
        \includegraphics[width=0.29\textwidth]{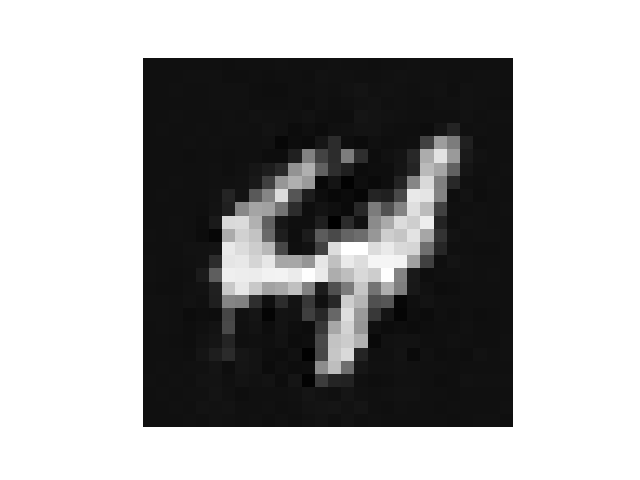} \\
        \makebox[0.29\textwidth]{\small Noisy data} &
        \makebox[0.29\textwidth]{\small Initialization} &
        \makebox[0.29\textwidth]{\small Minimizer}
    \end{tabular}
    \caption{Optimization results of denoising the number four (left) under a learned geodesic submanifold constraint. The proposed $\ell^2$-PG-IRD \cref{eq:iso-first-order-scheme} initialized at the data barycentre (center) converges on a 20-dimensional learned subspace to a noiseless solution (right) after only 10 iterations.}
    \label{fig:mnist-results-opt}
\end{figure}

In the learned pullback setting, we consider a denoising inverse problem on a learned 20-dimensional manifold equipped with an off-the-shelf, pretrained pullback metric \cite{diepeveen2025manifold}. As shown in \Cref{fig:mnist-results-opt}, $\ell^2$-PG-IRD converges to a nearly noiseless representation of the digit four. The convergence is again linear. However, unlike in the preceding toy examples, it is difficult to determine whether the objective is globally iso-convex, as for the blue manifold, or whether iso-convexity holds only locally, as for the orange manifold in \Cref{fig:opt-modeled-r2}. 


\paragraph{Improved conditioning}

To see how the Riemannian case can attain superior performance compared to solving the inverse problem without the manifold constraint\footnote{That is, without even assuming that having a data manifold prior will give qualitatively better solutions.}, it is worth highlighting the following two settings in more detail: overdetermined and underdetermined $\forward$.

For overdetermined $\forward$, we assume that $\dimIndB \geq \dimInd$ and $\operatorname{rank}(\mathbf{A}) = \dimInd$, i.e., least squares fitting. If $\dim (\manifold) < \dimInd$, we have  $\beta \geq \lambda_{\min} (\mathbf{A}^\top \mathbf{A}) >0$ and $\eta \leq \lambda_{\max} (\mathbf{A}^\top \mathbf{A})$. So, as long as $\beta + \delta >0$ and $\eta + \zeta \leq \lambda_{\max} (\mathbf{A}^\top \mathbf{A})$ we know that we can take a bigger step size and still converge linearly, but notably faster compared to the unconstraint case. Such a scenario arises when the the manifold $\manifold$ lies within a (higher dimensional) linear subspace that is orthogonal to the subspace spanned by the eigenvectors with the largest eigenvalues. In such a case $\eta$ can be significantly smaller than $\lambda_{\max} (\mathbf{A}^\top \mathbf{A})$, leading to improved performance provided the manifold geometry does not negate this advantage.

For underdetermined $\forward$, we assume that that $\dimIndB \leq \dimInd$ and $\operatorname{rank}(\mathbf{A}) = \dimIndB$. We note that this problem needs additional regularization in order to be a well-posed problem and is also the more realistic inverse problem setting -- as there are typically less measurements than variables to solve for. The same line of reasoning as above holds, where we now explicitly know that $\lambda_{\min} (\mathbf{A}^\top \mathbf{A}) = 0$. In other words, if $\beta +\delta >0$, the manifold constraint has turned the ill-posed problem into a well-posed one, while $\ell^2$-PG-IRD in \cref{eq:iso-first-order-scheme} enables us efficient computation of a minimizer. To give a concrete example of a scenario where this might happen, we consider compressed sensing. In particular, if we assume restricted isometry in the sense that there exists a $\varepsilon \in [0,1)$ such that for every $\Vector \neq \VectorB \in \manifold$
\begin{equation}
    (1 - \varepsilon) \distance^{\iso}_{\manifold}(\Vector, \VectorB)^2 \leq \|\mathbf{A} \log^{\iso}_{\Vector}(\VectorB)\|_2^2 \leq (1 + \varepsilon) \distance^{\iso}_{\manifold}(\Vector, \VectorB)^2,
\end{equation}
on the manifold we get $\beta = 1 - \varepsilon$ and $\eta = 1 + \varepsilon$. So as long as the potential non-convexity induced by the geometry of the manifold does not dominate, i.e., we have $1 - \varepsilon + \delta>0$, we can efficiently compute the solution to the inverse problem.

\section{Convergence from a vector field perspective}
\label{sec:vectorfield-convergence}

The classical Riemannian optimization literature offers another well-known route to convergence: establishing the monotonicity and Lipschitz continuity of a vector field $\tangentVector \in \vectorfield(\manifold)$ along which the iterates evolve in order to find a singularity. In this section, we extend these notions to the iso-Riemannian setting by introducing iso-monotonicity and iso-Lipschitzness. Before presenting the corresponding convergence result -- and considering the case of the projected gradient field $\tangentVector := \mathbf{P}_{(\cdot)} \nabla \function$ --, we examine these definitions in greater detail. Unlike in the previous section, we find it difficult to obtain meaningful results beyond the one-dimensional setting. On $(\Real,(\cdot,\cdot))$ and on one-dimensional manifolds, we obtain results analogous to those in \Cref{prop:equivalence-iso-convex-smooth-R,prop:implications-on-manifold}, and find that iso-monotonicity and iso-Lipschitzness are equivalent to iso-g-convexity and iso-$L$-smoothness, respectively. Beyond the one-dimensional setting, however, there is little more that can be said at present. In fact, we expect that imposing iso-monotonicity and iso-Lipschitzness in higher dimensions may be overly restrictive. We defer all proofs to \Cref{app:supplement-vector-field-convergence-proofs}.

More generally, the classical Riemannian equivalences do not necessarily carry over to the iso-Riemannian setting, which is what we aim to highlight in this section. In particular, although $\alpha$-strongly g-convex functions have $\alpha$-strongly monotone Riemannian gradients, and $L$-smooth functions have $L$-Lipschitz Riemannian gradients, the corresponding implications need not hold in the iso-Riemannian case. Thus, this section also highlights that iso-Riemannian optimization admits multiple ways of generalizing classical optimization schemes, with each approach yielding convergence results for a potentially different class of functions. 

Despite the weaker theoretical guarantees obtained here for solving problems of the form \cref{eq:mfld-opt-problem-init}, this framework is applicable to problems that are naturally formulated as finding zeros of vector fields. One such application is the construction of iso-Riemannian barycentres, which have not yet been proposed within the iso-Riemannian geometry framework \cite{diepeveen2025manifold}. The tools developed in the previous section cannot be applied to these problems. Nevertheless, our approach yields practical algorithms that extend beyond the scope of the current theory -- again, both for modeled and learned manifolds\footnote{In the latter case, we again use the off-the-shelf learned structure and the 20-dimensional manifold representation for MNIST from \cite{diepeveen2025manifold}.} -- and it readily supports downstream tasks such as iso-\(K\)-means clustering, which we discuss in more detail in \Cref{sec:iso-k-means}. We defer all details of the numerical experiments to \Cref{sec:app-iso-bary}.

Finally, we emphasize again that we develop our theory for the general case of a Riemannian manifold $(\manifold, (\cdot,\cdot))$ embedded in $\Real^\dimInd$ on which the iso-exponential map is both globally defined (i.e., on the entire tangent bundle) and a global diffeomorphism; as before, this choice keeps the notation simple and ensures that the definitions are well-posed. We remind the reader that the applications of interest are totally geodesic submanifolds of learned Euclidean pullback structures on the entire ambient space. 
For such spaces, these assumptions are satisfied, and consequently the iterates produced by the \emph{iso-Riemannian descent} (IRD) algorithm
\begin{equation}
    \Vector^{(\sumIndC+1)} := \exp^{\iso}_{\Vector^{(\sumIndC)}} (- r \tangentVector_{\Vector^{(\sumIndC)}}), \quad r>0, \;\Vector^{(0)}\in \manifold,
    \label{eq:ird}
\end{equation}
are always well-defined.

\subsection{Iso-monotone and iso-Lipschitz vector fields}
\label{sec:iso-monotone}

Here too, we begin by proposing a new notion of monotonicity, which will underlie subsequent analysis, and proceed in the same manner for Lipschitz continuity.

\begin{definition}[Iso-monotone vector fields]
\label{def:iso-monotone-vf}
    A vector field $\tangentVector \in \vectorfield(\manifold)$ is \emph{iso-monotone} on a Riemannian manifold $(\manifold, (\cdot, \cdot)$ embedded in $\Real^\dimInd$ on which the iso-exponential map is a global diffeomorphism,
    if for every $\Vector\neq \VectorB\in \manifold$
\begin{equation}
    (\tangentVector_{\VectorB} - \partransport^{\iso}_{\VectorB\leftarrow \Vector} \tangentVector_{\Vector}, \partransport^{\iso}_{\VectorB\leftarrow \Vector} \log^{\iso}_{\Vector} (\VectorB))_2 \geq 0,
\end{equation}
and \emph{$\alpha$-strongly iso-monotone} if
\begin{equation}
    (\tangentVector_{\VectorB} - \partransport^{\iso}_{\VectorB\leftarrow \Vector} \tangentVector_{\Vector}, \partransport^{\iso}_{\VectorB\leftarrow \Vector} \log^{\iso}_{\Vector} (\VectorB))_2 \geq \alpha \distance^{\iso}_{\manifold} (\Vector, \VectorB)^2, \quad \alpha \geq 0.
\end{equation}
\end{definition}

\begin{definition}[Iso-Lipschitz vector fields]
\label{def:iso-lipschitz-vf}
    A vector field $\tangentVector \in \vectorfield(\manifold)$ is \emph{iso-$L$-Lipschitz} on a Riemannian manifold $(\manifold, (\cdot, \cdot)$ embedded in $\Real^\dimInd$ on which the iso-exponential map is a global diffeomorphism, if for every $\Vector\neq \VectorB\in \manifold$ 
\begin{equation}
     \| \tangentVector_{\VectorB} - \partransport^{\iso}_{\VectorB\leftarrow \Vector}\tangentVector_{\Vector} \|_2 \leq L \distance^{\iso}_{\manifold}(\Vector, \VectorB), \quad L \geq 0.
\end{equation}
\end{definition}

As a direct sanity check of the definitions, Euclidean monotonicity on $\Real$ is equivalent to iso-monotonicity and a similar equivalence holds for Lipschitzness, as one would naturally expect. In addition, we have that iso-monotonicity and iso-g-convexity are equivalent on $\Real$ and vice versa for iso-$L$-Lipschitzness and iso-$L$-smoothness.

\begin{proposition}[Equivalences on $\Real$]
\label{prop:equivalence-iso-monotone-smooth-R}
    Let $\function:\Real \to \Real$ be a continuously differentiable real-valued function and consider any Riemannian structure $(\Real,(\cdot, \cdot))$. 

    Then:
    \begin{enumerate}[label=(\roman*)]
        \item \(\function'\) is (\(\alpha\)-strongly) monotone on \(\Real\) if and
only if it is (\(\alpha\)-strongly) iso-monotone on
\((\Real,(\cdot,\cdot))\).
        \item $\function'$ is $L$-Lipschitz on $\Real$ if and only if it is iso-$L$-Lipschitz on $(\Real, (\cdot, \cdot))$.
    \end{enumerate}
\end{proposition}

\begin{corollary}
     Let $\function:\Real \to \Real$ be a continuously differentiable real-valued function and consider any Riemannian structure $(\Real,(\cdot, \cdot))$. 

    Then:
    \begin{enumerate}[label=(\roman*)]
        \item $f$ is ($\alpha$-strongly) iso-g-convex on $(\Real, (\cdot, \cdot))$ if and only if $f'$ is ($\alpha$-strongly) iso-monotone on $(\Real, (\cdot, \cdot))$.
        \item $f$ is iso-$L$-smooth on $(\Real, (\cdot, \cdot))$ if and only if $f'$ is iso-$L$-Lipschitz on $(\Real, (\cdot, \cdot))$.
    \end{enumerate}
\end{corollary}


To verify iso-monotonicity and iso-Lipschitzness -- and to take a closer look at equivalence with iso-g-convexity and iso-$L$-smoothness -- there is a more convenient way, if the vector field is differentiable. However, to do that, we need the following result, which is an iso version of the Fundamental Theorem of Calculus for vector fields.

\begin{lemma}
\label{prop:taylor-iso}
Let $\tangentVector \in \vectorfield(\manifold)$ be a smooth vector field on a Riemannian manifold $(\manifold, (\cdot,\cdot))$ embedded in $\Real^\dimInd$ on which the iso-exponential map is both globally defined and a global diffeomorphism.

Then, for every $\Vector\neq \VectorB\in \manifold$ 
    \begin{equation}
        \tangentVector_{\VectorB} - \partransport^{\iso}_{\VectorB\leftarrow \Vector} \tangentVector_{\Vector} = \int_{0}^{1} \frac{\|\dot{\geodesic}_{\Vector, \VectorB}(t)\|_2}{\|\log_{\VectorB} (\Vector)\|_2} \partransport_{\VectorB\leftarrow \geodesic_{\Vector, \VectorB}(t)}   \nabla^{\iso}_{\dot{\geodesic}_{\Vector, \VectorB}(t)} \tangentVector_{(\cdot)} \; \mathrm{d}t.
        \label{eq:lem-taylor-iso}
    \end{equation}
\end{lemma}

Having \Cref{prop:taylor-iso} at our disposal, our more convenient criterion follows.

\begin{theorem}[Iso-monotone and iso-Lipschitz differentiable vector fields]
\label{thm:iso-monotone-impl}
    Let $\tangentVector \in \vectorfield(\manifold)$ be a vector field on a Riemannian manifold $(\manifold, (\cdot,\cdot))$ embedded in $\Real^\dimInd$ on which the iso-exponential map is both globally defined and a global diffeomorphism. The following implications hold:
    \begin{enumerate}[label=(\roman*)]
        \item If for every $\Vector\neq \VectorB\in \manifold$ 
        \begin{equation}
        (\partransport_{\VectorB\leftarrow \geodesic_{\Vector, \VectorB}(t)} \nabla^{\iso}_{\dot{\geodesic}_{\Vector, \VectorB}(t)}\tangentVector_{(\cdot)}, \partransport_{\VectorB\leftarrow \Vector} \log_{\Vector} (\VectorB))_2 \geq 0, \quad t \in [0,1],
        \label{eq:lem-iso-mono-impl}
    \end{equation}
        then $\tangentVector$ is iso-monotone on $(\manifold, (\cdot,\cdot))$.
    \item If for every $\Vector\neq \VectorB\in \manifold$ 
    \begin{equation}
        (\partransport_{\VectorB\leftarrow \geodesic_{\Vector, \VectorB}(t)} \nabla^{\iso}_{\dot{\geodesic}_{\Vector, \VectorB}(t)}\tangentVector_{(\cdot)}, \partransport_{\VectorB\leftarrow \Vector} \log_{\Vector} (\VectorB))_2 \geq \alpha  \|\log_{\VectorB}(\Vector)\|_2^2, \quad t\in [0,1],
        \label{eq:lem-iso-alpha-mono-impl}
    \end{equation}
    then $\tangentVector$ is $\alpha$-strongly iso-monotone on $(\manifold, (\cdot,\cdot))$.
    \item If for every $\Vector\neq \VectorB\in \manifold$ 
    \begin{equation}
        \| \partransport_{\VectorB\leftarrow \geodesic_{\Vector, \VectorB}(t)} \nabla^{\iso}_{\dot{\geodesic}_{\Vector, \VectorB}(t)}\tangentVector_{(\cdot)}\|_2 \leq L \|\log_{\VectorB}(\Vector)\|_2, \quad t\in [0,1],
        \label{eq:lem-iso-L-Lips-impl}
    \end{equation}
    then $\tangentVector$ is iso-$L$-Lipschitz on $(\manifold, (\cdot,\cdot))$.
    \end{enumerate}
    
\end{theorem}

We can now use the preceding result to examine more closely the relationship between Euclidean convex and $L$-smooth functions and iso-monotonicity and iso-$L$-Lipschitzness of their projected gradient fields beyond the one-dimensional case in \Cref{prop:equivalence-iso-monotone-smooth-R}. To that end, we begin by evaluating the iso-covariant derivative along a geodesic, as summarized in the following result, where we highlight that we are only able to do this on one-dimensional manifolds.

\begin{lemma}
\label{lem:second-deriv-along-iso-geo-vf}
    Let $\function:\Real^\dimInd\to \Real$ be a function that is twice continuously differentiable and let $(\manifold, (\cdot,\cdot))$ be a one-dimensional Riemannian manifold embedded in $\Real^\dimInd$ on which the iso-exponential map is both globally defined and a global diffeomorphism.

    Then, for every  $\Vector \neq \VectorB\in \manifold$
    \begin{align}
        \nabla^{\iso}_{\dot{\geodesic}_{\Vector, \VectorB}(t)} \mathbf{P}_{(\cdot)} & \nabla \function =
        \\
        &\frac{D_{\VectorC_t}^2 \function [\tangentVector_{\VectorC_t}, \tangentVector_{\VectorC_t}] - \bigl(\bigl(\mathbf{I} - \frac{\tangentVector_{\VectorC_t}}{\|\tangentVector_{\VectorC_t}\|_2} \otimes \frac{\tangentVector_{\VectorC_t}}{\|\tangentVector_{\VectorC_t}\|_2} \bigr) \nabla \function (\VectorC_t), \Gamma_{\VectorC_t} [\tangentVector_{\VectorC_t}, \tangentVector_{\VectorC_t}]\bigr)_2}{\|\tangentVector_{\VectorC_t}\|_2^2} \tangentVector_{\VectorC_t}, \nonumber
    \end{align}
    where $\tangentVector_{\VectorC_t} := \dot{\geodesic}_{\Vector,\VectorB}(t)$ and $ \VectorC_t = \geodesic_{\Vector,\VectorB}(t)$.
\end{lemma}

Next, using \Cref{thm:iso-monotone-impl} through \Cref{lem:second-deriv-along-iso-geo-vf}, it naturally follows that the projected gradients of Euclidean convex and $L$-smooth functions on one-dimensional manifolds are related to iso-monotone and iso-$L$-Lipschitz vector fields via a geometric correction; conversely, with an appropriate geometric correction, certain non-convex functions may have iso-monotone projected gradients. Relating this back to our prior result \Cref{prop:implications-on-manifold} we see once again that iso-monotonicity and iso-g-convexity are equivalent on one-dimensional manifolds and vice versa for iso-$L$-Lipschitzness and iso-$L$-smoothness.

\begin{proposition}[Implications on $\manifold$]
    \label{prop:implications-on-manifold-vf}
    Let $\function:\Real^\dimInd\to \Real$ be a function that is twice continuously differentiable and let $(\manifold, (\cdot,\cdot))$ be a one-dimensional Riemannian manifold embedded in $\Real^\dimInd$ on which the iso-exponential map is both globally defined and a global diffeomorphism. Furthermore, assume that 
    \begin{equation}
        \frac{D_{\Vector}^2 \function [\tangentVector_{\Vector}, \tangentVector_{\Vector}]}{\|\tangentVector_{\Vector}\|_2^2} \in [\beta, \eta],
    \end{equation}
    and
    \begin{equation}
    - \frac{( (\mathbf{I} - \frac{\tangentVector_{\Vector}}{\|\tangentVector_{\Vector}\|_2} \otimes \frac{\tangentVector_{\Vector}}{\|\tangentVector_{\Vector}\|_2})\nabla \function (\Vector), \Gamma_{\Vector} (\tangentVector_{\Vector}, \tangentVector_{\Vector}) )_2}{\|\tangentVector_{\Vector}\|_2^2} \in [\delta, \zeta],
    \label{eq:thm-geo-factor-vf}
    \end{equation}
    where $\Gamma_{\Vector}:\tangent_{\Vector} \manifold \times \tangent_{\Vector} \manifold \to \tangent_{\Vector} \Real^{\dimInd}$ is the Christoffel operator on $(\manifold, (\cdot, \cdot))$. Finally, assume that $\beta + \delta >0$. 
    
    Then, the following hold:
    \begin{enumerate}[label=(\roman*)]
        \item $\mathbf{P}_{(\cdot)}\nabla \function$ is $(\beta + \delta)$-strongly iso-monotone on $(\manifold, (\cdot,\cdot))$.
        \item $\mathbf{P}_{(\cdot)}\nabla \function$ is iso-$(\eta + \zeta)$-Lipschitz on $(\manifold, (\cdot,\cdot))$.
    \end{enumerate}
\end{proposition}

\begin{corollary}
\label{cor:equivalence-iso-monotone-iso-convex-1d}
     Let $\function:\Real^\dimInd\to \Real$ be a function that is twice continuously differentiable and let $(\manifold, (\cdot,\cdot))$ be a one-dimensional Riemannian manifold embedded in $\Real^\dimInd$ on which the iso-exponential map is both globally defined and a global diffeomorphism.

    Then:
    \begin{enumerate}[label=(\roman*)]
        \item $\function$ is ($\alpha$-strongly) iso-g-convex on $(\manifold, (\cdot,\cdot))$ if and only if $\mathbf{P}_{(\cdot)} \nabla \function$ is ($\alpha$-strongly) iso-monotone on $(\manifold, (\cdot,\cdot))$.
        \item $\function$ is iso-$L$-smooth on $(\manifold, (\cdot,\cdot))$ if and only if $\mathbf{P}_{(\cdot)} \nabla \function$ is iso-$L$-Lipschitz on $(\manifold, (\cdot,\cdot))$.
    \end{enumerate}
\end{corollary}

\begin{remark}
\label{rem:non-equiv-func-vf-defs}
    However, it is important to emphasize that the restriction to the one-dimensional setting is not merely technical. In higher dimensions, we do not expect differentiating through the projection operator to yield an analogous result. The main obstruction is that parallel transport does not, in general, preserve \(\ell^2\)-orthogonality. This issue is particularly pronounced for Riemannian structures inherited from a Euclidean pullback ambient structure, where one cannot expect parallel transport to preserve \(\ell^2\)-orthogonality. Consequently, there is a mismatch both in our ability to extend the argument in \Cref{prop:implications-on-manifold-vf} beyond one dimension and in the equivalence established in \Cref{cor:equivalence-iso-monotone-iso-convex-1d}. Thus, in general, the class of functions whose projected gradients are iso-monotone or iso-\(L\)-Lipschitz differs from the class of iso-\(g\)-convex and iso-\(L\)-smooth functions.

\end{remark}


Beyond vector fields that are projected gradients, there are others that fit our definitions for iso-monotone and iso-$L$-Lipschitz, which we showcase in the following two examples concerning iso-geodesic velocity fields and iso-logarithmic fields. However, we see a similar pattern, i.e., it is hard to establish any results beyond one-dimensional manifolds.

\begin{proposition}[Properties of iso-geodesic velocity fields]
\label{prop:iso-mono-geo}
   Let $(\manifold, (\cdot,\cdot))$ be a one-dimensional Riemannian manifold embedded in $\Real^\dimInd$ on which the iso-exponential map is both globally defined and a global diffeomorphism.
   
   Then, for every $\VectorD\neq \VectorE\in \manifold$ the vector field $\dot{\geodesic}^{\iso}_{\VectorD, \VectorE}$ is iso-monotone and iso-Lipschitz with $L=0$ on $(\geodesic_{\VectorD, \VectorE} ((0,1)), (\cdot,\cdot))$.
\end{proposition}

\begin{proposition}[properties of iso-logarithmic fields]
\label{prop:iso-mono-log-field}
    Let $(\manifold, (\cdot,\cdot))$ be a one-dimensional Riemannian manifold embedded in $\Real^\dimInd$ on which the iso-exponential map is both globally defined and a global diffeomorphism. 
    
    Then, for every $\VectorF\in \manifold$ the vector field 
    $
    -\log_{(\cdot)}^{\iso}(\VectorF) \in \vectorfield(\manifold)
    $ 
    is $1$-strongly iso-monotone and iso-$1$-Lipschitz.
\end{proposition}

\begin{remark}
\label{rem:iso-geo-iso-log-limitations}
    We note that the results in \Cref{prop:iso-mono-geo,prop:iso-mono-log-field} reminisce of monotonicity of geodesic velocity fields and 1-strong monotonicity and 1-Lipschitzness of logarithmic fields in the classical Riemannian sense, which renders them as sanity checks for \Cref{def:iso-monotone-vf,def:iso-lipschitz-vf}. That said, the discussion in \Cref{rem:non-equiv-func-vf-defs} extends here, i.e., we do not expect that these relations necessarily hold beyond one-dimensional manifolds for similar reasons.
\end{remark}

\subsection{Finding zeros of iso-monotone and iso-Lipschitz vector fields}
\label{sec:finding-mins-monotone}

Despite the limitations discussed in
\Cref{rem:non-equiv-func-vf-defs,rem:iso-geo-iso-log-limitations},
particularly the difficulty of verifying these conditions beyond
one-dimensional manifolds and their apparent incompatibility in higher
dimensions, we continue the analysis for completeness. At this point, we have all the ingredients needed to derive sufficient conditions for convergence of the iso-Riemannian descent scheme in \cref{eq:ird}, of which $\ell^2$-PG-IRD in \cref{eq:iso-first-order-scheme} is a direct special case.

Since our focus is on convergence of the iterates, we state the following result under assumptions of strong iso-monotonicity and iso-\(L\)-Lipschitzness. The result also characterizes conditions under which the scheme converges to a minimizer of a Euclidean convex and \(L\)-smooth function. This characterization is developed further in the subsequent corollary, which is notably more restrictive than \Cref{cor:iso-convex}. Together with the aforementioned limitations, this further supports the view that the function-based framework of \Cref{sec:function-convergence}, which yields \Cref{cor:iso-convex}, is likely to be of greater practical relevance.

\begin{theorem}[Convergence of iso-Riemannian descent]
    \label{thm:convergence-EGiRD}
    Let $\tangentVector \in \vectorfield(\manifold)$ be a $\alpha$-strongly iso-monotone and iso-$L$-Lipschitz vector field on a Riemannian manifold $(\manifold, (\cdot,\cdot))$ embedded in $\Real^\dimInd$ on which the iso-exponential map is both globally defined and a global diffeomorphism. Furthermore, assume that $\tangentVector$ has a singularity, i.e., there is a point $\bar{\Vector} \in \manifold$ that satisfies $\tangentVector_{\bar{\Vector}} = \mathbf{0}$. 
    
    Then, the sequence $\{\Vector^{(\sumIndC)}\}_{\sumIndC=0}^\infty \in \manifold$ defined by the algorithm
    \[
        \Vector^{(\sumIndC+1)}=\exp^{\iso}_{\Vector^{(\sumIndC)}}(-r \mathbf{P}_{\Vector^{(\sumIndC)}} \nabla f (\Vector^{(\sumIndC)})), \quad \Vector^{(0)} \in \manifold,
    \]
    where $0 < r < \frac{2 \alpha}{L^2}$ converges linearly in \(\ell^2\) to $\bar{\Vector} \in \manifold$.
\end{theorem}

\begin{corollary}
\label{cor:iso-monotone}
    Under the assumptions of \Cref{prop:implications-on-manifold-vf} -- notably assuming that $\beta + \delta >0$ and that the manifold $\manifold$ is one-dimensional --, suppose additionally that there exists a point $\bar{\Vector} \in \manifold$ satisfying $\mathbf{P}_{\bar{\Vector}} \nabla f(\bar{\Vector}) = \mathbf{0}.$
    
    Then, the sequence $\{\Vector^{(\sumIndC)}\}_{\sumIndC=0}^\infty \in \manifold$ defined by the algorithm
    \[
        \Vector^{(\sumIndC+1)}=\exp^{\iso}_{\Vector^{(\sumIndC)}}(-r \mathbf{P}_{\Vector^{(\sumIndC)}} \nabla f (\Vector^{(\sumIndC)})), \quad \Vector^{(0)} \in \manifold,
    \]
    where \(0 < r \leq \frac{2 (\beta + \delta)}{(\eta + \zeta)^2}\) converges linearly in \(\ell^2\) to $\bar{\Vector} \in \manifold$.    
\end{corollary}



\subsection{Implications for iso-barycentres and applications}
\label{sec:iso-bary}

Finally, on a more constructive note, and to highlight a promising direction for future work, we illustrate how the zero-finding formulation for vector fields can be used beyond the scope of the present analysis. In particular, we demonstrate how iso-Riemannian descent (IRD) enables the efficient computation of iso-Riemannian barycentres, and how these barycentres can improve downstream tasks such as clustering compared with na\"ive approaches based on Euclidean or standard Riemannian barycentres.

\subsubsection{Existence and uniqueness of iso-barycentres}

We begin by proposing a notion of the iso-Riemannian barycentre or \emph{iso-barycentre}. Recall that the Riemannian barycentre $\bar{\Vector} \in \manifold$, with respect to $(\manifold,(\cdot,\cdot))$, can -- at least locally -- be characterized either as the minimizer of the following optimization problem or equivalently as the zero of its Riemannian gradient field:
\begin{equation}
    \bar{\Vector} := \operatorname{argmin}_{\Vector\in \manifold} \frac{1}{2\dataPointNum} \sum_{\sumIndA=1}^\dataPointNum \distance_{\manifold}(\Vector, \Vector^\sumIndA)^2, \quad \Leftrightarrow \quad \frac{1}{\dataPointNum} \sum_{\sumIndA=1}^\dataPointNum \log_{\bar{\Vector}} (\Vector^\sumIndA) = \mathbf{0}.
\end{equation}

When replacing the above manifold mappings by iso-versions, the above equivalence can -- even locally -- not be expected to hold anymore. In this work, we will use the following definition based on the vector field characterization.

\begin{definition}[Iso-Riemannian barycentre] 
\label{def:iso-bary}
    A point $\bar{\Vector}\in \manifold$ is called an \emph{iso-Riemannian barycentre} or \emph{iso-barycentre} of a data set $\{\Vector^\sumIndA\}_{\sumIndA=1}^\dataPointNum \subset \manifold$ on a Riemannian manifold $(\manifold, (\cdot,\cdot))$ embedded in $\Real^\dimInd$ on which the iso-exponential map is both globally defined and a global diffeomorphism, if it satisfies 
    \begin{equation}
        \frac{1}{\dataPointNum} \sum_{\sumIndA=1}^\dataPointNum \log^{\iso}_{\bar{\Vector}} (\Vector^\sumIndA) = \mathbf{0}.
        \label{eq:iso-bary-problem}
    \end{equation}
\end{definition}

The motivation for \Cref{def:iso-bary} is illustrated by the following two examples: the iso-barycentre of two points lies at the midpoint of the iso-geodesic connecting them, ensuring that the iso-distances to both endpoints -- measured as the $\ell^2$-arclengths along the curve -- are equal, as one would naturally expect, and the iso-barycentre of a data set on $\Real$ is the Euclidean mean. 

\begin{proposition}[Iso-geodesic midpoint]
\label{prop:two-point-iso-bary}
    Let $(\manifold, (\cdot,\cdot))$ be a Riemannian manifold embedded in $\Real^\dimInd$ on which the iso-exponential map is both globally defined and a global diffeomorphism.
    
    Then, the iso-barycentre $\bar{\Vector}$ on $(\manifold, (\cdot,\cdot))$ of any two points $\Vector^1, \Vector^2\in \manifold$ on a  satisfies $\bar{\Vector} = \geodesic^{\iso}_{\Vector^1, \Vector^2} (\frac{1}{2})$.
\end{proposition}

\begin{proposition}[Iso-barycentres on $\Real$]
\label{prop:iso-bary-R}
    Consider any Riemannian structure $(\Real,(\cdot, \cdot))$.

    Then, the iso-barycentre $\bar{x}$ on $(\Real,(\cdot, \cdot))$  of any data set $\{x^\sumIndA\}_{\sumIndA}^\dataPointNum\in \Real$ satisfies $\bar{x} = \frac{1}{\dataPointNum}\sum_{\sumIndA=1}^\dataPointNum x^\sumIndA$. 
\end{proposition}

In more general settings, \Cref{def:iso-bary} raises several preliminary questions even before considering the application of iso-Riemannian descent to compute an iso-barycentre. In particular, both existence and uniqueness remain unclear. While the motivating example in \Cref{prop:two-point-iso-bary} guarantees the existence of at least one iso-barycentre, without additional assumptions there is no guarantee of uniqueness.

Regarding existence, as in the case of the Riemannian barycentre, it depends jointly on the data set and the underlying Riemannian structure.

\begin{theorem}[Existence of iso-barycentres]
\label{thm:existence-iso-barycentre}
    Let $(\manifold, (\cdot,\cdot))$ be a Riemannian manifold embedded in $\Real^\dimInd$ on which the iso-exponential map is both globally defined and a global diffeomorphism. Furthermore, let $\{\Vector^\sumIndA\}_{\sumIndA=1}^\dataPointNum \subset \manifold$ be a set of points and assume that there is a geodesically convex submanifold $\manifold'\subset\manifold$ containing these points. 
    
    Then, there exists at least one $\bar{\Vector} \in \manifold$ that solves (\ref{eq:iso-bary-problem}).
\end{theorem}

A natural route to uniqueness is to demonstrate $\alpha$-strong iso-monotonicity on $\manifold$. In general, however, this is difficult to quantify, and the same holds for Lipschitz continuity -- which is required to show convergence of iso-Riemannian descent. Nevertheless, in the one-dimensional case, explicit iso-monotonicity and iso-Lipschitz constants of the iso-barycentre field can in fact be obtained.

\begin{theorem}[Properties of 1D iso-barycentre fields]
\label{thm:iso-mono-bary-field-1d}
     Let $(\manifold, (\cdot,\cdot))$ be a Riemannian manifold embedded in $\Real^\dimInd$ on which the iso-exponential map is both globally defined and a global diffeomorphism. 
     
     Then, for any set of points $\{\Vector^\sumIndA\}_{\sumIndA=1}^\dataPointNum \subset \manifold$, the vector field 
    $
    -\frac{1}{\dataPointNum} \sum_{\sumIndA=1}^\dataPointNum 
    \log_{(\cdot)}^{\iso}(\Vector^\sumIndA)  \in \vectorfield(\manifold)
    $
    is $1$-strongly iso-monotone and iso-$1$-Lipschitz.
\end{theorem}

\begin{corollary}
\label{thm:iso-mono-bary-field-1d-uniqueness}
    Under the assumptions in \Cref{thm:iso-mono-bary-field-1d}, there is a unique $\bar{\Vector} \in \manifold$ that solves (\ref{eq:iso-bary-problem}).
\end{corollary}

\subsubsection{Computing iso-barycentres}

\begin{figure}[h!]
    \centering
    \setlength{\tabcolsep}{2pt} 
    \renewcommand{\arraystretch}{0} 
    \begin{tabular}{@{}c ccc@{}}
        \makebox[0pt][r]{\raisebox{35pt}{\rotatebox[origin=c]{90}{\small\shortstack{river\\data set}}}} &
        \includegraphics[width=0.28\textwidth]{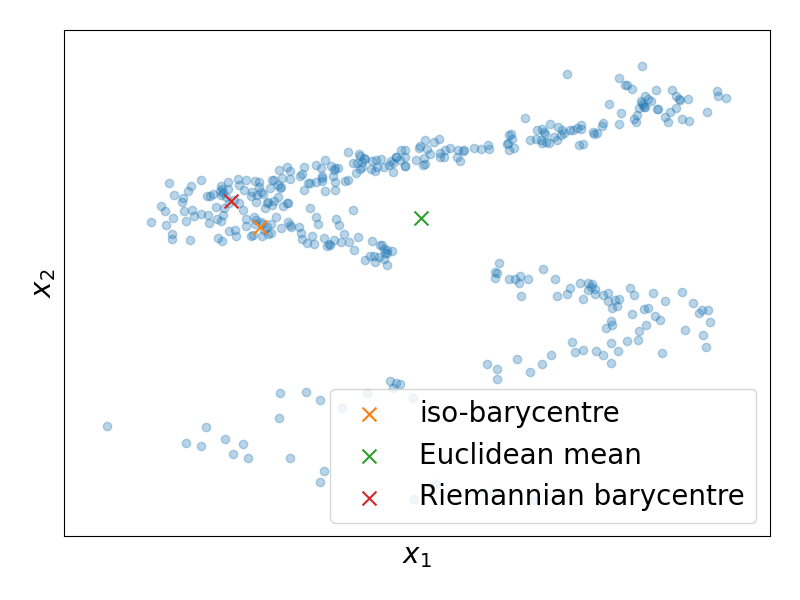} &
        \includegraphics[width=0.28\textwidth]{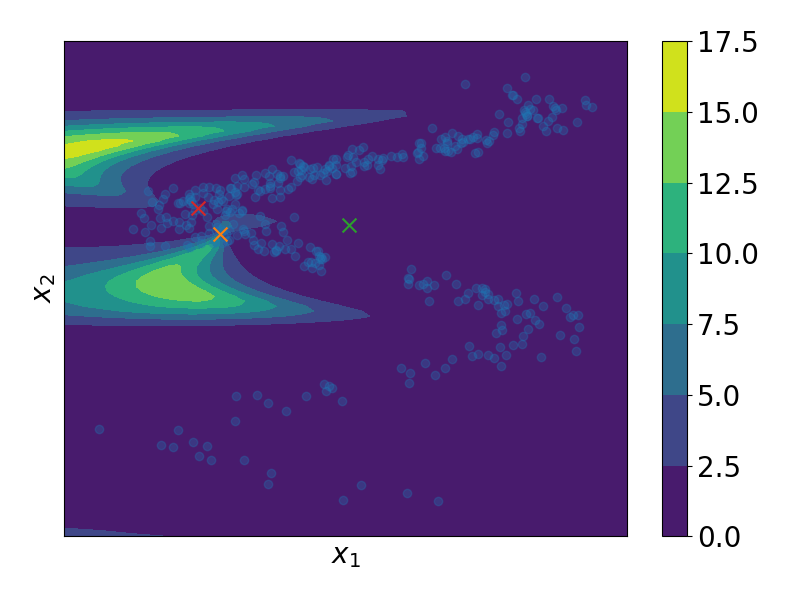} &
        \includegraphics[width=0.28\textwidth]{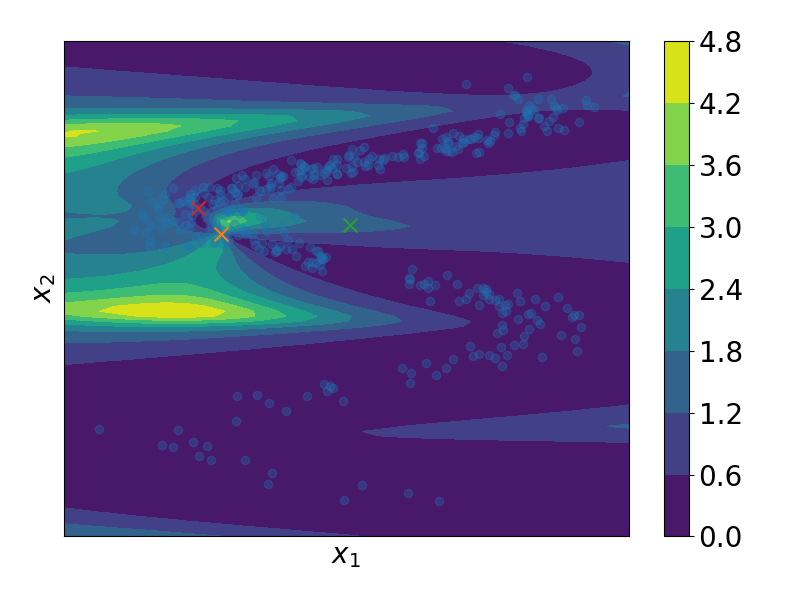} \\
        \makebox[0pt][r]{\raisebox{35pt}{\rotatebox[origin=c]{90}{\small\shortstack{spiral\\data set}}}} &
        \includegraphics[width=0.28\textwidth]{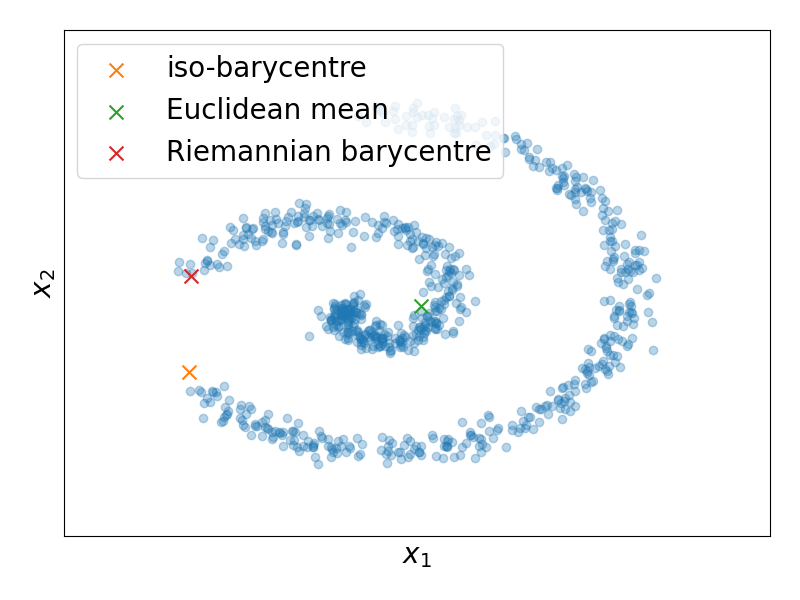} &
        \includegraphics[width=0.28\textwidth]{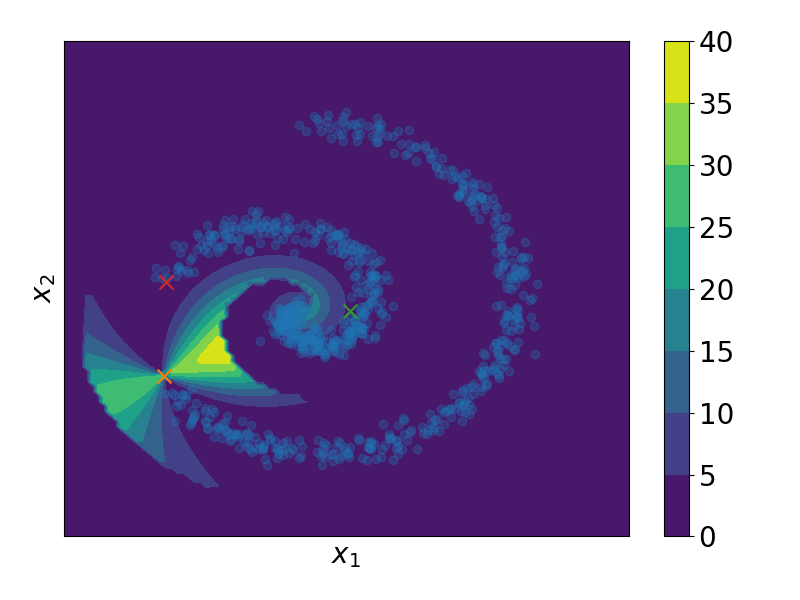} &
        \includegraphics[width=0.28\textwidth]{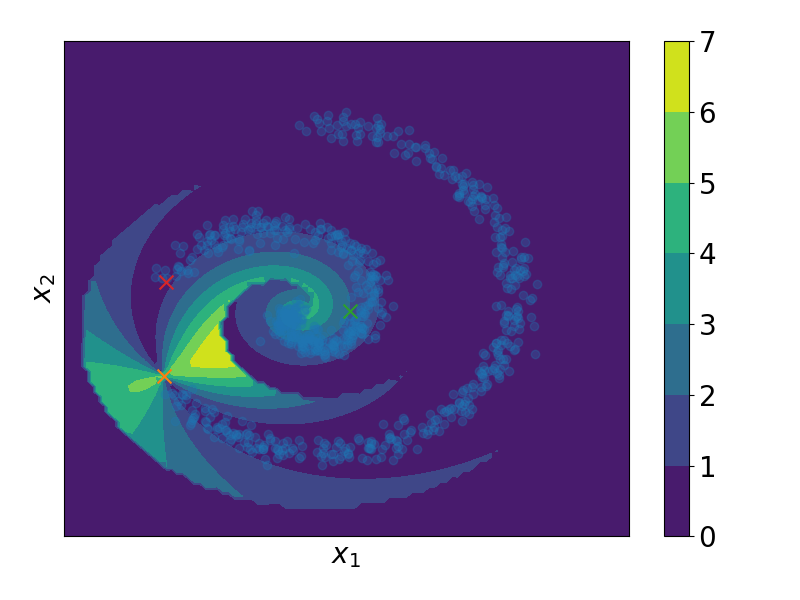} \\
        & \makebox[0.28\textwidth]{\small Barycentres}
        & \makebox[0.28\textwidth]{\small Iso-monotonicity ratio}
        & \makebox[0.28\textwidth]{\small Iso-Lipschitz ratio}
    \end{tabular}
    \caption{Barycentre results on the synthetic river and spiral data sets under modeled pullback structures. Both the Riemannian and iso-barycentre (computed with \Cref{alg:iso-bary-line-search}) yield more interpretable data representations -- i.e., lying within the data set -- compared to the Euclidean mean. Moreover, based on the observed iso-monotonicity and iso-Lipschitz ratios, we generally expect the parameter $\alpha$ to be much smaller than 1, while $L$ tends to be much larger.}
    \label{fig:synthetic-results-barys}
\end{figure}

\Cref{thm:iso-mono-bary-field-1d} aligns with the behavior of the classical Riemannian barycentre field, exhibiting the same strong monotonicity and Lipschitz constants along a one-dimensional manifold. The situation in higher dimensions, however, is less clear. For example, unlike the 1D case, non-zero curvature of the underlying Riemannian structure is expected to alter the behavior, much as in the classical Riemannian barycentre setting. In addition to that, even in the zero-curvature case -- for instance, under Euclidean pullback structures we focus on in this work -- the iso-monotonicity and iso-Lipschitz constants do not necessarily remain close to 1 in practice. This is illustrated in the $\Real^2$-valued river and spiral data sets in \Cref{fig:synthetic-results-barys}, where under modeled pullback structures on $\manifold := \Real^2$ (see \Cref{sec:app-iso-bary} for details) the minimum of the \emph{iso-monotonicity ratio} and the maximum of the \emph{iso-Lipschitz ratio} from the iso-barycentre\footnote{That is, normalized strong iso-monotonicity from \Cref{def:iso-monotone-vf} and iso-Lipschitzness from \Cref{def:iso-lipschitz-vf} where we use that the iso-barycentre field vanishes at the iso-barycentre.}
\begin{equation}
    \Vector\mapsto \frac{(\frac{1}{\dataPointNum} \sum_{\sumIndA=1}^\dataPointNum \log^{\iso}_{\Vector} (\Vector^\sumIndA), \partransport^{\iso}_{\Vector\leftarrow \bar{\Vector}} \log^{\iso}_{\bar{\Vector}} (\Vector))_2}{\distance^{\iso}_{\manifold} (\bar{\Vector}, \Vector)^2} \quad \text{and} \quad \Vector\mapsto \frac{\| \frac{1}{\dataPointNum} \sum_{\sumIndA=1}^\dataPointNum \log^{\iso}_{\Vector} (\Vector^\sumIndA)\|_2}{\distance^{\iso}_{\manifold}(\bar{\Vector}, \Vector)},
\end{equation}
give upper and lower bounds respectively to the constants $\alpha$ and $L$. In particular, neither $\alpha$ nor $L$ is close to 1\footnote{It is also good to note that despite the positive iso-monotonicity ratios for both data sets and pullback structures, it is unclear whether this can be expected in general or whether there exist cases where it becomes negative -- indicating that the iso-barycentre field would be non-$\alpha$-strongly iso-monotone in general.}.

\begin{remark}
    It is worth emphasizing the distinction in setting. Rather than computing barycentres on a strictly lower-dimensional manifold, we consider barycentres on the full ambient space equipped with a nontrivial Riemannian structure, which is likewise covered by the proposed framework.
\end{remark}

Building on these observations, we introduce \Cref{alg:iso-bary-line-search} for computing iso-barycentres in the general setting. The algorithm combines iso-Riemannian descent with a line-search procedure and is initialized at the Riemannian barycentre. This initialization ensures that the iterates begin in the geodesically convex submanifold \(\manifold'\) containing the data set, as motivated by the existence result for iso-barycentres established in \Cref{thm:existence-iso-barycentre}.

\begin{figure}[h!]
    \centering
    \setlength{\tabcolsep}{2pt}
    \renewcommand{\arraystretch}{0}
    \begin{tabular}{ccc}
        \includegraphics[width=0.29\textwidth]{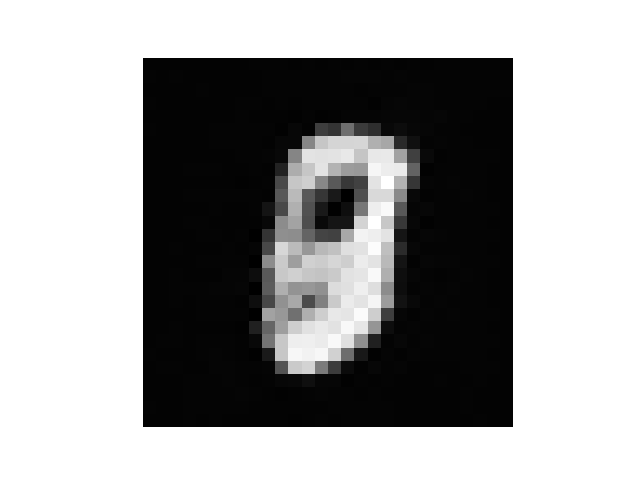} &
        \includegraphics[width=0.29\textwidth]{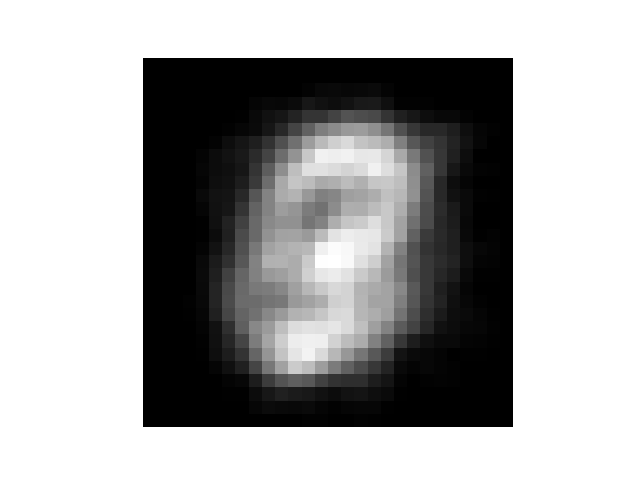} &
        \includegraphics[width=0.29\textwidth]{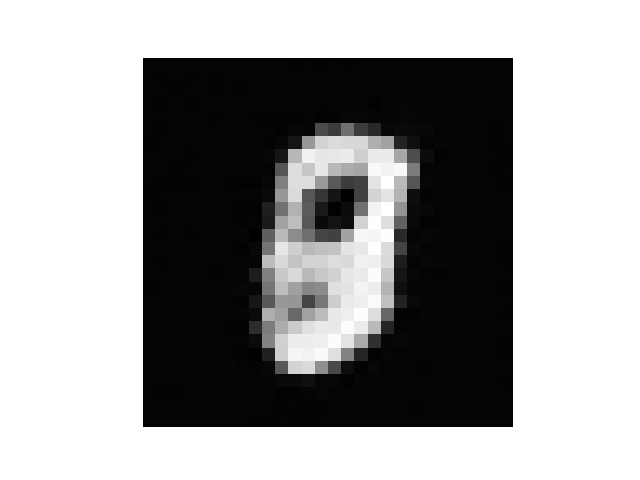} \\
        \makebox[0.29\textwidth]{\small Iso-barycentre} &
        \makebox[0.29\textwidth]{\small Euclidean mean} &
        \makebox[0.29\textwidth]{\small Riemannian barycentre}
    \end{tabular}
    \caption{Barycentre results for MNIST under different geometries. Both the Riemannian and iso-barycentre (\Cref{alg:iso-bary-line-search}) give very similar and interpretable data representations -- looking somewhat like an 8 -- whereas the Euclidean mean just yields a blurry mix. Notably, the iso-barycentre and the Riemannian barycentre are visually very similar.}
    \label{fig:mnist-results-barys}
\end{figure}

Empirically, \Cref{alg:iso-bary-line-search} enables stable computation of iso-barycentres both on the modeled pullback structures underlying the synthetic data sets in \Cref{fig:synthetic-results-barys} and on real-world data. In particular, \Cref{fig:mnist-results-barys} demonstrates its application to MNIST on a learned \(20\)-dimensional manifold equipped with the off-the-shelf pretrained pullback metric from \cite{diepeveen2025manifold}.

\begin{figure}[h!]
    \centering
    \setlength{\tabcolsep}{2pt} 
    \renewcommand{\arraystretch}{0} 
    \begin{tabular}{@{}c ccc@{}}
        \makebox[0pt][r]{\raisebox{35pt}{\rotatebox[origin=c]{90}{\small\shortstack{river\\data set}}}} &
        \includegraphics[width=0.28\textwidth]{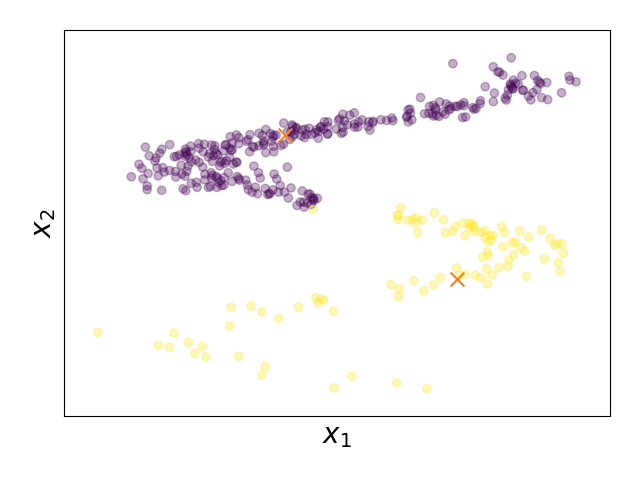} &
        \includegraphics[width=0.28\textwidth]{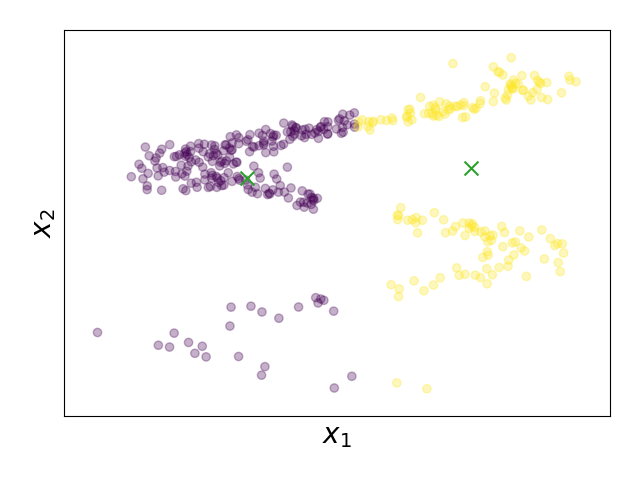} &
        \includegraphics[width=0.28\textwidth]{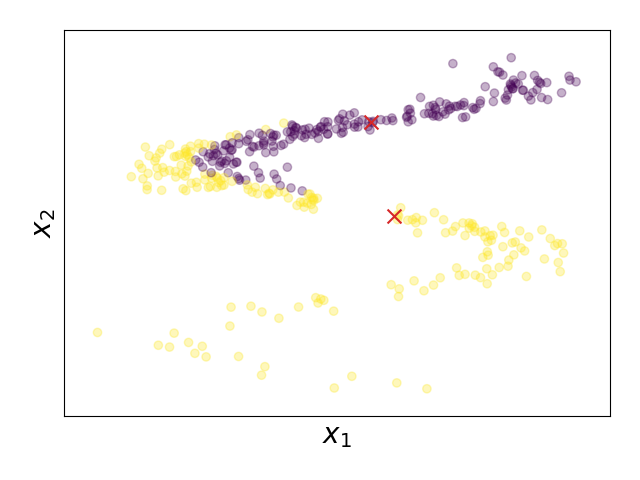} \\
        \makebox[0pt][r]{\raisebox{35pt}{\rotatebox[origin=c]{90}{\small\shortstack{spiral\\data set}}}} &
        \includegraphics[width=0.28\textwidth]{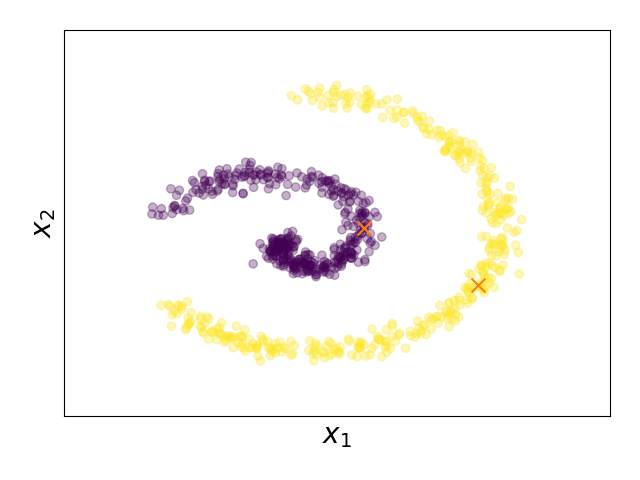} &
        \includegraphics[width=0.28\textwidth]{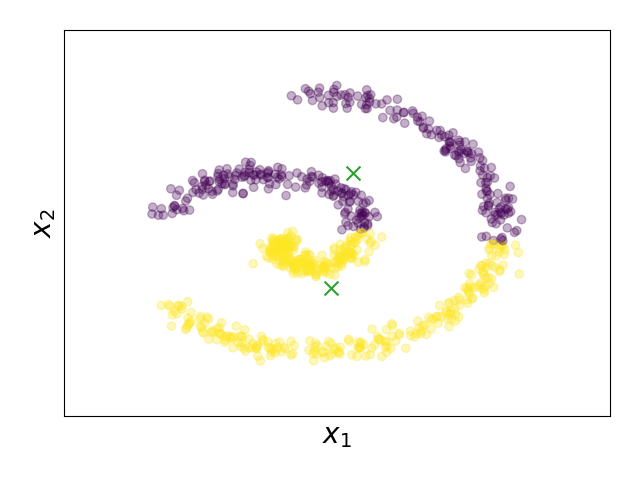} &
        \includegraphics[width=0.28\textwidth]{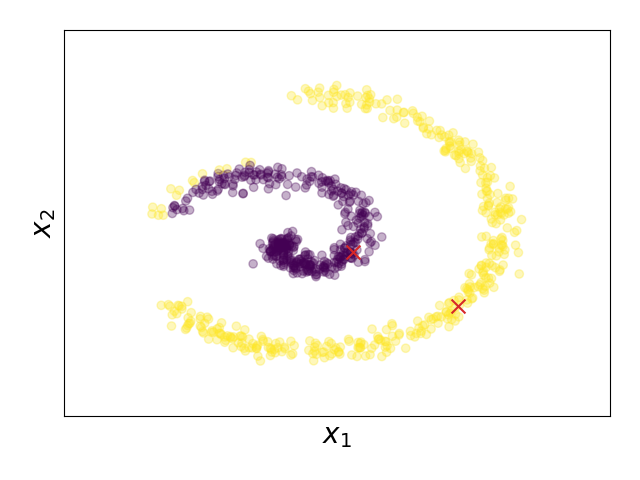} \\
        & \makebox[0.28\textwidth]{\small Iso-K-means}
        & \makebox[0.28\textwidth]{\small Euclidean K-means}
        & \makebox[0.28\textwidth]{\small Riemannian K-means}
    \end{tabular}
    \caption{Clustering results on the synthetic river and spiral data sets under modeled pullback structures. Both the Riemannian and iso-barycentre-based methods yield more interpretable data representations compared to the Euclidean mean. However, only iso-K-means clustering (\Cref{alg:iso-riem-kmeans}) also gives meaningful clusters. That is, in contrast to Riemannian K-means, which gives odd clusters due to the fact that the Riemannian distance does not encode proximity in an $\ell^2$-sense.}
    \label{fig:synthetic-results-k-means}
\end{figure}

\begin{remark}
\label{rem:little-difference}
Although the Riemannian and iso-barycentres in \Cref{fig:synthetic-results-barys,fig:mnist-results-barys} lie in similar locations, making their relative representational quality difficult to assess directly, the distinction becomes clearer in downstream tasks. In particular, \Cref{fig:synthetic-results-k-means} shows that iso-barycentres yield superior clustering performance; see \Cref{sec:iso-k-means} for details.
\end{remark}

\begin{remark}
\label{rem:num-error-pile-up}
    For all three data sets the convergence to the iso-barycentre using \Cref{alg:iso-bary-line-search} is linear, as expected from the one-dimensional case in \Cref{thm:convergence-EGiRD}.
    However, for the high-dimensional MNIST data set we notice that the combination of single precision of the data structures in \texttt{PyTorch} -- in which the pretrained Riemannian structure by \cite{diepeveen2025manifold} was implemented -- and the numerical approximations for the iso-manifold mappings -- as proposed in \cite{diepeveen2025manifold} -- affect the convergence behavior. In particular, numerical errors eventually become significant enough to disrupt the optimization process, so that the search direction may no longer be a descent direction. As a result, the line search fails and the algorithm stalls when too high an accuracy is sought. That said, solving these numerical issues is beyond the scope of this work.
\end{remark}

\begin{algorithm}[h!]
\caption{Iso-Riemannian barycentre through iso-Riemannian descent (IRD) with line search}
\label{alg:iso-bary-line-search}
\begin{algorithmic}
\REQUIRE  Data points $\{\Vector^{\sumIndA}\}_{\sumIndA=1}^\dataPointNum \subset \manifold$, initial step size $r_0 > 0$, reduction parameter $c \in (0,1)$, iteration counter $k = 0$
\STATE Initial point $\Vector^{(0)} \in \operatorname{argmin}_{\Vector\in \manifold}  \frac{1}{2\dataPointNum} \sum_{\sumIndA=1}^\dataPointNum \distance_{\manifold}(\Vector, \Vector^\sumIndA)^2$ \hfill (Riemannian barycentre)
\WHILE{not converged}
    \STATE $\tangentVector_{\Vector^{(k)}} \gets -\frac{1}{\dataPointNum} \sum_{\sumIndA=1}^\dataPointNum 
    \log_{\Vector^{(k)}}^{\iso}(\Vector^\sumIndA)$
    \STATE $r \gets r_0$
    \REPEAT
        \STATE $\Vector_{\text{trial}} \gets \exp^{\iso}_{\Vector^{(k)}} (- r \tangentVector_{\Vector^{(k)}})$
        \STATE $\tangentVector_{\Vector_{\text{trial}}} \gets -\frac{1}{\dataPointNum} \sum_{\sumIndA=1}^\dataPointNum 
    \log_{\Vector_{\text{trial}}}^{\iso}(\Vector^\sumIndA)$
        \IF{ $\|\tangentVector_{\Vector_{\text{trial}}}\|_2 > \|\tangentVector_{\Vector^{(k)}}\|_2$ }
            \STATE $r \gets c r$ 
        \ENDIF
    \UNTIL{ $\|\tangentVector_{\Vector_{\text{trial}}}\|_2 < \|\tangentVector_{\Vector^{(k)}}\|_2$  }
    \STATE $\Vector^{(k+1)} \gets \Vector_{\text{trial}}$
    \STATE $k \gets k+1$
\ENDWHILE
\end{algorithmic}
\end{algorithm}

\section{Conclusions}
\label{sec:conclusions}

This work develops a first-order optimization framework for iso-Riemannian geometry, motivated in particular by the setting of optimization on Riemannian manifolds learned from data. In such settings, a direct application of classical Riemannian optimization theory can be inadequate: objectives with favorable Euclidean convexity or \(L\)-smoothness properties may fail to exhibit corresponding geodesic properties with respect to the learned Riemannian geometry, while the Riemannian gradient may produce a poorly conditioned or otherwise unsuitable search direction. By instead working with the iso-connection and projected Euclidean gradient directions, we obtain optimization schemes that are adapted both to the pullback structure underlying the learned geometry and to the structure of the objective.

We considered two routes to convergence. From a function-based perspective, we introduced iso-g-convexity and iso-\(L\)-smoothness, related strong iso-g-convexity to a PL-type condition, and established general convergence guarantees for the proposed first-order scheme. This framework addresses the conditioning and search-direction issues arising in the Riemannian formulation. From a vector-field perspective, we introduced iso-monotonicity and iso-Lipschitzness and established analogous results in the one-dimensional setting. Although this perspective currently yields a more limited theory, it applies naturally to problems posed as finding zeros of vector fields, such as the iso-Riemannian barycentre problem, which need not admit a direct function-based formulation. 

A central conclusion is that the function- and vector-field-based routes, which are closely connected in the classical Levi-Civita setting, give rise to genuinely different convergence theories in the iso-Riemannian setting. In particular, the classical implications between geodesic convexity and monotonicity of the Riemannian gradient, as well as between smoothness and Lipschitz continuity of the Riemannian gradient, do not generally carry over to their iso-Riemannian analogues. Consequently, extending a classical Riemannian optimization result to the iso-Riemannian setting requires care: the appropriate generalization depends on whether the problem is most naturally described through an objective function or through a vector field. Our results suggest that the function-based perspective is often the more tractable route for optimization in higher dimensions, whereas the vector-field perspective remains essential for applications that are intrinsically formulated as root-finding problems.

Several directions for future work remain. First, it would be valuable to develop a higher-dimensional convergence theory for iso-monotone and iso-Lipschitz vector fields, or to identify weaker conditions that retain the flexibility needed for vector-field-based applications. Second, the present framework is formulated using the iso-exponential map. Replacing this map with an appropriate notion of iso-Riemannian retraction -- which would be advantageous in terms of runtime -- would require both a definition of retractions compatible with iso-Riemannian geometry and a corresponding extension of the convergence theory. Such an extension lies beyond the scope of this work, but we expect that the underlying principles developed here should extend beyond the iso-exponential-map setting.

Finally, several algorithmic extensions are of interest. Higher-order methods may alleviate some of the limitations associated with step-size selection and improve practical performance on challenging geometries. Developing nonsmooth iso-Riemannian optimization methods is another important direction, including methods for iso-medians, viewed as analogues of Riemannian medians, and inverse problems with nonsmooth objectives. Addressing these questions will likely require iso-Riemannian counterparts of proximal methods, duality theory, or related nonsmooth optimization tools.

\subsubsection*{Acknowledgments} 

MW was partially supported by an Alfred P. Sloan Research Fellowship in
Mathematics, NSF Awards CBET-2112085 and DMS-2406905, and the Army Research Office under Grant No. W911NF-26-1-A060.
The views and conclusions contained in this document are those of the authors and should not be interpreted as representing the official policies, either expressed or implied, of the Army Research Office or the U.S. Government. The U.S. Government is authorized to reproduce and distribute reprints for Government purposes notwithstanding any copyright notation herein.

Some of the numerical computations in this paper were run on the FASRC Cannon cluster supported by the FAS Division of Science Research Computing Group at Harvard.

\bibliographystyle{plain}
\bibliography{bibliography}

\appendix
\newpage
\section{Related work}
\label{sec:related-work}
To the best of our knowledge, no prior work has investigated optimization in the context of iso‑Riemannian geometry. That said, the intrinsic perspective underlying iso‑Riemannian geometry closely echoes the extensive literature on optimization over abstract manifolds, where the geometry of the space is taken as a fundamental structural constraint. This suggests that a unified Riemannian framework, rather than the predominantly chart‑based approaches common in the manifold learning literature, has the potential to alleviate many of the practical and theoretical challenges that arise when optimizing under a data-driven Riemannian geometry.

In the remainder of this section we review two related lines of work. First, we recall approaches for optimization on abstract manifolds, which operate purely in the intrinsic setting. Second, we discuss optimization methods tailored to learned data manifolds. 

\paragraph{Optimization on abstract manifolds}
The original motivation for optimization on manifolds was to exploit lower‑dimensional structure in order to make optimization both more efficient and more well‑posed. This motivation remains closely aligned with our goals, particularly in the context of solving inverse problems where learned data manifolds capture intrinsic low‑dimensional structure. For example, directions transverse to such manifolds often come with very large Lipschitz constants when using learned regularizers, which in turn forces conventional optimization methods to take prohibitively small step sizes. By rephrasing optimization  with respect to the intrinsic geometry of a manifold, one can avoid these difficulties.

Adopting a Riemannian structure has distinct advantages. Compared to early projection‑based methods \cite{luenberger1972gradient}, intrinsic Riemannian optimization no longer requires iterated projections back onto the manifold. Compared to chart‑based approaches, the Riemannian setting enables global reasoning without dependence on local coordinate neighborhoods, thereby sidestepping geometric distortions. Crucially, this makes it possible to analyze optimization in terms of geodesic convexity, a natural analogue of convexity in Euclidean spaces but one that is valid at a global scale.

Over the past three decades, these insights have led to a rapidly expanding literature that systematically generalizes classical optimization algorithms to the Riemannian setting \cite{absil2008optimization,boumal2023introduction}. Early work extended gradient descent, Newton’s method, and conjugate gradient techniques to manifolds (see, e.g., \cite{smith1994optimization,udriste1994convex}). This line of research was followed by generalizations of quasi‑Newton methods \cite{huang2015broyden,huang2018riemannian}, and second‑order methods such as trust region frameworks \cite{absil2007trust} and adaptive regularization with cubics \cite{agarwal2021adaptive}. Beyond smooth optimization, substantial progress has been made in extending non‑smooth and constrained methods to the manifold setting. For non‑smooth optimization, this includes subgradient methods \cite{bergmann2024riemannian,ferreira1998subgradient,hoseini2023proximal}, proximal mappings \cite{ferreira2002proximal}, and splitting schemes \cite{adachi2022riemannian,bacak2014computing,bergmann2024difference,bergmann2021fenchel,bergmann2016parallel,diepeveen2021inexact,souza2015proximal}. For constrained optimization, Riemannian analogues of the augmented Lagrangian method \cite{liu2020simple}, interior point methods \cite{lai2024riemannian}, and the Frank–Wolfe algorithm \cite{weber2023riemannian,weber2021projection} have been developed. A particularly notable recent direction is the development of landing methods \cite{ablin2022fast}, originally proposed for optimization over orthogonal matrices. These methods combine elements of projection and intrinsic strategies, providing new algorithmic flexibility. Subsequent extensions have broadened their scope to stochastic and variance‑reduced settings \cite{ablin2024infeasible} as well as to the Stiefel manifold \cite{goyens2025geometric,vary2024optimization}.

Despite these advances, the focus of the Riemannian optimization community has primarily been on well‑understood matrix manifolds, where closed‑form expressions for manifold mappings -- or first-order approximations thereof -- are readily available, and where the objective functions exhibits convexity in the Riemannian sense. This leaves open important questions regarding settings where the objective function carries Euclidean structure, e.g., strong convexity or Lipschitz continuity of its gradients in the ambient space, while the optimization proceeds along a learned Riemannian geometry. In principle, local convergence guarantees from the abstract manifold literature carry over, but the precise interaction between Euclidean notions of regularity and Riemannian intrinsic geometry is much less clear. It is exactly in this gap that our work is situated.

\paragraph{Optimization on learned data manifolds}

Attempts to bring ideas from Riemannian optimization into the setting of learned data manifolds have largely proceeded by constraining optimization to linear Euclidean subspaces or to chart‑based parameterizations of non-linear manifolds.

One prominent direction has been the study of functions whose variation is effectively confined to a low‑dimensional linear subspace. In such cases, methods have been proposed that first extract this subspace before performing optimization \cite{cartis2024learning,cosson2023low}. Because the optimization remains linear, structural properties inherited from the ambient problem -- such as strong convexity and gradient Lipschitz continuity -- transfer in a relatively straightforward manner, allowing convergence guarantees from unconstrained Euclidean optimization to carry over. However, this setting differs from the one we consider here. Our focus is on situations where the geometry is dictated by the data itself rather than by the function of interest, and where the relevant structure cannot be assumed to be linear.

A second strand of research considers solving inverse problems over non-linear manifolds but does so through global chart‑based representations, typically using the latent space of generative models. In this paradigm, inverse problems are formulated as optimization problems constrained to the latent space, with the neural network serving as the chart that maps latent variables to data space. Seminal work \cite{hand2018global} established that, under certain conditions for random networks, the optimization landscape for linear inverse problems with Euclidean fidelity terms is benign, and related results extend to specific non-linear inverse problems such as phase retrieval \cite{hand2018phase}. For trained networks, however, this favorable landscape disappears: in general, inverting neural networks is NP‑hard \cite{lei2019inverting}, spurring development of specialized solvers for inversion with generative priors \cite{daras2021intermediate}. Chart‑based optimization inherits additional significant drawbacks. Besides convexity, properties such as Lipschitz continuity of gradients are not preserved under neural charting either -- especially when activations such as ReLU can inflate Lipschitz constants drastically. Consequently, even basic gradient methods require heavy tuning and often lack generality, making this line of work unsatisfactory for the broader scope of optimization problems we are interested in.

A more geometrically motivated approach begins with local estimates of the manifold structure itself. Methods using local polynomial fits \cite{shustin2022manifold} (based on \cite{sober2020manifold}) or quadratic approximations from sampled data \cite{robinett2025manifold} (based on \cite{sritharan2021computing}) provide explicit local charts and use the Riemannian geometry inherited from the Euclidean metric on the ambient space to guide optimization. Compared to latent‑space optimization, this line of work does account for intrinsic geometry, but it still faces critical difficulties. Chief among these is the dependence of smoothness constants on the choice of charts. Convexity and Lipschitz gradients in the Euclidean ambient space do not translate robustly into convexity in the learned charts and lead to distorted Lipschitz constants, once again forcing step sizes to be unnecessarily small. In addition, this dependency on the chart undermines the possibility of a general and interpretable optimization framework as one needs to change charts during optimization.

A major obstacle to practical optimization on data manifolds has been the lack of reliable global manifold representations equipped with explicit manifold mappings. Recent work has begun to address this issue. In particular, the pullback constructions of \cite{diepeveen2024pulling,diepeveen2025scorebased} yield practical, scalable global manifold models with closed-form manifold mappings, while the iso-Riemannian extension of \cite{diepeveen2025manifold} further accounts for distortions induced by non-constant Euclidean speed. These developments provide the geometric ingredients that were previously unavailable for optimization beyond chart-based or purely local methods; it is exactly in this gap that our work is situated.

\section{Supplementary materials to \Cref{sec:function-convergence}}
\label{app:supplement-function-convergence}

\subsection{Proofs of the results in \Cref{sec:iso-g-convex,sec:finding-mins-g-convex}}
\label{app:supplement-function-convergence-proofs}

\paragraph{Proofs of the results in \Cref{sec:iso-g-convex}}

\begin{proof}[Proof of \Cref{prop:equivalence-iso-convex-smooth-R}]
    First, we note that on $\Real$ all iso-manifold mappings reduce to simpler expressions, i.e., for $x,y \in \Real$ and $\xi_x \in \tangent_x \Real$ we have
    \begin{align}
        \geodesic^{\iso}_{x,y}(t) &= (1 - t) x + t y,\\
        \log^{\iso}_{x} (y) &= y - x,\\
        \exp^{\iso}_{x} (\xi_x) &= x + \xi_x,\\
        \distance^{\iso}_{\Real}(x,y) & = |x - y|,\\
        \partransport^{\iso}_{y \leftarrow x}\xi_x &= \xi_x.
    \end{align}
    So for showing (i) we note that
    \begin{equation}
         \function(y) \geq \function(x) + (\nabla \function(x), \log^{\iso}_{x}(y) )_2 + \frac{\alpha}{2} \distance^{\iso}_{\Real} (x, y)^2 
    \end{equation}
    is equivalent to 
    \begin{equation}
        \function(y) \geq \function(x) +  \function'(x) (y-x) + \frac{\alpha}{2} (x-y)^2.
    \end{equation}
    In other words, ($\alpha$-strong) iso-g-convexity of $f$ implies ($\alpha$-strong) convexity and vice versa, which proves the first claim.

    Similarly, for showing (ii) we also have that
    \begin{equation}
         \function(y) \leq \function(x) + (\nabla \function(x), \log^{\iso}_{x}(y) )_2 + \frac{L}{2} \distance^{\iso}_{\Real} (x, y)^2 
    \end{equation}
    is equivalent to 
    \begin{equation}
        \function(y) \leq \function(x) +  \function'(x) (y-x) + \frac{L}{2} (x-y)^2.
    \end{equation}
    So, iso-$L$-smoothness of $f$ implies $L$-smoothness of $f$ and vice versa, which proves the second claim.
\end{proof}

\begin{proof}[Proof of \Cref{thm:iso-convex-smooth-impl}]
    First, we define $\phi:[0,1]\to \Real$ as
    \begin{equation}
        \phi (t) := \function (\geodesic^{\iso}_{\Vector,\VectorB}(t)),
    \end{equation}
    and consider its Taylor approximation around $t=0$. That is, we have
    \begin{equation}
        \phi(1) = \phi(0) + \phi'(0) + \int_0^1 (1-t) \phi''(t) \; \mathrm{d}t.
    \end{equation}
    Here, we can evaluate the $\phi(1)$, $\phi(0)$, and $\phi'(0)$ directly. Indeed,
    \begin{equation}
        \phi(1) = \function(\geodesic^{\iso}_{\Vector,\VectorB}(1)) = \function(\VectorB), \quad \phi(0) = \function(\geodesic^{\iso}_{\Vector,\VectorB} (0)) =  \function(\Vector),
    \end{equation}
    and
    \begin{equation}
        \phi(0)' = (\nabla \function(\Vector), \dot{\geodesic}^{\iso}_{\Vector,\VectorB}(0))_2 = (\nabla \function(\Vector), \log^{\iso}_{\Vector}(\VectorB) )_2.
    \end{equation}
    So we have
    \begin{equation}
        \function(\VectorB) = \function(\Vector) + (\nabla \function(\Vector), \log^{\iso}_{\Vector}(\VectorB) )_2 + \int_0^1 (1-t) \phi''(t) \; \mathrm{d}t.
    \end{equation}
    
    Next, we use our assumptions to bound the integral and to prove the three claims.
    For (i), choose any $\Vector \neq \VectorB\in \manifold$. Then, we have $\phi''(t) = \frac{\mathrm{d}^2}{\mathrm{d}t^2} \function (\geodesic^{\iso}_{\Vector,\VectorB}(t)) \geq 0$ by assumption and have
    \begin{equation}
        \function(\VectorB) \geq \function(\Vector) + (\nabla \function(\Vector), \log^{\iso}_{\Vector}(\VectorB) )_2,
    \end{equation}
    which proves the claim.
    
    Similarly for (ii), choose any $\Vector \neq \VectorB\in \manifold$. Then, we have $\phi''(t) = \frac{\mathrm{d}^2}{\mathrm{d}t^2} \function (\geodesic^{\iso}_{\Vector,\VectorB}(t)) \geq \alpha \distance^{\iso}_{\manifold} (\Vector, \VectorB)^2$ by assumption and have
    \begin{multline}
        \function(\VectorB) \geq \function(\Vector) + (\nabla \function(\Vector), \log^{\iso}_{\Vector}(\VectorB) )_2 + \alpha \distance^{\iso}_{\manifold} (\Vector, \VectorB)^2 \int_0^1 (1 - t) \; \mathrm{d}t \\
        = \function(\Vector) + (\nabla \function(\Vector), \log^{\iso}_{\Vector}(\VectorB) )_2 + \frac{\alpha}{2} \distance^{\iso}_{\manifold} (\Vector, \VectorB)^2
    \end{multline}
    which proves the claim.
    
    Finally for (iii), choose any $\tangentVector_\Vector\in\tangent_\Vector\manifold$. Then, defining $\VectorB := \exp^{\iso}_\Vector(\tangentVector_\Vector)$ we have $\phi''(t) = \frac{\mathrm{d}^2}{\mathrm{d}t^2} \function (\geodesic^{\iso}_{\Vector,\VectorB}(t)) \leq L \distance^{\iso}_{\manifold} (\Vector, \VectorB)^2$ by assumption and have that
    \begin{multline}
        \function(\exp^{\iso}_\Vector(\tangentVector_\Vector)) = \function(\VectorB) \leq \function(\Vector) + (\nabla \function(\Vector), \log^{\iso}_{\Vector}(\VectorB) )_2 + L \distance^{\iso}_{\manifold} (\Vector, \VectorB)^2 \int_0^1 (1 - t) \; \mathrm{d}t \\
        = \function(\Vector) + (\nabla \function(\Vector), \tangentVector_\Vector )_2 + \frac{L}{2} \|\tangentVector_\Vector\|_2^2,
    \end{multline}
    which proves the claim.
\end{proof}

\begin{proof}[Proof of \Cref{lem:second-deriv-along-iso-geo}]
    The claim follows from direct evaluation. First, we have
    \begin{equation}
    \frac{\mathrm{d}}{\mathrm{d}t} \function (\geodesic^{\iso}_{\Vector,\VectorB}(t)) = (\nabla \function (\geodesic^{\iso}_{\Vector,\VectorB}(t)), \dot{\geodesic}^{\iso}_{\Vector,\VectorB}(t))_2,
    \end{equation}
    which subsequently gives
    \begin{equation}
    \frac{\mathrm{d}^2}{\mathrm{d}t^2} \function (\geodesic^{\iso}_{\Vector,\VectorB}(t)) = D_{\geodesic^{\iso}_{\Vector,\VectorB}(t)}^2 \function [\dot{\geodesic}^{\iso}_{\Vector,\VectorB}(t), \dot{\geodesic}^{\iso}_{\Vector,\VectorB}(t)] + (\nabla \function (\geodesic^{\iso}_{\Vector,\VectorB}(t)), \ddot{\geodesic}^{\iso}_{\Vector,\VectorB}(t))_2.
    \end{equation}
    We know from \cref{eq:iso-geo-speed,eq:iso-distance} that 
    \begin{equation}
    \dot{\geodesic}^{\iso}_{\Vector,\VectorB}(t) = \frac{\distance^{\iso}_{\manifold} (\Vector,\VectorB)}{\|\dot{\geodesic}_{\Vector,\VectorB}(\tau_{\Vector,\VectorB}(t))\|_2} \dot{\geodesic}_{\Vector,\VectorB}(\tau_{\Vector,\VectorB}(t))
    \end{equation}
    So
    \begin{equation}
    \ddot{\geodesic}^{\iso}_{\Vector,\VectorB}(t) = \frac{\distance^{\iso}_{\manifold} (\Vector,\VectorB)^2} {\|\dot{\geodesic}_{\Vector,\VectorB}(\tau_{\Vector,\VectorB}(t))\|_2^2} \Bigl( \ddot{\geodesic}_{\Vector,\VectorB}(\tau_{\Vector,\VectorB}(t)) - \frac{(\dot{\geodesic}_{\Vector,\VectorB}(\tau_{\Vector,\VectorB}(t)), \ddot{\geodesic}_{\Vector,\VectorB}(\tau_{\Vector,\VectorB}(t)))}{\|\dot{\geodesic}_{\Vector,\VectorB}(\tau_{\Vector,\VectorB}(t))\|_2^2} \dot{\geodesic}_{\Vector,\VectorB}(\tau_{\Vector,\VectorB}(t)) \Bigr)
    \label{eq:thm-step}
    \end{equation}
    Next, 
    \begin{multline}
    \frac{(\dot{\geodesic}_{\Vector,\VectorB}(\tau_{\Vector,\VectorB}(t)), \ddot{\geodesic}_{\Vector,\VectorB}(\tau_{\Vector,\VectorB}(t)))_2}{\|\dot{\geodesic}_{\Vector,\VectorB}(\tau_{\Vector,\VectorB}(t))\|_2^2} \dot{\geodesic}_{\Vector,\VectorB}(\tau_{\Vector,\VectorB}(t)) \\
    = \frac{\dot{\geodesic}_{\Vector,\VectorB}(\tau_{\Vector,\VectorB}(t))}{\|\dot{\geodesic}_{\Vector,\VectorB}(\tau_{\Vector,\VectorB}(t))\|_2} \otimes \frac{\dot{\geodesic}_{\Vector,\VectorB}(\tau_{\Vector,\VectorB}(t))}{\|\dot{\geodesic}_{\Vector,\VectorB}(\tau_{\Vector,\VectorB}(t))\|_2} \ddot{\geodesic}_{\Vector,\VectorB}(\tau_{\Vector,\VectorB}(t))
    \label{eq:thm-step-p1}
    \end{multline}
    and
    \begin{equation}
    \ddot{\geodesic}_{\Vector,\VectorB}(\tau_{\Vector,\VectorB}(t)) = - \Gamma_{\geodesic_{\Vector,\VectorB}(\tau_{\Vector,\VectorB}(t))} [\dot{\geodesic}_{\Vector,\VectorB}(\tau_{\Vector,\VectorB}(t)), \dot{\geodesic}_{\Vector,\VectorB}(\tau_{\Vector,\VectorB}(t))]
    \label{eq:thm-step-p2}
    \end{equation}

    Substituting \cref{eq:thm-step-p1,eq:thm-step-p2} into \cref{eq:thm-step} gives
    \begin{equation}
    \ddot{\geodesic}^{\iso}_{\Vector,\VectorB}(t) = - \Bigl(\mathbf{I} - \frac{\dot{\geodesic}^{\iso}_{\Vector,\VectorB}(t)}{\|\dot{\geodesic}^{\iso}_{\Vector,\VectorB}(t)\|_2} \otimes \frac{\dot{\geodesic}^{\iso}_{\Vector,\VectorB}(t)}{\|\dot{\geodesic}^{\iso}_{\Vector,\VectorB}(t)\|_2} \Bigr) \Gamma_{\geodesic^{\iso}_{\Vector,\VectorB}(t)} [\dot{\geodesic}^{\iso}_{\Vector,\VectorB}(t), \dot{\geodesic}^{\iso}_{\Vector,\VectorB}(t)]
    \end{equation}
    so that
    \begin{multline}
    (\nabla \function (\geodesic^{\iso}_{\Vector,\VectorB}(t)), \ddot{\geodesic}^{\iso}_{\Vector,\VectorB}(t))_2 \\
    =  - \Biggl(\Bigl(\mathbf{I} - \frac{\dot{\geodesic}^{\iso}_{\Vector,\VectorB}(t)}{\|\dot{\geodesic}^{\iso}_{\Vector,\VectorB}(t)\|_2} \otimes \frac{\dot{\geodesic}^{\iso}_{\Vector,\VectorB}(t)}{\|\dot{\geodesic}^{\iso}_{\Vector,\VectorB}(t)\|_2} \Bigr) \nabla \function (\geodesic^{\iso}_{\Vector,\VectorB}(t)), \Gamma_{\geodesic^{\iso}_{\Vector,\VectorB}(t)} [\dot{\geodesic}^{\iso}_{\Vector,\VectorB}(t), \dot{\geodesic}^{\iso}_{\Vector,\VectorB}(t)]\Biggr)_2,
    \end{multline}
    from which we conclude that
    \begin{multline}
        \frac{\mathrm{d}^2}{\mathrm{d}t^2} \function (\geodesic^{\iso}_{\Vector,\VectorB}(t)) = D_{\geodesic^{\iso}_{\Vector,\VectorB}(t)}^2 \function [\dot{\geodesic}^{\iso}_{\Vector,\VectorB}(t), \dot{\geodesic}^{\iso}_{\Vector,\VectorB}(t)] 
        \\
        - \Biggl(\Bigl(\mathbf{I} - \frac{\dot{\geodesic}^{\iso}_{\Vector,\VectorB}(t)}{\|\dot{\geodesic}^{\iso}_{\Vector,\VectorB}(t)\|_2} \otimes \frac{\dot{\geodesic}^{\iso}_{\Vector,\VectorB}(t)}{\|\dot{\geodesic}^{\iso}_{\Vector,\VectorB}(t)\|_2} \Bigr) \nabla \function (\geodesic^{\iso}_{\Vector,\VectorB}(t)), \Gamma_{\geodesic^{\iso}_{\Vector,\VectorB}(t)} [\dot{\geodesic}^{\iso}_{\Vector,\VectorB}(t), \dot{\geodesic}^{\iso}_{\Vector,\VectorB}(t)]\Biggr)_2\\
        \overset{\distance^{\iso}_{\manifold} (\Vector,\VectorB) = \|\dot{\geodesic}^{\iso}_{\Vector,\VectorB}(t)\|_2}{=} \Biggl( \frac{D_{\geodesic^{\iso}_{\Vector,\VectorB}(t)}^2 \function [\dot{\geodesic}^{\iso}_{\Vector,\VectorB}(t), \dot{\geodesic}^{\iso}_{\Vector,\VectorB}(t)]}{\|\dot{\geodesic}^{\iso}_{\Vector,\VectorB}(t)\|_2^2} \\
        - \frac{\Bigl(\bigl(\mathbf{I} - \frac{\dot{\geodesic}^{\iso}_{\Vector,\VectorB}(t)}{\|\dot{\geodesic}^{\iso}_{\Vector,\VectorB}(t)\|_2} \otimes \frac{\dot{\geodesic}^{\iso}_{\Vector,\VectorB}(t)}{\|\dot{\geodesic}^{\iso}_{\Vector,\VectorB}(t)\|_2} \bigr) \nabla \function (\geodesic^{\iso}_{\Vector,\VectorB}(t)), \Gamma_{\geodesic^{\iso}_{\Vector,\VectorB}(t)} [\dot{\geodesic}^{\iso}_{\Vector,\VectorB}(t), \dot{\geodesic}^{\iso}_{\Vector,\VectorB}(t)]\Bigr)_2}{\|\dot{\geodesic}^{\iso}_{\Vector,\VectorB}(t)\|_2^2}\Biggr) \distance^{\iso}_{\manifold} (\Vector,\VectorB)^2
    \end{multline}
\end{proof}

\begin{proof}[Proof of \Cref{prop:implications-on-manifold}]
    By \Cref{lem:second-deriv-along-iso-geo}, we know that 
    \begin{multline}
        \frac{\mathrm{d}^2}{\mathrm{d}t^2} \function (\geodesic^{\iso}_{\Vector,\VectorB}(t)) = \Biggl( \frac{D_{\geodesic^{\iso}_{\Vector,\VectorB}(t)}^2 \function [\dot{\geodesic}^{\iso}_{\Vector,\VectorB}(t), \dot{\geodesic}^{\iso}_{\Vector,\VectorB}(t)]}{\|\dot{\geodesic}^{\iso}_{\Vector,\VectorB}(t)\|_2^2} \\
        - \frac{\Bigl(\bigl(\mathbf{I} - \frac{\dot{\geodesic}^{\iso}_{\Vector,\VectorB}(t)}{\|\dot{\geodesic}^{\iso}_{\Vector,\VectorB}(t)\|_2} \otimes \frac{\dot{\geodesic}^{\iso}_{\Vector,\VectorB}(t)}{\|\dot{\geodesic}^{\iso}_{\Vector,\VectorB}(t)\|_2} \bigr) \nabla \function (\geodesic^{\iso}_{\Vector,\VectorB}(t)), \Gamma_{\geodesic^{\iso}_{\Vector,\VectorB}(t)} [\dot{\geodesic}^{\iso}_{\Vector,\VectorB}(t), \dot{\geodesic}^{\iso}_{\Vector,\VectorB}(t)]\Bigr)_2}{\|\dot{\geodesic}^{\iso}_{\Vector,\VectorB}(t)\|_2^2}\Biggr) \distance^{\iso}_{\manifold} (\Vector,\VectorB)^2.
    \end{multline}
    So, by assumption, the function $\function$ is $(\beta + \delta)$-strongly iso-g-convex by \Cref{thm:iso-convex-smooth-impl} (ii) and $(\eta+ \zeta)$-Lipschitz by \Cref{thm:iso-convex-smooth-impl} (iii), which proves the claim.
\end{proof}

\paragraph{Proofs of results in \Cref{sec:finding-mins-g-convex}}

\begin{proof}[Proof of \Cref{lem:iso-g-convex-implies-PL}]
    Let $\bar{\Vector}\in\manifold$ be a minimizer of $\function$ on $\manifold$. Then by strong iso-g-convexity we have for every $\Vector\in \manifold$   
    \begin{equation}
        \function(\bar{\Vector}) \geq \function(\Vector) + (\nabla \function(\Vector), \log^{\iso}_{\Vector}(\bar{\Vector}) )_2 + \frac{\alpha}{2} \distance^{\iso}_{\manifold} (\Vector, \bar{\Vector})^2,
    \end{equation}
    which can be rewritten into
    \begin{equation}
    \function(\Vector) - \function(\bar{\Vector}) \leq - (\nabla \function(\Vector), \log^{\iso}_{\Vector}(\bar{\Vector}) )_2 - \frac{\alpha}{2} \distance^{\iso}_{\manifold} (\Vector, \bar{\Vector})^2.
    \end{equation}
    Because $\log^{\iso}_{\Vector}(\bar{\Vector}) \in \tangent_\Vector\manifold$ we have
    \begin{equation}
    (\nabla \function(\Vector), \log^{\iso}_{\Vector}(\bar{\Vector}) )_2 = (\nabla \function(\Vector), \mathbf{P}_{\Vector}\log^{\iso}_{\Vector}(\bar{\Vector}) )_2 = (\mathbf{P}_{\Vector}\nabla \function(\Vector), \log^{\iso}_{\Vector}(\bar{\Vector}) )_2.
    \end{equation}
    In addition, we have that
    \begin{equation}
    |(\mathbf{P}_{\Vector}\nabla \function(\Vector), \log^{\iso}_{\Vector}(\bar{\Vector}) )_2| \leq \|\mathbf{P}_{\Vector}\nabla \function(\Vector)\|_2 \|\log^{\iso}_{\Vector}(\bar{\Vector}) )_2\|_2 = \|\mathbf{P}_{\Vector}\nabla \function(\Vector)\|_2 \distance^{\iso}_{\manifold} (\Vector, \bar{\Vector}).
    \end{equation}
    So we have that
    \begin{equation}
        \function(\Vector) - \function(\bar{\Vector}) \leq \|\mathbf{P}_{\Vector}\nabla \function(\Vector)\|_2 \distance^{\iso}_{\manifold} (\Vector, \bar{\Vector}) - \frac{\alpha}{2} \distance^{\iso}_{\manifold} (\Vector, \bar{\Vector})^2. 
    \end{equation}
    Next, note that the mapping
    \begin{equation}
    r\mapsto \|\mathbf{P}_{\Vector}\nabla \function(\Vector)\|_2 r - \frac{\alpha}{2} r^2
    \end{equation}
    is maximized for $r = \frac{1}{\alpha}\|\mathbf{P}_{\Vector}\nabla \function(\Vector)\|_2$, from which we conclude that
    \begin{equation}
    \|\mathbf{P}_{\Vector}\nabla \function(\Vector)\|_2 \distance^{\iso}_{\manifold} (\Vector, \bar{\Vector}) - \frac{\alpha}{2} \distance^{\iso}_{\manifold} (\Vector, \bar{\Vector})^2 \leq \frac{1}{2\alpha} \|\mathbf{P}_{\Vector}\nabla \function(\Vector)\|_2^2,
    \end{equation}
    which gives the PL property
    \begin{equation}
        \function(\Vector) - \function(\bar{\Vector}) \leq \frac{1}{2\alpha} \|\mathbf{P}_{\Vector}\nabla \function(\Vector)\|_2^2
    \end{equation}
    as claimed.
\end{proof}

\begin{proof}[Proof of \Cref{thm:conv-pg-ird}]
    Using iso-$L$-smoothness for $\tangentVector_\Vector := - r \mathbf{P}_{\Vector^{(\sumIndC)}} \nabla f (\Vector^{(\sumIndC)})$ we have
    \begin{equation}
        \function(\Vector^{(\sumIndC+1)}) \leq \function(\Vector^{(\sumIndC)}) - r (\nabla \function(\Vector^{(\sumIndC)}), \mathbf{P}_{\Vector^{(\sumIndC)}} \nabla f (\Vector^{(\sumIndC)}) )_2 + \frac{r^2 L}{2}\|\mathbf{P}_{\Vector^{(\sumIndC)}} \nabla f (\Vector^{(\sumIndC)})\|_2^2.
    \end{equation}
    Since
    \begin{equation}
    (\nabla \function(\Vector^{(\sumIndC)}), \mathbf{P}_{\Vector^{(\sumIndC)}} \nabla f (\Vector^{(\sumIndC)}) )_2 = (\mathbf{P}_{\Vector^{(\sumIndC)}}\nabla \function(\Vector^{(\sumIndC)}), \mathbf{P}_{\Vector^{(\sumIndC)}} \nabla f (\Vector^{(\sumIndC)}) )_2 = \|\mathbf{P}_{\Vector^{(\sumIndC)}} \nabla f (\Vector^{(\sumIndC)})\|_2^2,
    \end{equation}
    we have that
    \begin{equation}
        \function(\Vector^{(\sumIndC+1)}) \leq \function(\Vector^{(\sumIndC)}) - r \|\mathbf{P}_{\Vector^{(\sumIndC)}} \nabla f (\Vector^{(\sumIndC)})\|_2^2 + \frac{r^2L}{2}\|\mathbf{P}_{\Vector^{(\sumIndC)}} \nabla f (\Vector^{(\sumIndC)})\|_2^2
    \end{equation}
    and because $0< r \leq \frac{1}{L} \Rightarrow L \leq \frac{1}{r}$ it follows that
    \begin{equation}
    \function(\Vector^{(\sumIndC+1)}) \leq \function(\Vector^{(\sumIndC)}) - \frac{r}{2} \|\mathbf{P}_{\Vector^{(\sumIndC)}} \nabla f (\Vector^{(\sumIndC)})\|_2^2 .
    \end{equation}
    Subsequently, since $\function$ is PL by \Cref{lem:iso-g-convex-implies-PL}, we have
    \begin{equation}
    \function(\Vector^{(\sumIndC+1)}) \leq \function(\Vector^{(\sumIndC)}) - \frac{r}{2} (2 \alpha (\function(\Vector^{(\sumIndC)}) - \function^*)) = \function(\Vector^{(\sumIndC)}) - \alpha r (\function(\Vector^{(\sumIndC)}) - \function^*),
    \end{equation}
    which can be rewritten as
    \begin{equation}
    \function(\Vector^{(\sumIndC+1)}) - \function^* \leq (1 - \alpha r) (\function(\Vector^{(\sumIndC)}) - \function^* ),
    \end{equation}
    from which follows that
    \begin{equation}
    \function(\Vector^{(\sumIndC)}) - \function^* \leq (1 - \alpha r)^k (\function(\Vector^{(0)}) - \function^* ),
    \end{equation}
    and which implies that $\function(\Vector^{(\sumIndC)}) \to \function^*$.

    To show convergence of iterates, we invoke strong iso-g-convexity once more. In particular, we have 
    \begin{equation}
    \function(\Vector^{(\sumIndC)}) \geq \function(\bar{\Vector}) + (\nabla \function(\bar{\Vector}), \log^{\iso}_{\bar{\Vector}}(\Vector^{(\sumIndC)}) )_2 + \frac{\alpha}{2} \distance^{\iso}_{\manifold} (\bar{\Vector}, \Vector^{(\sumIndC)})^2
    \end{equation}
    in which the second term on the right hand side vanishes because
    \begin{equation}
    (\nabla \function(\bar{\Vector}), \log^{\iso}_{\bar{\Vector}}(\Vector^{(\sumIndC)}) )_2 = (\nabla \function(\bar{\Vector}), \mathbf{P}_{\bar{\Vector}}\log^{\iso}_{\bar{\Vector}}(\Vector^{(\sumIndC)}) )_2 = (\mathbf{P}_{\bar{\Vector}}\nabla \function(\bar{\Vector}), \log^{\iso}_{\bar{\Vector}}(\Vector^{(\sumIndC)}))_2 = 0.
    \end{equation}
    In other words, we have
    \begin{equation}
    \function(\Vector^{(\sumIndC)}) - \function^* \geq \frac{\alpha}{2} \distance^{\iso}_{\manifold} (\bar{\Vector}, \Vector^{(\sumIndC)})^2
    \end{equation}
    
    Finally, we have
    \begin{multline}
    \|\bar{\Vector} - \Vector^{(\sumIndC)} \|_2^2 = \Bigl(\inf_{\geodesic_{\bar{\Vector}, \Vector^{(\sumIndC)}}} \int_{0}^1 \|\dot{\geodesic}_{\bar{\Vector}, \Vector^{(\sumIndC)}} (t)\|_2 \; \mathrm{d}t \Bigr)^2 \leq \distance^{\iso}_{\manifold} (\bar{\Vector}, \Vector^{(\sumIndC)})^2 \\
    \leq \frac{2}{\alpha} (\function(\Vector^{(\sumIndC)}) - \function^*) \leq \frac{2}{\alpha} (1 - \alpha r)^k (\function(\Vector^{(0)}) - \function^* )
    \end{multline}
    which proves the claim that $\Vector^{(\sumIndC)}\to \bar{\Vector}$ in $\ell^2$.
\end{proof}

\subsection{Details to numerical experiments in \Cref{sec:inverse-problems-appl}}
\label{app:numerics-function-convergence}

\subsubsection{Banana pullback geometry}
The banana diffeomorphism $\diffeo_{\text{banana}} : \Real^2 \to \Real^2$ is defined as
\begin{equation}
    \diffeo_{\text{banana}}(\Vector) := \left( \Vector_1 - a \Vector_2^2 - z,\, \Vector_2 \right),
\end{equation}
where $a \in \Real$ is the shear parameter and $z \in \Real$ is an offset. Its inverse is
\begin{equation}
    \diffeo_{\text{banana}}^{-1}(\VectorB) := \left( \VectorB_1 + a \VectorB_2^2 + z,\, \VectorB_2 \right).
\end{equation}
In the experiments we used $a = 1/9$ and $z= 0.0$.

\subsubsection{Additional results}
\begin{figure}[h!]
    \centering
    \setlength{\tabcolsep}{2pt}
    \renewcommand{\arraystretch}{0}
    \begin{tabular}{cc}
        \includegraphics[width=0.48\textwidth]{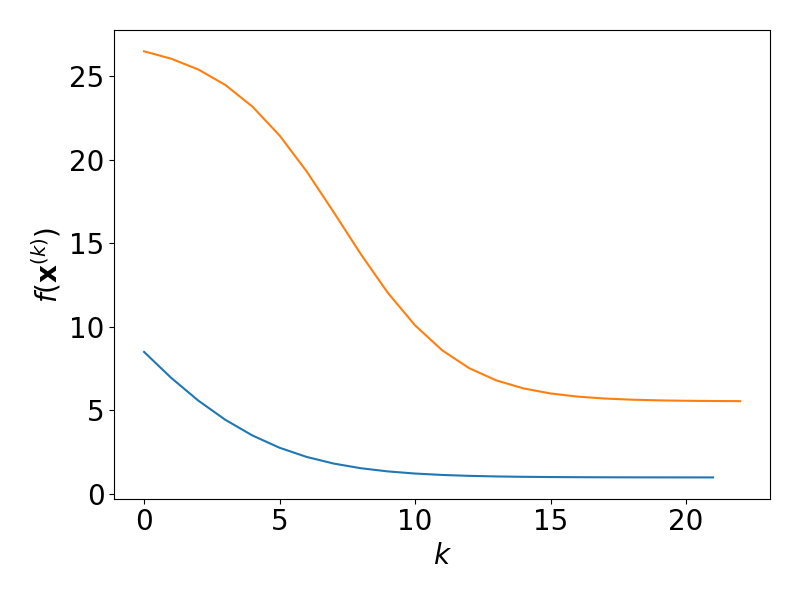} &
        \includegraphics[width=0.48\textwidth]{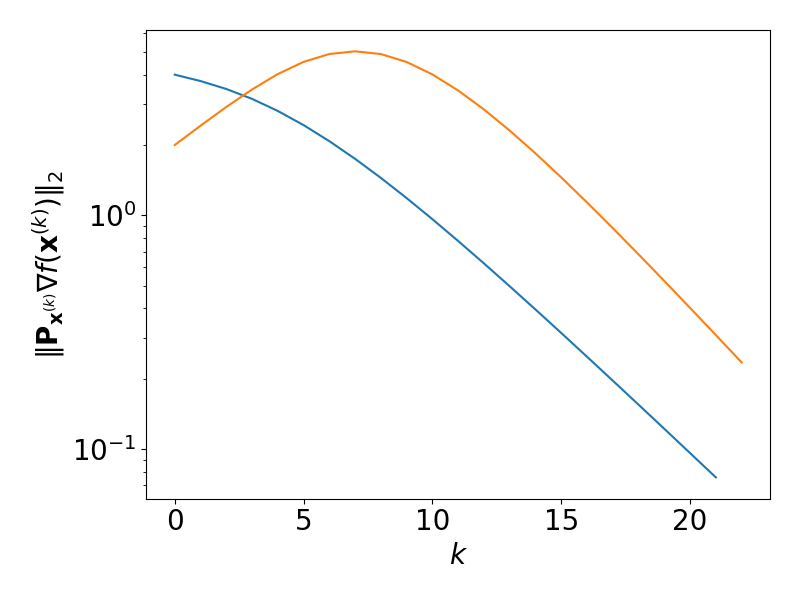} 
    \end{tabular}
    \caption{Convergence results on the synthetic data manifold. We observe monotone decay of the loss (left) and local linear convergence of the projected gradient norm (right), which is in line with our theory.}
    \label{fig:conv-results-opt-toy}
\end{figure}

\begin{figure}[h!]
    \centering
    \setlength{\tabcolsep}{2pt}
    \renewcommand{\arraystretch}{0}
    \begin{tabular}{cc}
        \includegraphics[width=0.48\textwidth]{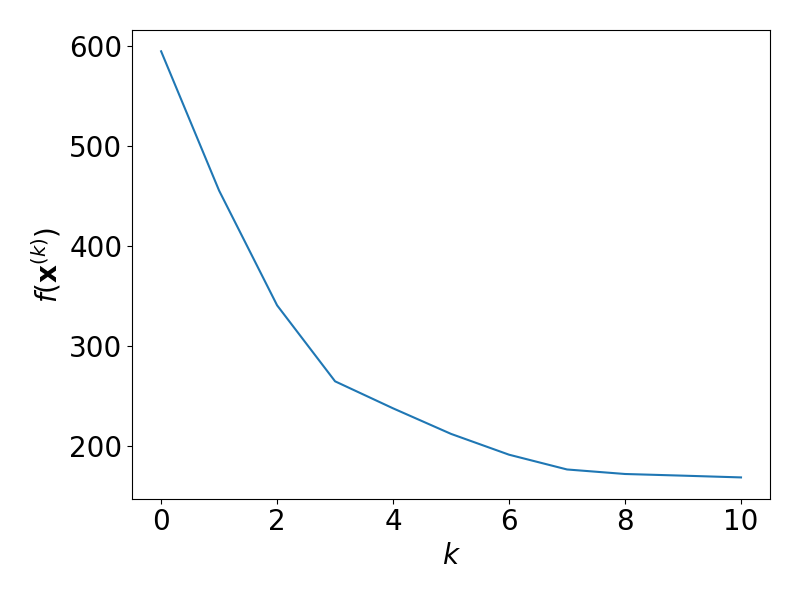} &
        \includegraphics[width=0.48\textwidth]{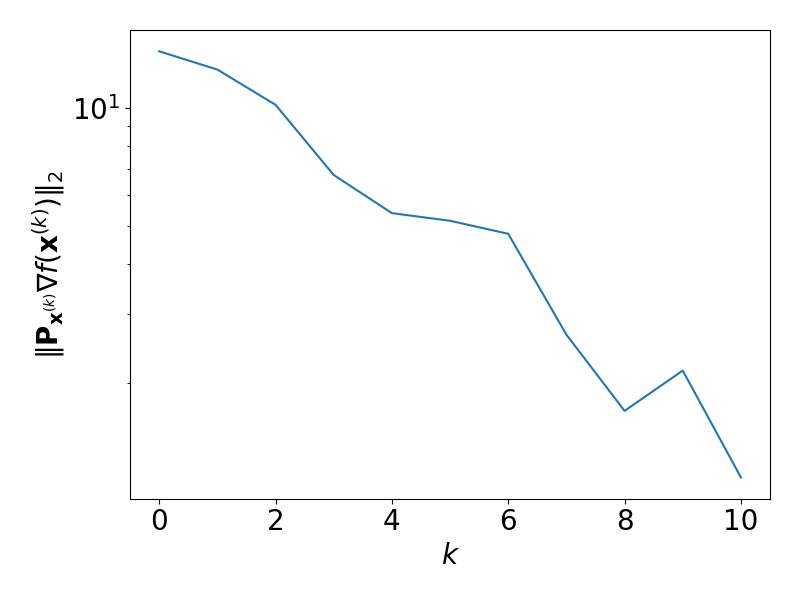} 
    \end{tabular}
    \caption{Convergence results on the learned mnist data manifold. We observe monotone decay of the loss (left) and perturbed linear convergence of the projected gradient norm (right), which is most likely due to numerical errors.}
    \label{fig:conv-results-opt-mnist}
\end{figure}

The convergence plots of \cref{eq:iso-first-order-scheme} for the banana (with step size $r=0.1$) and MNIST (with step size $r=1$) data are shown in \Cref{fig:conv-results-opt-toy,fig:conv-results-opt-mnist} and are in line with our expectations. That is, we attain linear convergence in 
$
k\mapsto \|\mathbf{P}_{\Vector^{(\sumIndC)}} \nabla f (\Vector^{(\sumIndC)})\|_2$. The river and spiral experiments were run until $\frac{\|\mathbf{P}_{\Vector^{(\sumIndC)}} \nabla f (\Vector^{(\sumIndC)})\|_2}{f(\Vector^{(0)})} < 10^{-2}$, whereas for MNIST the algorithm was terminated once 
$\frac{\|\mathbf{P}_{\Vector^{(\sumIndC)}} \nabla f (\Vector^{(\sumIndC)})\|_2}{f(\Vector^{(0)})} < 2 \cdot 10^{-3}$

\section{Supplementary material to \Cref{sec:vectorfield-convergence}}
\label{app:vectorfield-convergence}

\subsection{Proofs of the results in \Cref{sec:iso-monotone,sec:finding-mins-monotone}}
\label{app:supplement-vector-field-convergence-proofs}

\paragraph{Proofs of the results in \Cref{sec:iso-monotone}}

\begin{proof}[Proof of \Cref{prop:equivalence-iso-monotone-smooth-R}]
    First, we note that on $\Real$ all iso-manifold mappings reduce to simpler expressions, i.e., for $x,y \in \Real$ and $\xi_x \in \tangent_x \Real$ we have
    \begin{align}
        \geodesic^{\iso}_{x,y}(t) &= (1 - t) x + t y,\\
        \log^{\iso}_{x} (y) &= y - x,\\
        \exp^{\iso}_{x} (\xi_x) &= x + \xi_x,\\
        \distance^{\iso}_{\Real}(x,y) & = |x - y|,\\
        \partransport^{\iso}_{y \leftarrow x}\xi_x &= \xi_x.
    \end{align}
    So for showing (i) we note that
    \begin{equation}
        (f'(y) - \partransport^{\iso}_{y\leftarrow x} f'(x), \partransport^{\iso}_{y\leftarrow x} \log^{\iso}_{x} (y))_2 \geq \alpha \distance^{\iso}_{\Real} (x, y)^2 
    \end{equation}
    is equivalent to
    \begin{equation}
         (f'(y) -  f'(x), y - x)_2 \geq \alpha (x-y)^2.
    \end{equation}
    In other words, ($\alpha$-strong) iso-monotonicity of $f'$ implies monotonicity and vice versa, which proves the first claim.

    Similarly, for showing (ii) we also have that
    \begin{equation}
        |f'(y) - \partransport^{\iso}_{y\leftarrow x} f'(x)| \leq L \distance^{\iso}_{\Real} (x, y)
    \end{equation}
    is equivalent to
    \begin{equation}
        |f'(y) -  f'(x)| \leq L |x - y|.
    \end{equation}
    So, $L$-iso-Lipschitzness of $f'$ implies $L$-Lipschitzness of $f'$ and vice versa, which proves the second claim.
\end{proof}

\begin{proof}[Proof of \Cref{prop:taylor-iso}]
    Choose any $\tangentVectorB_{\VectorB} \in \tangent_{\VectorB} \Real^\dimInd$ and define $\phi: [0,1]\to \Real$ as
    \begin{equation}
        \phi(t) := (\frac{\|\dot{\geodesic}_{\Vector, \VectorB}(t)\|_2}{\|\log_{\VectorB} (\Vector)\|_2} \tangentVector_{\geodesic_{\Vector, \VectorB}(t)}, \tangentVectorB_{\geodesic_{\Vector, \VectorB}(t)})_{\geodesic_{\Vector, \VectorB}(t)},
    \end{equation}
    where $\tangentVectorB_{\geodesic_{\Vector, \VectorB}(t)} := \partransport_{\geodesic_{\Vector, \VectorB}(t) \leftarrow \VectorB} \tangentVectorB_{\VectorB}$.

    By the fundamental theorem of calculus we have that
    \begin{equation}
        \phi(1) - \phi(0) = \int_{0}^{1} \frac{\mathrm{d}}{\mathrm{d}t} \phi(t) \; \mathrm{d}t.
        \label{eq:lem-taylor-iso-simple}
    \end{equation}
    For the left hand side of (\ref{eq:lem-taylor-iso-simple}) we have
    \begin{multline}
        \phi(0) = (\frac{\|\dot{\geodesic}_{\Vector, \VectorB}(0)\|_2}{\|\log_{\VectorB} (\Vector)\|_2} \tangentVector_{\geodesic_{\Vector, \VectorB}(0)}, \tangentVectorB_{\geodesic_{\Vector, \VectorB}(0)})_{\geodesic_{\Vector, \VectorB}(0)} 
        = (\frac{\|\log_{\Vector} (\VectorB)\|_2}{{\|\log_{\VectorB} (\Vector)\|_2}} \tangentVector_{\Vector}, \partransport_{\Vector \leftarrow \VectorB} \tangentVectorB_{\VectorB})_{\Vector} \\
        = (\frac{\|\log_{\Vector} (\VectorB)\|_2}{{\|\log_{\VectorB} (\Vector)\|_2}} \partransport_{\VectorB \leftarrow \Vector} \tangentVector_{\Vector}, \tangentVectorB_{\VectorB})_{\VectorB} 
        = (\partransport^{\iso}_{\VectorB \leftarrow \Vector} \tangentVector_{\Vector}, \tangentVectorB_{\VectorB})_{\VectorB},
        \label{eq:lem-taylor-iso-simple-1}
    \end{multline}
    and
    \begin{multline}
        \phi(1) = (\frac{\|\dot{\geodesic}_{\Vector, \VectorB}(1)\|_2}{\|\log_{\VectorB} (\Vector)\|_2}  \tangentVector_{\geodesic_{\Vector, \VectorB}(1)}, \tangentVectorB_{\geodesic_{\Vector, \VectorB}(1)})_{\geodesic_{\Vector, \VectorB}(1)} \\
        = (\frac{\|-\log_{\VectorB} (\Vector)\|_2}{\|\log_{\VectorB} (\Vector)\|_2} \tangentVector_{\VectorB}, \tangentVectorB_{\VectorB})_{\VectorB} = ( \tangentVector_{\VectorB}, \tangentVectorB_{\VectorB})_{\VectorB}.
        \label{eq:lem-taylor-iso-simple-2}
    \end{multline}
    Then, for the right hand side of (\ref{eq:lem-taylor-iso-simple}) we have
    \begin{multline}
        \int_{0}^{1} \frac{\mathrm{d}}{\mathrm{d}t} \phi(t) \; \mathrm{d}t 
        =  \int_{0}^{1} \frac{\mathrm{d}}{\mathrm{d}t} (\frac{\|\dot{\geodesic}_{\Vector, \VectorB}(t)\|_2}{\|\log_{\VectorB} (\Vector)\|_2} \tangentVector_{\geodesic_{\Vector, \VectorB}(t)}, \tangentVectorB_{\geodesic_{\Vector, \VectorB}(t)})_{\geodesic_{\Vector, \VectorB}(t)} \; \mathrm{d}t \\
        =  \int_{0}^{1}  \dot{\geodesic}_{\Vector, \VectorB}(t) (\frac{\|\dot{\geodesic}_{\Vector, \VectorB}(t)\|_2}{\|\log_{\VectorB} (\Vector)\|_2} \tangentVector_{(\cdot)}, \tangentVectorB_{(\cdot)})_{(\cdot)} \; \mathrm{d}t \\
        \overset{\nabla_{\dot{\geodesic}_{\Vector, \VectorB}(t)} \tangentVectorB_{(\cdot)}=\mathbf{0} }{=}  \int_{0}^{1}  \frac{1}{\|\log_{\VectorB} (\Vector)\|_2} (\nabla_{\dot{\geodesic}_{\Vector, \VectorB}(t)} \|\dot{\geodesic}_{\Vector, \VectorB}\|_2 \tangentVector_{(\cdot)}, \tangentVectorB_{\geodesic_{\Vector, \VectorB}(t)})_{\geodesic_{\Vector, \VectorB}(t)} \; \mathrm{d}t \\
        =  \int_{0}^{1} \frac{1}{\|\log_{\VectorB} (\Vector)\|_2} (\nabla_{\dot{\geodesic}_{\Vector, \VectorB}(t)}\|\partransport_{(\cdot)\leftarrow \geodesic_{\Vector, \VectorB}(t)}\dot{\geodesic}_{\Vector, \VectorB} (t)\|_2\tangentVector_{(\cdot)}, \tangentVectorB_{\geodesic_{\Vector, \VectorB}(t)})_{\geodesic_{\Vector, \VectorB}(t)} \; \mathrm{d}t\\
        =  \int_{0}^{1}  (\frac{\|\dot{\geodesic}_{\Vector, \VectorB}(t)\|_2}{\|\log_{\VectorB} (\Vector)\|_2} \nabla^{\iso}_{\dot{\geodesic}_{\Vector, \VectorB}(t)}\tangentVector_{(\cdot)}, \partransport_{\geodesic_{\Vector, \VectorB}(t) \leftarrow \VectorB} \tangentVectorB_{\VectorB})_{\geodesic_{\Vector, \VectorB}(t)} \; \mathrm{d}t \\
        =  \int_{0}^{1}  (\frac{\|\dot{\geodesic}_{\Vector, \VectorB}(t)\|_2}{\|\log_{\VectorB} (\Vector)\|_2} \partransport_{\VectorB \leftarrow \geodesic_{\Vector, \VectorB}(t)} \nabla^{\iso}_{\dot{\geodesic}_{\Vector, \VectorB}(t)}\tangentVector_{(\cdot)}, \tangentVectorB_{\VectorB})_{\VectorB} \; \mathrm{d}t.
        \label{eq:lem-taylor-iso-simple-3}
    \end{multline}
    Since $\tangentVectorB_{\VectorB}$ was arbitrary, the claim (\ref{eq:lem-taylor-iso}) follows from substituting (\ref{eq:lem-taylor-iso-simple-1}), (\ref{eq:lem-taylor-iso-simple-2}) and (\ref{eq:lem-taylor-iso-simple-3}) into (\ref{eq:lem-taylor-iso-simple}).
\end{proof}

\begin{proof}[Proof of \Cref{thm:iso-monotone-impl}]
    The claim (ii) follows from direct computation. That is, by \Cref{prop:taylor-iso}
    \begin{multline}
        (\tangentVector_{\VectorB} - \partransport^{\iso}_{\VectorB\leftarrow \Vector} \tangentVector_{\Vector}, \partransport^{\iso}_{\VectorB\leftarrow \Vector} \log^{\iso}_{\Vector} (\VectorB))_2 \\
         = \int_{0}^{1} ( \frac{\|\dot{\geodesic}_{\Vector, \VectorB}(t)\|_2}{\|\log_{\VectorB} (\Vector)\|_2} \partransport_{\VectorB\leftarrow \geodesic_{\Vector, \VectorB}(t)}   \nabla^{\iso}_{\dot{\geodesic}_{\Vector, \VectorB}(t)} \tangentVector_{(\cdot)}, \frac{\|\log_{\Vector} (\VectorB)\|_2}{\|\log_{\VectorB} (\Vector)\|_2} \partransport_{\VectorB\leftarrow \Vector} \frac{\distance^{\iso}_{\manifold}(\Vector, \VectorB)}{\|\log_{\Vector} (\VectorB)\|_2} \log_{\Vector} (\VectorB) ) \; \mathrm{d}t \\
        = \frac{\distance^{\iso}_{\manifold}(\Vector, \VectorB)}{\|\log_{\VectorB} (\Vector)\|_2^2} \int_{0}^{1} \|\dot{\geodesic}_{\Vector, \VectorB}(t)\|_2 (  \partransport_{\VectorB\leftarrow \geodesic_{\Vector, \VectorB}(t)}   \nabla^{\iso}_{\dot{\geodesic}_{\Vector, \VectorB}(t)} \tangentVector_{(\cdot)}, \partransport_{\VectorB\leftarrow \Vector} \log_{\Vector} (\VectorB) ) \; \mathrm{d}t\\
        \overset{(\ref{eq:lem-iso-alpha-mono-impl})}{\geq}  \frac{\distance^{\iso}_{\manifold}(\Vector, \VectorB)}{\|\log_{\VectorB} (\Vector)\|_2^2} \int_{0}^{1} \|\dot{\geodesic}_{\Vector, \VectorB}(t)\|_2 ( \alpha  \|\log_{\VectorB}(\Vector)\|_2^2 ) \; \mathrm{d}t \\
         = \alpha \distance^{\iso}_{\manifold}(\Vector, \VectorB) \int_{0}^{1} \|\dot{\geodesic}_{\Vector, \VectorB}(t)\|_2 \; \mathrm{d}t  =  \alpha \distance^{\iso}_{\manifold} (\Vector, \VectorB)^2.
    \end{multline}
    The claim (i) follows analogously.

    Next, the claim (iii) also follows from direct computation. That is, by \Cref{prop:taylor-iso}
    \begin{multline}
        \| \tangentVector_{\VectorB} - \partransport^{\iso}_{\VectorB\leftarrow \Vector}\tangentVector_{\Vector} \|_2
        \leq  \int_{0}^{1} \|\frac{\|\dot{\geodesic}_{\Vector, \VectorB}(t)\|_2}{\|\log_{\VectorB} (\Vector)\|_2} \partransport_{\VectorB\leftarrow \geodesic_{\Vector, \VectorB}(t)}   \nabla^{\iso}_{\dot{\geodesic}_{\Vector, \VectorB}(t)} \tangentVector_{(\cdot)}\|_2 \; \mathrm{d}t \\
        = \frac{1}{\|\log_{\VectorB} (\Vector)\|_2} \int_{0}^{1} \|\dot{\geodesic}_{\Vector, \VectorB}(t)\|_2  \| \partransport_{\VectorB\leftarrow \geodesic_{\Vector, \VectorB}(t)}   \nabla^{\iso}_{\dot{\geodesic}_{\Vector, \VectorB}(t)} \tangentVector_{(\cdot)}\|_2 \; \mathrm{d}t \\
        \overset{(\ref{eq:lem-iso-L-Lips-impl})}{\leq} \frac{1}{\|\log_{\VectorB} (\Vector)\|_2} \int_{0}^{1} \|\dot{\geodesic}_{\Vector, \VectorB}(t)\|_2 (L \|\log_{\VectorB} (\Vector)\|_2) \; \mathrm{d}t \\
        = L \int_{0}^{1} \|\dot{\geodesic}_{\Vector, \VectorB}(t)\|_2 \; \mathrm{d}t= L \distance^{\iso}_{\manifold} (\Vector, \VectorB).
    \end{multline}
\end{proof}

\begin{proof}[Proof of \Cref{lem:second-deriv-along-iso-geo-vf}]
    The claim follows from direct computation where we use that we can write 
    \begin{equation}
        \mathbf{P}_{\geodesic_{\Vector, \VectorB}(t)} = \frac{1}{\|\dot{\geodesic}_{\Vector, \VectorB}(t)\|_2^2} \dot{\geodesic}_{\Vector, \VectorB}(t) \otimes \dot{\geodesic}_{\Vector, \VectorB}(t), \quad t\in [0,1].
    \end{equation}
    That is, we have
    \begin{multline}
        \nabla^{\iso}_{\dot{\geodesic}_{\Vector, \VectorB}(t)}\mathbf{P}_{(\cdot)} \nabla f = \frac{1}{\|\dot{\geodesic}_{\Vector, \VectorB}(t)\|_2}  \nabla_{\dot{\geodesic}_{\Vector, \VectorB}(t)} \Bigl(\|\dot{\geodesic}_{\Vector, \VectorB}\|_2 \mathbf{P}_{(\cdot)} \nabla f \Bigr) \\
        = \frac{1}{\|\dot{\geodesic}_{\Vector, \VectorB}(t)\|_2}  \nabla_{\dot{\geodesic}_{\Vector, \VectorB}(t)} \Bigl( \frac{(\nabla f, \dot{\geodesic}_{\Vector, \VectorB})_2}{\|\dot{\geodesic}_{\Vector, \VectorB}\|_2} \dot{\geodesic}_{\Vector, \VectorB} \Bigr) \\
        = \frac{1}{\|\dot{\geodesic}_{\Vector, \VectorB}(t)\|_2}  \frac{\mathrm{d}}{\mathrm{d}t} \Bigl( \frac{(\nabla f (\geodesic_{\Vector, \VectorB} (t)), \dot{\geodesic}_{\Vector, \VectorB} (t))_2}{\|\dot{\geodesic}_{\Vector, \VectorB} (t)\|_2}  \Bigr) \dot{\geodesic}_{\Vector, \VectorB} (t) \\
        =  \frac{(\nabla f (\geodesic_{\Vector, \VectorB} (t)), \ddot{\geodesic}_{\Vector, \VectorB} (t))_2 + (D_{\geodesic_{\Vector, \VectorB} (t)} \nabla f[\dot{\geodesic}_{\Vector, \VectorB} (t)], \dot{\geodesic}_{\Vector, \VectorB} (t))_2}{\|\dot{\geodesic}_{\Vector, \VectorB} (t)\|_2^2} \dot{\geodesic}_{\Vector, \VectorB} (t)\\
         -  \frac{(\nabla f (\geodesic_{\Vector, \VectorB} (t)), \dot{\geodesic}_{\Vector, \VectorB} (t))_2 (\dot{\geodesic}_{\Vector, \VectorB} (t), \ddot{\geodesic}_{\Vector, \VectorB} (t))_2}{\|\dot{\geodesic}_{\Vector, \VectorB} (t)\|_2^4} \dot{\geodesic}_{\Vector, \VectorB} (t)\\
        = \Bigl( \frac{D_{\Vector}^2 f [\dot{\geodesic}_{\Vector, \VectorB}(t),\dot{\geodesic}_{\Vector, \VectorB}(t)]}{\|\dot{\geodesic}_{\Vector, \VectorB}(t)\|_2^2} 
        + \frac{( (\mathbf{I} - \mathbf{P}_{\geodesic_{\Vector, \VectorB}(t)})\nabla f (\geodesic_{\Vector, \VectorB}(t)), \ddot{\geodesic}_{\Vector, \VectorB} (t) )_2}{\|\dot{\geodesic}_{\Vector, \VectorB} (t)\|_2^2} \Bigr) \dot{\geodesic}_{\Vector, \VectorB}(t)\\
        = \Bigl( \frac{D_{\Vector}^2 f [\dot{\geodesic}_{\Vector, \VectorB}(t),\dot{\geodesic}_{\Vector, \VectorB}(t)]}{\|\dot{\geodesic}_{\Vector, \VectorB}(t)\|_2^2} 
        \\
        - \frac{( (\mathbf{I} - \mathbf{P}_{\geodesic_{\Vector, \VectorB}(t)})\nabla f (\geodesic_{\Vector, \VectorB}(t)), \Gamma_{\geodesic_{\Vector, \VectorB}(t)}(\dot{\geodesic}_{\Vector, \VectorB} (t), \dot{\geodesic}_{\Vector, \VectorB} (t)) )_2}{\|\dot{\geodesic}_{\Vector, \VectorB} (t)\|_2^2} \Bigr) \dot{\geodesic}_{\Vector, \VectorB}(t).
    \end{multline}
\end{proof}

\begin{proof}[Proof of \Cref{prop:implications-on-manifold-vf}]
    By \Cref{lem:second-deriv-along-iso-geo-vf}, we know that 
    \begin{multline}
        \nabla^{\iso}_{\dot{\geodesic}_{\Vector, \VectorB}(t)} \mathbf{P}_{(\cdot)} \nabla \function = \Biggl( \frac{D_{\geodesic_{\Vector,\VectorB}(t)}^2 \function [\dot{\geodesic}_{\Vector,\VectorB}(t), \dot{\geodesic}_{\Vector,\VectorB}(t)]}{\|\dot{\geodesic}_{\Vector,\VectorB}(t)\|_2^2} \\
        - \frac{\bigl(\bigl(\mathbf{I} - \frac{\dot{\geodesic}_{\Vector,\VectorB}(t)}{\|\dot{\geodesic}_{\Vector,\VectorB}(t)\|_2} \otimes \frac{\dot{\geodesic}_{\Vector,\VectorB}(t)}{\|\dot{\geodesic}_{\Vector,\VectorB}(t)\|_2} \bigr) \nabla \function (\geodesic_{\Vector,\VectorB}(t)), \Gamma_{\geodesic_{\Vector,\VectorB}(t)} [\dot{\geodesic}_{\Vector,\VectorB}(t), \dot{\geodesic}_{\Vector,\VectorB}(t)]\bigr)_2}{\|\dot{\geodesic}_{\Vector,\VectorB}(t)\|_2^2} \Biggr) \dot{\geodesic}_{\Vector,\VectorB}(t).
    \end{multline}
    So, by assumption, the function $\function$ is $(\beta + \delta)$-strongly iso-g-convex by \Cref{thm:iso-monotone-impl} (ii) and $(\eta+ \zeta)$-Lipschitz by \Cref{thm:iso-monotone-impl} (iii), which proves the claim.
\end{proof}

\begin{proof}[Proof of \Cref{prop:iso-mono-geo}]
    By (i) and (iii) in \Cref{thm:iso-monotone-impl} it suffices to show that for every $\Vector \neq \VectorB \in \geodesic_{\VectorD, \VectorE}$
    \begin{equation}
        \nabla^{\iso}_{\dot{\geodesic}_{\Vector, \VectorB}(t)}\dot{\geodesic}^{\iso}_{\VectorD, \VectorE} = \mathbf{0}, \quad t \in [0,1],
        \label{eq:prop-iso-mono-geo-claim}
    \end{equation}
    because then
    \begin{equation}
        (\partransport_{\VectorB\leftarrow \geodesic_{\Vector, \VectorB}(t)} \nabla^{\iso}_{\dot{\geodesic}_{\Vector, \VectorB}(t)}\dot{\geodesic}^{\iso}_{\VectorD, \VectorE}, \partransport_{\VectorB\leftarrow \Vector} \log_{\Vector} (\VectorB))_2 = 0,
    \end{equation}
    which proves the claim.

    The claim (\ref{eq:prop-iso-mono-geo-claim}) follows from direct computation, where we assume without loss of generality that $\Vector$ is closest to $\VectorD$ and $\VectorB$ closest to $\VectorE$. That is,
    \begin{align}
        \nabla^{\iso}_{\dot{\geodesic}_{\Vector, \VectorB}(t)}\dot{\geodesic}^{\iso}_{\VectorD, \VectorE} &= \frac{1}{\|\dot{\geodesic}_{\Vector, \VectorB}(t)\|_2} \nabla_{\dot{\geodesic}_{\Vector, \VectorB}(t)} \|\partransport_{(\cdot) \leftarrow \Vector} \dot{\geodesic}_{\Vector, \VectorB}(t)\|_2 \dot{\geodesic}^{\iso}_{\VectorD, \VectorE}\\
        &= \frac{1}{\|\dot{\geodesic}_{\Vector, \VectorB}(t)\|_2} \nabla_{\dot{\geodesic}_{\Vector, \VectorB}(t)} \|\dot{\geodesic}_{\Vector, \VectorB}\|_2 \frac{\int_{0}^{1} \|\dot{\geodesic}_{\VectorD, \VectorE}(s)\|_2 \; \mathrm{d}s}{\|\dot{\geodesic}_{\VectorD, \VectorE}\|_2} \dot{\geodesic}_{\VectorD, \VectorE}\\
        &\overset{\dot{\geodesic}_{\VectorD, \VectorE}\mid_{\geodesic_{\Vector, \VectorB}} = \frac{\distance_{\manifold}(\VectorE, \VectorE)}{\distance_{\manifold}(\Vector, \VectorB)} \dot{\geodesic}_{\Vector, \VectorB}}{=} \frac{\int_{0}^{1} \|\dot{\geodesic}_{\VectorD, \VectorE}(s)\|_2 \; \mathrm{d}s}{\|\dot{\geodesic}_{\Vector, \VectorB}(t)\|_2} \nabla_{\dot{\geodesic}_{\Vector, \VectorB}(t)} \dot{\geodesic}_{\Vector, \VectorB} = \mathbf{0}.
    \end{align}
\end{proof}

\begin{proof}[Proof of \Cref{prop:iso-mono-log-field}]
    By (ii) and (iii) in \Cref{thm:iso-monotone-impl} it suffices to show that for every $\Vector\neq \VectorB\in \manifold$ and $t\in [0,1]$
    \begin{equation}
        \nabla^{\iso}_{\dot{\geodesic}_{\Vector, \VectorB}(t)} \log_{(\cdot)}^{\iso}(\VectorF) = - \dot{\geodesic}_{\Vector, \VectorB}(t),
        \label{eq:prop-iso-mono-bary-field-claim}
    \end{equation}
    because then
    \begin{equation}
        (\partransport_{\VectorB\leftarrow \geodesic_{\Vector, \VectorB}(t)} \nabla^{\iso}_{\dot{\geodesic}_{\Vector, \VectorB}(t)} - \log_{(\cdot)}^{\iso}(\VectorF), \partransport_{\VectorB\leftarrow \Vector} \log_{\Vector} (\VectorB))_{2} = \|\log_{\VectorB} (\Vector)\|_2^2,
    \end{equation}
    and
    \begin{equation}
        \|\partransport_{\VectorB\leftarrow \geodesic_{\Vector, \VectorB}(t)} \nabla^{\iso}_{\dot{\geodesic}_{\Vector, \VectorB}(t)} - \log_{(\cdot)}^{\iso}(\VectorF)\|_2 = \|\log_{\VectorB} (\Vector)\|_2,
    \end{equation}
    which proves the claim.

    There are three possible cases: (1) $\Vector$ is in between $\VectorB$ and $\VectorF$, (2)  $\VectorB$ is in between $\Vector$ and $\VectorF$, and (3) $\VectorF$ is in between  $\Vector$ and $\VectorB$. We will show the claim (\ref{eq:prop-iso-mono-bary-field-claim}) for case (1), which follows in several steps. Cases (2) and (3) follow analogously.

    First, we will evaluate the terms
    \begin{multline}
        \nabla^{\iso}_{\dot{\geodesic}_{\Vector, \VectorB}(t)} \log_{(\cdot)}^{\iso}(\VectorF)\\
         = \frac{1}{\|\dot{\geodesic}_{\Vector, \VectorB}(t)\|_2} \nabla_{\dot{\geodesic}_{\Vector, \VectorB}(t)} \|\partransport_{(\cdot) \leftarrow \geodesic_{\Vector, \VectorB}(t)} \dot{\geodesic}_{\Vector, \VectorB}(t)\|_2 \frac{\int_{0}^{1}\|\dot{\geodesic}_{(\cdot), \VectorF} (s)\|_2 \; \mathrm{d}s}{\|\log_{(\cdot)} (\VectorF)\|_2} \log_{(\cdot)}(\VectorF)\\
        = \frac{1}{\|\dot{\geodesic}_{\Vector, \VectorB}(t)\|_2} \frac{\mathrm{d}}{\mathrm{d}t} \Bigl( \|\dot{\geodesic}_{\Vector, \VectorB}(t)\|_2 \frac{\int_{0}^{1}\|\dot{\geodesic}_{\geodesic_{\Vector, \VectorB}(t), \VectorF} (s)\|_2 \; \mathrm{d}s}{\|\log_{\geodesic_{\Vector, \VectorB}(t)} (\VectorF)\|_2}\Bigr) \log_{\geodesic_{\Vector, \VectorB}(t)}(\VectorF) \\
         \qquad \qquad \qquad \qquad \qquad \qquad \qquad+ \frac{\int_{0}^{1}\|\dot{\geodesic}_{\geodesic_{\Vector, \VectorB}(t), \VectorF} (s)\|_2 \; \mathrm{d}s}{\|\log_{\geodesic_{\Vector, \VectorB}(t)} (\VectorF)\|_2} \nabla_{\dot{\geodesic}_{\Vector, \VectorB}(t)} \log_{(\cdot)}(\VectorF) \\
        = \frac{1}{\|\dot{\geodesic}_{\Vector, \VectorB}(t)\|_2} \frac{\mathrm{d}}{\mathrm{d}t} \Bigl( \|\dot{\geodesic}_{\Vector, \VectorB}(t)\|_2 \frac{\int_{0}^{1}\|\dot{\geodesic}_{\geodesic_{\Vector, \VectorB}(t), \VectorF} (s)\|_2 \; \mathrm{d}s}{\|\log_{\geodesic_{\Vector, \VectorB}(t)} (\VectorF)\|_2}\Bigr) \log_{\geodesic_{\Vector, \VectorB}(t)}(\VectorF)\\
        - \frac{\int_{0}^{1}\|\dot{\geodesic}_{\geodesic_{\Vector, \VectorB}(t), \VectorF} (s)\|_2 \; \mathrm{d}s}{\|\log_{\geodesic_{\Vector, \VectorB}(t)} (\VectorF)\|_2} \dot{\geodesic}_{\Vector, \VectorB}(t).
        \label{eq:prop-iso-bary-field}
    \end{multline}
    
    To further simplify the above expression we will simplify 
    \begin{equation}
        \frac{\int_{0}^{1}\|\dot{\geodesic}_{\geodesic_{\Vector, \VectorB}(t), \VectorF} (s)\|_2 \; \mathrm{d}s}{\|\log_{\geodesic_{\Vector, \VectorB}(t)} (\VectorF)\|_2}.
        \label{eq:iso-log-prefactor}
    \end{equation}
    To do that we will use that 
    \begin{equation}
        \geodesic_{\geodesic_{\Vector, \VectorB}(t), \VectorF} (s) = \geodesic_{\VectorF, \VectorB} \Bigl((1-s)\Bigl(\frac{\distance_{\manifold}(\VectorF, \Vector)}{\distance_{\manifold}(\VectorF, \VectorB)} + \frac{\distance_{\manifold}(\Vector, \VectorB)}{\distance_{\manifold}(\VectorF, \VectorB)}t\Bigr)\Bigr) = \geodesic_{\VectorF, \VectorB} ((1-s)(c + (1-c)t)),
    \end{equation}
    for $c := \frac{\distance_{\manifold}(\VectorF, \Vector)}{\distance_{\manifold}(\VectorF, \VectorB)}$, from which follows that
    \begin{equation}
        \dot{\geodesic}_{\geodesic_{\Vector, \VectorB}(t), \VectorF} (s) = - (c + (1-c)t) \dot{\geodesic}_{\VectorF, \VectorB} ((1-s)(c + (1-c)t)).
        \label{eq:prop-iso-bary-field-geo-rewrite}
    \end{equation}
    Substituting (\ref{eq:prop-iso-bary-field-geo-rewrite}) into (\ref{eq:iso-log-prefactor}) gives
    \begin{multline}
        \frac{\int_{0}^{1}\|\dot{\geodesic}_{\geodesic_{\Vector, \VectorB}(t), \VectorF} (s)\|_2 \; \mathrm{d}s}{\|\log_{\geodesic_{\Vector, \VectorB}(t)} (\VectorF)\|_2} = \frac{\int_{0}^{1}\|\dot{\geodesic}_{\VectorF, \VectorB} ((1-s)(c + (1-c)t))\|_2 \; \mathrm{d}s}{\|\dot{\geodesic}_{\VectorF, \VectorB} (c + (1-c)t)\|_2} \\
        =  \frac{1}{\|\dot{\geodesic}_{\VectorF, \VectorB} (c + (1-c)t)\|_2} \frac{\int_{0}^{c + (1-c)t}\|\dot{\geodesic}_{\VectorF, \VectorB} (s)\|_2 \; \mathrm{d}s}{c + (1-c)t}\\
        = \frac{\distance_{\manifold}(\Vector, \VectorB)}{\distance_{\manifold}(\VectorF, \VectorB)} \frac{1}{\|\dot{\geodesic}_{\Vector, \VectorB} (t)\|_2} \frac{\int_{0}^{c + (1-c)t}\|\dot{\geodesic}_{\VectorF, \VectorB} (s)\|_2 \; \mathrm{d}s}{c + (1-c)t} \\
        = \frac{1-c}{\|\dot{\geodesic}_{\Vector, \VectorB} (t)\|_2} \frac{\int_{0}^{c + (1-c)t}\|\dot{\geodesic}_{\VectorF, \VectorB} (s)\|_2 \; \mathrm{d}s}{c + (1-c)t}.
        \label{eq:prop-iso-bary-field-ratio-rewrite}
    \end{multline}
    
    Next, from (\ref{eq:prop-iso-bary-field-ratio-rewrite}) we can evaluate the first term in (\ref{eq:prop-iso-bary-field})
    \begin{multline}
        \frac{1}{\|\dot{\geodesic}_{\Vector, \VectorB}(t)\|_2} \frac{\mathrm{d}}{\mathrm{d}t} \Bigl( \|\dot{\geodesic}_{\Vector, \VectorB}(t)\|_2 \frac{\int_{0}^{1}\|\dot{\geodesic}_{\geodesic_{\Vector, \VectorB}(t), \VectorF} (s)\|_2 \; \mathrm{d}s}{\|\log_{\geodesic_{\Vector, \VectorB}(t)} (\VectorF)\|_2}\Bigr) \log_{\geodesic_{\Vector, \VectorB}(t)}(\VectorF)  \\
        = \frac{1-c}{\|\dot{\geodesic}_{\Vector, \VectorB}(t)\|_2}  \frac{\mathrm{d}}{\mathrm{d}t} \Bigl( \frac{\int_{0}^{c + (1-c)t}\|\dot{\geodesic}_{\VectorF, \VectorB} (s)\|_2 \; \mathrm{d}s}{c + (1-c)t}\Bigr) \log_{\geodesic_{\Vector, \VectorB}(t)}(\VectorF) \\
        = \frac{1-c}{\|\dot{\geodesic}_{\Vector, \VectorB}(t)\|_2}  \Bigl( \frac{\|\dot{\geodesic}_{\Vector, \VectorB} (t)\|_2}{c + (1-c)t} - \frac{(1-c)}{c + (1-c)t}\frac{\int_{0}^{c + (1-c)t}\|\dot{\geodesic}_{\VectorF, \VectorB} (s)\|_2 \; \mathrm{d}s}{c + (1-c)t} \Bigr) \log_{\geodesic_{\Vector, \VectorB}(t)}(\VectorF) \\
        = - \frac{1-c}{\|\dot{\geodesic}_{\Vector, \VectorB}(t)\|_2}  \Bigl( \frac{\|\dot{\geodesic}_{\Vector, \VectorB} (t)\|_2}{c + (1-c)t} - \frac{(1-c)}{c + (1-c)t}\frac{\int_{0}^{c + (1-c)t}\|\dot{\geodesic}_{\VectorF, \VectorB} (s)\|_2 \; \mathrm{d}s}{c + (1-c)t} \Bigr) \\
        \hfill \times (c + (1-c)t) \dot{\geodesic}_{\VectorF, \VectorB} (c + (1-c)t)\\
        = - \frac{1}{\|\dot{\geodesic}_{\Vector, \VectorB}(t)\|_2}  \Bigl(\|\dot{\geodesic}_{\Vector, \VectorB} (t)\|_2 - (1-c)\frac{\int_{0}^{c + (1-c)t}\|\dot{\geodesic}_{\VectorF, \VectorB} (s)\|_2 \; \mathrm{d}s}{c + (1-c)t} \Bigr)  \dot{\geodesic}_{\Vector, \VectorB} (t)\\
        = - \dot{\geodesic}_{\Vector, \VectorB} (t) + \frac{\int_{0}^{1}\|\dot{\geodesic}_{\geodesic_{\Vector, \VectorB}(t), \VectorF} (s)\|_2 \; \mathrm{d}s}{\|\log_{\geodesic_{\Vector, \VectorB}(t)} (\VectorF)\|_2} \dot{\geodesic}_{\Vector, \VectorB}(t).
        \label{eq:prop-iso-bary-field-ratio-rewrite-term-1}
    \end{multline}
    
    Finally, substituting (\ref{eq:prop-iso-bary-field-ratio-rewrite-term-1}) into (\ref{eq:prop-iso-bary-field}) we find that
    \begin{multline}
        \nabla^{\iso}_{\dot{\geodesic}_{\Vector, \VectorB}(t)} \log_{(\cdot)}^{\iso}(\VectorF) \\
        = - \dot{\geodesic}_{\Vector, \VectorB} (t) + \frac{\int_{0}^{1}\|\dot{\geodesic}_{\geodesic_{\Vector, \VectorB}(t), \VectorF} (s)\|_2 \; \mathrm{d}s}{\|\log_{\geodesic_{\Vector, \VectorB}(t)} (\VectorF)\|_2} \dot{\geodesic}_{\Vector, \VectorB}(t) - \frac{\int_{0}^{1}\|\dot{\geodesic}_{\geodesic_{\Vector, \VectorB}(t), \VectorF} (s)\|_2 \; \mathrm{d}s}{\|\log_{\geodesic_{\Vector, \VectorB}(t)} (\VectorF)\|_2} \dot{\geodesic}_{\Vector, \VectorB}(t)\\
        = - \dot{\geodesic}_{\Vector, \VectorB} (t).
    \end{multline}
    So (\ref{eq:prop-iso-mono-bary-field-claim}) holds, which proves the claim.
\end{proof}

\paragraph{Proofs of the results in \Cref{sec:finding-mins-monotone}}

\begin{proof}[Proof of \Cref{thm:convergence-EGiRD}]
    We will show that for every $\sumIndC$ we have that
    \begin{equation}
        \|\tangentVector_{\Vector^{(\sumIndC+1)}}\|_2 \leq q \|\tangentVector_{\Vector^{(\sumIndC)}}\|_2, \quad 0<q<1,
        \label{eq:thm-convergence-claim}
    \end{equation}
    because then we have that $\lim_{\sumIndC\to\infty} \|\tangentVector_{\Vector^{(\sumIndC)}}\|_2 = 0$ and 
    \begin{multline}
        \|\Vector^{(\sumIndC)} - \Vector^{(\sumIndC+1)} \|_2 = \inf_{\geodesic_{\Vector^{(\sumIndC)}, \Vector^{(\sumIndC+1)}}} \int_{0}^1 \|\dot{\geodesic}_{\Vector^{(\sumIndC)}, \Vector^{(\sumIndC+1)}} (t)\|_2 \; \mathrm{d}t 
        \leq \distance^{\iso}_{\manifold} (\Vector^{(\sumIndC)}, \Vector^{(\sumIndC+1)})  \\
        = \distance^{\iso}_{\manifold} (\Vector^{(\sumIndC)}, \exp^{\iso}_{\Vector^{(\sumIndC)}}(-r \tangentVector_{\Vector^{(\sumIndC)}})) = \|r\tangentVector_{\Vector^{(\sumIndC)}}\|_2
        \leq r q^k\|\tangentVector_{\Vector^{(0)}}\|_2,
    \end{multline}
    which tells us that the sequence $\{\Vector^{(\sumIndC)}\}$ converges to some point $\bar{\Vector}'$ satisfying $\tangentVector_{\bar{\Vector}'}=\mathbf{0}$, because $\lim_{\sumIndC\to\infty} \tangentVector_{\Vector^{(\sumIndC)}}=\mathbf{0}$ and $\tangentVector$ is iso-Lipschitz continuous. In addition, the convergence is linear as
    \begin{equation}
        \|\Vector^{(\sumIndC)} - \bar{\Vector}'\|_2 \leq 
        \sum_{\sumIndC' = \sumIndC}^\infty \|\Vector^{(\sumIndC')} - \Vector^{(\sumIndC'+1)} \|_2 \leq \sum_{\sumIndC' = \sumIndC}^\infty r q^{k'}\|\tangentVector_{\Vector^{(0)}}\|_2 = r \frac{q^k}{1-q}\|\tangentVector_{\Vector^{(0)}}\|_2,
    \end{equation}
    and by 
    \begin{multline}
       0 =  (\tangentVector_{\bar{\Vector}'}, \partransport^{\iso}_{\bar{\Vector}' \leftarrow \bar{\Vector}} \log_{\bar{\Vector}} (\bar{\Vector}'))_2 
        \\
        = ( \tangentVector_{\bar{\Vector}'} - \partransport^{\iso}_{\bar{\Vector}' \leftarrow \bar{\Vector}} \tangentVector_{\bar{\Vector}},  \partransport^{\iso}_{\bar{\Vector}' \leftarrow \bar{\Vector}} \log_{\bar{\Vector}} (\bar{\Vector}') )_2
        \geq \alpha \|\log_{\bar{\Vector}}(\bar{\Vector}')\|_2 \distance^{\iso}_{\manifold} (\bar{\Vector}, \bar{\Vector}') \geq 0,
    \end{multline}
    we must have that $\bar{\Vector}' = \bar{\Vector}$, which proves the claim.

    To show the claim (\ref{eq:thm-convergence-claim}) we first rewrite
    \begin{multline}
        \|\tangentVector_{\Vector^{(\sumIndC+1)}}\|_2^2 = \|\tangentVector_{\Vector^{(\sumIndC+1)}} - \partransport^{\iso}_{\Vector^{(\sumIndC+1)}\leftarrow \Vector^{(\sumIndC)}} \tangentVector_{\Vector^{(\sumIndC)}} + \partransport^{\iso}_{\Vector^{(\sumIndC+1)}\leftarrow \Vector^{(\sumIndC)}} \tangentVector_{\Vector^{(\sumIndC)}} \|_2^2 \\
        = \|\tangentVector_{\Vector^{(\sumIndC+1)}} - \partransport^{\iso}_{\Vector^{(\sumIndC+1)}\leftarrow \Vector^{(\sumIndC)}} \tangentVector_{\Vector^{(\sumIndC)}} \|_2^2 + \| \partransport^{\iso}_{\Vector^{(\sumIndC+1)}\leftarrow \Vector^{(\sumIndC)}} \tangentVector_{\Vector^{(\sumIndC)}} \|_2^2 \\
        + 2 (\tangentVector_{\Vector^{(\sumIndC+1)}} - \partransport^{\iso}_{\Vector^{(\sumIndC+1)}\leftarrow \Vector^{(\sumIndC)}} \tangentVector_{\Vector^{(\sumIndC)}} , \partransport^{\iso}_{\Vector^{(\sumIndC+1)}\leftarrow \Vector^{(\sumIndC)}} \tangentVector_{\Vector^{(\sumIndC)}} )_2
        \label{eq:thm-convergence-claim-rewrite}
    \end{multline}
    and bound the three terms in (\ref{eq:thm-convergence-claim-rewrite}).

    First note that for every $\sumIndC \in \Natural$
    \begin{equation}
        (\tangentVector_{\Vector^{(\sumIndC+1)}} - \partransport^{\iso}_{\Vector^{(\sumIndC+1)}\leftarrow \Vector^{(\sumIndC)}} \tangentVector_{\Vector^{(\sumIndC)}}, \partransport^{\iso}_{\Vector^{(\sumIndC+1)}\leftarrow \Vector^{(\sumIndC)}} \log^{\iso}_{\Vector^{(\sumIndC)}} (\Vector^{(\sumIndC+1)}))_2 \geq \alpha  \distance^{\iso}_{\manifold} (\Vector^{(\sumIndC)}, \Vector^{(\sumIndC+1)})^2,
    \end{equation}
    by $\alpha$-strong iso-monotonicity, which boils down to
    \begin{multline}
        (\tangentVector_{\Vector^{(\sumIndC+1)}} - \partransport^{\iso}_{\Vector^{(\sumIndC+1)}\leftarrow \Vector^{(\sumIndC)}} \tangentVector_{\Vector^{(\sumIndC)}}, \partransport^{\iso}_{\Vector^{(\sumIndC+1)}\leftarrow \Vector^{(\sumIndC)}} - r \tangentVector_{\Vector^{(\sumIndC)}} )_2 \geq \alpha \|r \tangentVector_{\Vector^{(\sumIndC)}}\|_2^2\\
        \Leftrightarrow \quad (\tangentVector_{\Vector^{(\sumIndC+1)}} - \partransport^{\iso}_{\Vector^{(\sumIndC+1)}\leftarrow \Vector^{(\sumIndC)}} \tangentVector_{\Vector^{(\sumIndC)}}, \partransport^{\iso}_{\Vector^{(\sumIndC+1)}\leftarrow \Vector^{(\sumIndC)}} \tangentVector_{\Vector^{(\sumIndC)}} )_2 \leq -\alpha r \| \tangentVector_{\Vector^{(\sumIndC)}}\|_2^2
        \label{eq:thm-convergence-iso-monotone}
    \end{multline}
    by substituting $\log^{\iso}_{\Vector^{(\sumIndC)}} (\Vector^{(\sumIndC+1)}) = \log^{\iso}_{\Vector^{(\sumIndC)}} (\exp^{\iso}_{\Vector^{(\sumIndC)}}(-r \tangentVector_{\Vector^{(\sumIndC)}})) 
    =- r \tangentVector_{\Vector^{(\sumIndC)}}$, and similarly note that by iso-$L$-Lipschitz continuity
    \begin{equation}
        \| \tangentVector_{\Vector^{(\sumIndC+1)}} - \partransport^{\iso}_{\Vector^{(\sumIndC+1)}\leftarrow \Vector^{(\sumIndC)}}\tangentVector_{\Vector} \|_2 \leq L \distance^{\iso}_{\manifold}(\Vector^{(\sumIndC)}, \Vector^{(\sumIndC+1)}) = L r \|\tangentVector_{\Vector^{(\sumIndC)}} \|_2.
        \label{eq:thm-convergence-iso-Lipschitz}
    \end{equation}

    In other words, (\ref{eq:thm-convergence-iso-monotone}) and (\ref{eq:thm-convergence-iso-Lipschitz}) give us upper bounds for the first and third term in (\ref{eq:thm-convergence-claim-rewrite}). We can also rewrite the second term in (\ref{eq:thm-convergence-claim-rewrite}) as
    \begin{equation}
        \| \partransport^{\iso}_{\Vector^{(\sumIndC+1)}\leftarrow \Vector^{(\sumIndC)}} \tangentVector_{\Vector^{(\sumIndC)}} \|_2^2 = \| \frac{\|- r \tangentVector_{\Vector^{(\sumIndC)}}\|_2}{\|\partransport_{\Vector^{(\sumIndC+1)}\leftarrow \Vector^{(\sumIndC)}} -r \tangentVector_{\Vector^{(\sumIndC)}}\|_2} \partransport_{\Vector^{(\sumIndC+1)}\leftarrow \Vector^{(\sumIndC)}} \tangentVector_{\Vector^{(\sumIndC)}} \|_2^2 = \| \tangentVector_{\Vector^{(\sumIndC)}} \|_2^2.
        \label{eq:thm-convergence-claim-2}
    \end{equation}

    Substituting (\ref{eq:thm-convergence-iso-monotone}), (\ref{eq:thm-convergence-iso-Lipschitz}) and (\ref{eq:thm-convergence-claim-2}) into (\ref{eq:thm-convergence-claim-rewrite}) gives us that
    \begin{equation}
        \|\tangentVector_{\Vector^{(\sumIndC+1)}}\|_2^2 \leq L^2 r^2 \| \tangentVector_{\Vector^{(\sumIndC)}} \|_2^2 + \| \tangentVector_{\Vector^{(\sumIndC)}} \|_2^2 - 2 \alpha r \| \tangentVector_{\Vector^{(\sumIndC)}} \|_2^2 = (1 + r(L^2r - 2 \alpha )) \| \tangentVector_{\Vector^{(\sumIndC)}} \|_2^2.
    \end{equation}
    So we have
    \begin{equation}
        \|\tangentVector_{\Vector^{(\sumIndC+1)}}\|_2 \leq q \|\tangentVector_{\Vector^{(\sumIndC)}}\|_2,
    \end{equation}
    where
    \begin{equation}
        q:= (1 + r(L^2 r - 2 \alpha))^{\frac{1}{2}} \in (0,1),
    \end{equation}
    because $0 < r < \frac{2 \alpha}{L^2}$, which proves the claim (\ref{eq:thm-convergence-claim}).
\end{proof}


\paragraph{Proofs of the results in \Cref{sec:iso-bary}}

\begin{proof}[Proof of \Cref{prop:two-point-iso-bary}]
    We need to show that $\geodesic^{\iso}_{\Vector^1, \Vector^2} (\frac{1}{2})$ solves for $\bar{\Vector}$ in the equation
    \begin{multline}
        \frac{1}{2} \log_{\bar{\Vector}}^{\iso}(\Vector^1) + \frac{1}{2} \log_{\bar{\Vector}}^{\iso}(\Vector^2) = \mathbf{0} \\
        \Leftrightarrow \quad \frac{1}{2} \frac{\int_{0}^{1}\|\dot{\geodesic}_{\bar{\Vector}, \Vector^1} (s)\|_2 \; \mathrm{d}s}{\|\log_{\bar{\Vector}} (\Vector^1)\|_2}\log_{\bar{\Vector}}(\Vector^1) + \frac{1}{2} \frac{\int_{0}^{1}\|\dot{\geodesic}_{\bar{\Vector}, \Vector^2} (s)\|_2 \; \mathrm{d}s}{\|\log_{\bar{\Vector}} (\Vector^2)\|_2} \log_{\bar{\Vector}}(\Vector^2) = \mathbf{0},
    \end{multline}
    which is also equivalent to
    \begin{equation}
        \frac{\int_{0}^{1}\|\dot{\geodesic}_{\bar{\Vector}, \Vector^1} (s)\|_2 \; \mathrm{d}s}{\|\log_{\bar{\Vector}} (\Vector^1)\|_2}\log_{\bar{\Vector}}(\Vector^1) = - \frac{\int_{0}^{1}\|\dot{\geodesic}_{\bar{\Vector}, \Vector^2} (s)\|_2 \; \mathrm{d}s}{\|\log_{\bar{\Vector}} (\Vector^2)\|_2} \log_{\bar{\Vector}}(\Vector^2).
        \label{eq:prop-two-point-iso-bary-cond}
    \end{equation}
    
    We first note that the expression (\ref{eq:prop-two-point-iso-bary-cond}) tells us that we must have that the iso-barycentre $\bar{\Vector}$ lives on the geodesic $\geodesic_{\Vector^1, \Vector^2}$. So it remains to check where on the geodesic $\bar{\Vector}$ lives. 

    To do that, we see by taking $\ell^2$-norms on both sides of (\ref{eq:prop-two-point-iso-bary-cond}) that
    \begin{equation}
        \int_{0}^{1}\|\dot{\geodesic}_{\bar{\Vector}, \Vector^1} (s)\|_2 \; \mathrm{d}s = \int_{0}^{1}\|\dot{\geodesic}_{\bar{\Vector}, \Vector^2} (s)\|_2 \; \mathrm{d}s
    \end{equation}
    must hold, which means that $\bar{\Vector}$ lies at the point on the geodesic that has passed the same amount of arc length (in the $\ell^2$ sense) after leaving $\Vector^1$ as it still needs to pass to reach $\Vector^2$. Since iso-geodesics have constant $\ell^2$-speed, we must have that $\bar{\Vector}$ lies at the iso-geodesic midpoint, i.e.,  $\bar{\Vector} = \geodesic^{\iso}_{\Vector^1, \Vector^2} (\frac{1}{2})$, which proves the claim.
\end{proof}

\begin{proof}[Proof of \Cref{prop:iso-bary-R}]
    First, we note that on $\Real$ the iso-logarithmic mapping reduce to simpler expressions, i.e., for $x,y \in \Real$ we have
    \begin{equation}
        \log^{\iso}_{x} (y) = y - x.
    \end{equation}
    From the above expression it follows directly that
    \begin{equation}
        \frac{1}{\dataPointNum} \sum_{\sumIndA=1}^\dataPointNum \log^{\iso}_{\bar{x}} (x^\sumIndA) = 0 \quad \Leftrightarrow \quad \bar{x} = \frac{1}{\dataPointNum}\sum_{\sumIndA=1}^\dataPointNum x^\sumIndA,
    \end{equation}
    which proves the claim.
\end{proof}

\begin{proof}[Proof of \Cref{thm:existence-iso-barycentre}]
    By assumption, there exists a geodesically convex submanifold
    \(\manifold' \subset \manifold\) containing the data points and $\manifold'$ is even strongly geodesically convex by assumptions on the iso-exponential map. 
    
    Let \(\Vector \in \partial \manifold'\). Since
    \(\Vector^\sumIndA \in \manifold'\) for each
    \(\sumIndA = 1,\ldots,\dataPointNum\), strong geodesic convexity implies that
    \begin{equation}
        \log_{\Vector}^{\iso}(\Vector^\sumIndA)
    \end{equation}
    points into \(\manifold'\). Consequently, the vector field
    \begin{equation}
        \Vector \to
        \frac{1}{\dataPointNum}
        \sum_{\sumIndA=1}^{\dataPointNum}
        \log_{\Vector}^{\iso}(\Vector^\sumIndA)
    \end{equation}
    also points into \(\manifold'\) along \(\partial\manifold'\). Since this
    vector field is smooth on \(\manifold\), the Poincar\'e-Hopf theorem implies
    that it has at least one zero in \(\manifold'\). This proves the claim.
\end{proof}

\begin{proof}[Proof of \Cref{thm:iso-mono-bary-field-1d}]
    The proof follows directly from the proof of \Cref{prop:iso-mono-log-field}. In particular, we have for any $\Vector\neq \VectorB\in \manifold$ and $t\in [0,1]$
    \begin{equation}
        \nabla^{\iso}_{\dot{\geodesic}_{\Vector, \VectorB}(t)} \log_{(\cdot)}^{\iso}(\Vector^\sumIndA) = - \dot{\geodesic}_{\Vector, \VectorB}(t), \quad \sumIndA = 1, \ldots, \dataPointNum,
    \end{equation}
    from which follows that
    \begin{equation}
        (\partransport_{\VectorB\leftarrow \geodesic_{\Vector, \VectorB}(t)} \nabla^{\iso}_{\dot{\geodesic}_{\Vector, \VectorB}(t)} -\frac{1}{\dataPointNum} \sum_{\sumIndA=1}^\dataPointNum 
    \log_{(\cdot)}^{\iso}(\Vector^\sumIndA), \partransport_{\VectorB\leftarrow \Vector} \log_{\Vector} (\VectorB))_{2} = \|\log_{\VectorB} (\Vector)\|_2^2,
    \end{equation}
    and
    \begin{equation}
        \|\partransport_{\VectorB\leftarrow \geodesic_{\Vector, \VectorB}(t)} \nabla^{\iso}_{\dot{\geodesic}_{\Vector, \VectorB}(t)} -\frac{1}{\dataPointNum} \sum_{\sumIndA=1}^\dataPointNum 
    \log_{(\cdot)}^{\iso}(\Vector^\sumIndA)\|_2 = \|\log_{\VectorB} (\Vector)\|_2,
    \end{equation}
    which proves the claim. 
\end{proof}

\subsection{Iso-K-means clustering}
\label{sec:iso-k-means}

To further highlight the practical benefits of iso-barycentres -- which, as demonstrated in \Cref{fig:synthetic-results-barys}, can yield results comparable to the Riemannian barycentre -- we introduce an iso-Riemannian variant of K-means clustering. This is achieved by adapting Lloyd's algorithm, the standard approach for K-means, resulting in \Cref{alg:iso-riem-kmeans}.

\begin{remark}
    Similarly to how the iso-barycentre is defined, iso-K-means is not described as the minimizer of a particular function, but instead is specified by analogy to an alternative characterization of the standard K-means method. It is important to note that because of this, we do not inherit convergence of \Cref{alg:iso-riem-kmeans} to a local minimizer from Lloyd's algorithm. In particular, convergence of \Cref{alg:iso-riem-kmeans} is still open, but beyond the scope of this work.
\end{remark}


\begin{algorithm} 
\caption{Iso-K-means clustering through Lloyd's algorithm with iso-barycentres} 
\label{alg:iso-riem-kmeans} 
\begin{algorithmic} 
\REQUIRE Data points $\{\Vector^{\sumIndA}\}_{\sumIndA=1}^\dataPointNum \subset \manifold$, number of clusters $K$, iteration counter $k = 0$
\STATE $\clusterLabel^{(0)}, \centroid^{1,(0)},\dots,\centroid^{K,(0)} \in \operatorname{argmin}_{\clusterLabel\in \{1,\dots,K\}^\dataPointNum, \centroid^1, \ldots, \centroid^K \in \manifold} \sum_{\sumIndA=1, \sumIndB = 1}^{\dataPointNum, K} \mathds{1}_{\{\clusterLabel_{\sumIndA} = \sumIndB \}} \distance_{\manifold} (\Vector^\sumIndA, \centroid^{\sumIndB})^2$ \hfill (Riemannian K-means)
\WHILE{not converged}
    \STATE $\clusterLabel_\sumIndA^{(k+1)} \gets \arg\min_{j \in \{1,\dots,K\}} d^{\iso}_{\manifold}(\Vector^{\sumIndA},\centroid^{\sumIndB,(k)}) \quad \forall \sumIndA=1,\dots,\dataPointNum$
    \STATE $\centroid^{\sumIndB,(k+1)} \gets \centroid^{\sumIndB}$, where $\centroid^{\sumIndB}$ is the iso-barycentre of $ \{\Vector^\sumIndA \mid \clusterLabel_\sumIndA^{(k+1)} = \sumIndB\} \quad \forall \sumIndB=1,\dots,K$ \hfill (\Cref{alg:iso-bary-line-search})
    \STATE $k \gets k+1$
\ENDWHILE
\end{algorithmic} 
\end{algorithm}

Contrary to the results for barycentres as discussed in \Cref{rem:little-difference}, there is a more significant discrepancy in terms of performance for the downstream task of clustering. That is, \Cref{fig:synthetic-results-k-means} illustrates that only iso-K-means is able to find both correct clusters while maintaining meaningful centroids, which can be attributed to the fact that the underlying pullback metric does not encode our intuition on $\ell^2$-proximity, which deteriorates clustering performance when not corrected for (through iso-Riemannian geometry). We refer the reader to \Cref{sec:app-iso-bary} for details, including the stopping criterion used.

\begin{remark}
\label{rem:roundig-errors}
    Because numerical errors accumulate in the implementation described by \cite{diepeveen2025manifold} and used in this work -- as noted in \Cref{rem:num-error-pile-up} -- \Cref{alg:iso-riem-kmeans} stalls during the barycentre update step on the MNIST data set, preventing us from presenting results here. As also indicated in \Cref{rem:num-error-pile-up}, addressing these numerical difficulties is left for future work.
\end{remark}

\subsection{Details to numerical experiments in \Cref{sec:iso-bary}}
\label{sec:app-iso-bary}

\subsubsection{River and spiral pullback geometries}

\paragraph{The river diffeomorphism}

The river diffeomorphism $\diffeo_{\text{river}} : \Real^2 \to \Real^2$ is defined as
\begin{equation}
    \diffeo_{\text{river}}(\Vector) := \left( \Vector_1 - \beta \sin(\Vector_2),\; \sinh(\eta \Vector_2) \right),
\end{equation}
where $\beta, \eta > 0$ are fixed parameters. The inverse diffeomorphism is
\begin{equation}
    \diffeo_{\text{river}}^{-1}(\VectorB) := \left( \VectorB_1 + \beta \sin\left( \frac{\operatorname{arcsinh}(\VectorB_2)}{\eta} \right),\; \frac{\operatorname{arcsinh}(\VectorB_2)}{\eta} \right).
\end{equation}
In the experiments we used $\beta = 5$ and $\eta = 0.25$.

\paragraph{The spiral diffeomorphism}

The spiral diffeomorphism $\diffeo_{\text{spiral}} : \Real^2 \to \Real^2$ is defined by the forward transformation
\begin{equation}
    \diffeo_{\text{spiral}}(\Vector) := \left( \frac{R}{\beta},\, (\phi - \frac{R}{\beta}) \bmod 2\pi \right),
\end{equation}
where $R := \sqrt{\Vector_1^2 + \Vector_2^2}$ and $\phi := \angle\big((1,0), (\Vector_1, \Vector_2)\big)$, and $\beta> 0$ is a fixed parameter. The inverse transformation is given by
\begin{equation}
    \diffeo_{\text{spiral}}^{-1}(\mathbf{p}) := \beta r \begin{pmatrix} \cos(r+\theta) \\ \sin(r+\theta) \end{pmatrix}
\end{equation}
where $r := \mathbf{p}_1$, $\theta := \mathbf{p}_2$. In the experiments we used $\beta = 0.25$.

\subsubsection{Additional results}

\begin{figure}[h!]
    \centering
    \setlength{\tabcolsep}{2pt}
    \renewcommand{\arraystretch}{0}
    \begin{tabular}{ccc}
        \includegraphics[width=0.48\textwidth]{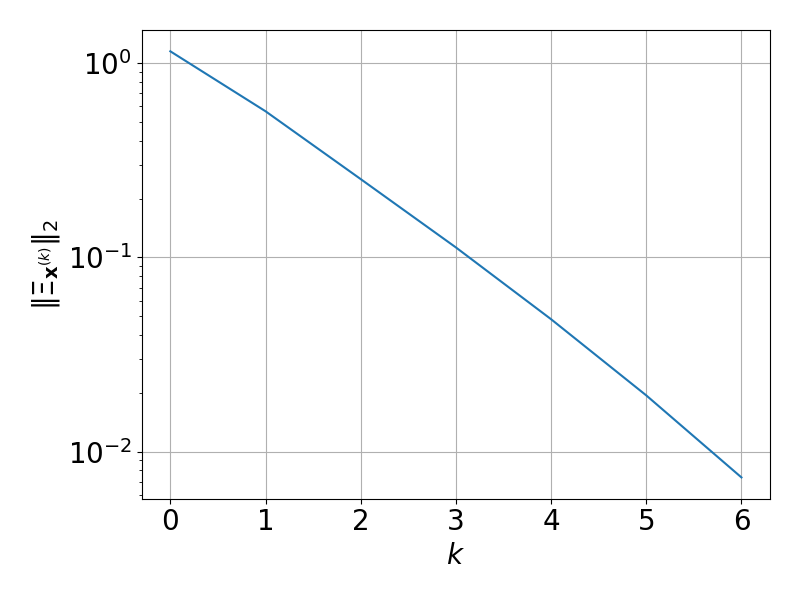} &
        \includegraphics[width=0.48\textwidth]{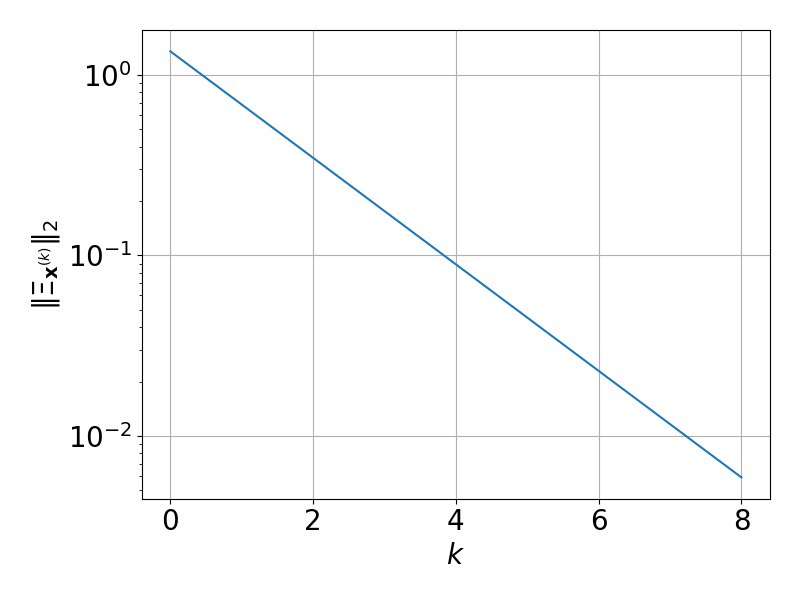}  \\
        \makebox[0.29\textwidth]{\small River data set} &
        \makebox[0.29\textwidth]{\small Spiral data set}
    \end{tabular}
    \caption{Linear convergence to the iso-barycentre of the two synthetic data sets, which is in line with our theory.}
    \label{fig:conv-results-barys}
\end{figure}

\paragraph{Iso-barycentre convergence}
The convergence plots of \Cref{alg:iso-bary-line-search} for the river, spiral and MNIST data sets are shown in \Cref{fig:conv-results-barys} and are in line with our expectations. That is, we attain linear convergence in 
$$
k\mapsto \|\tangentVector_{\Vector^{(k)}}\|_2 := \|\frac{1}{\dataPointNum} \sum_{\sumIndA=1}^\dataPointNum \log^{\iso}_{\Vector^{(k)}} (\Vector^\sumIndA)\|_2.
$$ The river and spiral experiments were run until $\|\tangentVector_{\Vector^{(k)}}\|_2 < 10^{-2}$, whereas for MNIST the algorithm was terminated once $\|\tangentVector_{\Vector^{(k)}}\|_2 < 10^{-1}$.


\begin{figure}[h!]
    \centering
    \setlength{\tabcolsep}{2pt}
    \renewcommand{\arraystretch}{0}
    \begin{tabular}{c}
        \includegraphics[width=0.48\textwidth]{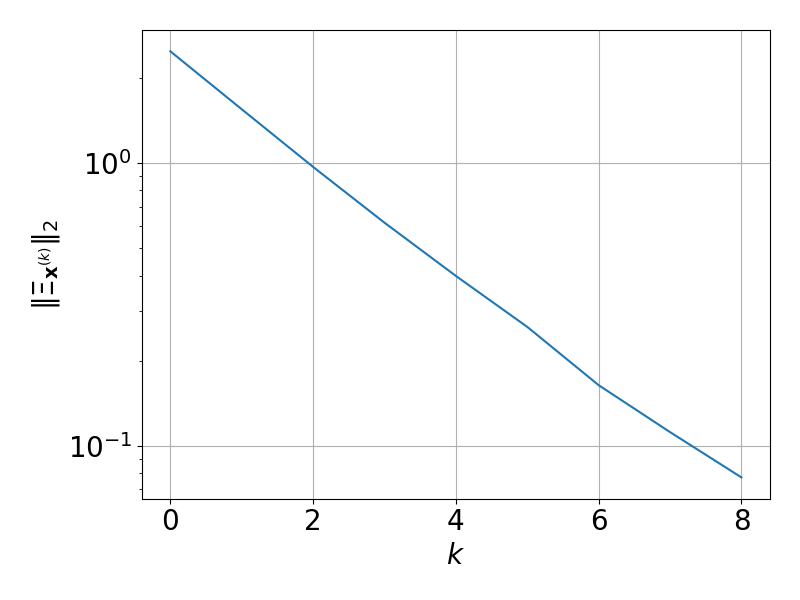} \\
        \makebox[0.29\textwidth]{\small MNIST data set}
    \end{tabular}
    \caption{Linear convergence to the iso-barycentre on the MNIST data set, which is in line with our theory.}
    \label{fig:conv-results-barys-mnist}
\end{figure}

\paragraph{Iso-K-means convergence}
All versions of K-means were run until the norm of the difference between the old and the new centroids satisfies $$\sqrt{\sum_{\sumIndB=1}^K \|\centroid^{\sumIndB,(k+1)} - \centroid^{\sumIndB,(k)}\|_2^2} < 10^{-4}.$$

\end{document}